\documentclass[12pt]{article}

\usepackage{geometry}
\usepackage{amsfonts,amssymb,amsmath,amsthm}
\usepackage{graphicx}
\usepackage{float}
\usepackage[title]{appendix}

\newcommand{\RR}{\mathbb R}
\newcommand{\CC}{\mathbb C}
\newcommand{\ZZ}{\mathbb Z}

\renewcommand{\Im}{\mathop{\mathrm{Im}}}

\newcommand{\sign}{\mathop{\mathrm{sign}}}
\newtheorem{theorem}{Theorem}
\newtheorem{lemma}{Lemma}
\newtheorem{remark}{Remark}
\newtheorem{example}{Example}
\newtheorem{definition}{Definition}
\newtheorem{corollary}{Corollary}
\def \GTNN {{Gr^{\mbox{\tiny TNN}} (k,n)}}
\newcommand\mycom[2]{\genfrac{}{}{0pt}{}{#1}{#2}}

\begin{document}

\title{Parametrization of positroid cells by Dubrovin-Natanzon Cartier divisors on maximal Mumford curves.\thanks{
The work of S. Abenda was supported by HORIZON-MSCA-2022-SE-01-01 CaLIGOLA,  COST Action CaLISTA CA21109, GNFM-INdAM,  INFN projects MMNLP and GAST.
The work of P.G. Grinevich was performed at the Steklov International Mathematical Center and supported by the Ministry of Science and Higher Education of the Russian Federation (agreement no. 075-15-2025-303).}}
\author{Simonetta Abenda\thanks{Dipartimento di Matematica, Universit\`a di Bologna and INFN Sezione di Bologna, Italy; e-mail: simonetta.abenda@unibo.it}, Petr G. Grinevich\thanks{Steklov Mathematical Institute of Russian Academy of Sciences, Moscow, Russia;
L.D. Landau Institute for theoretical Physics, of Russian Academy of Sciences  Russia;
Lomonosov Moscow State University, Russia; grinev@mi-ras.ru}}

\maketitle

\abstract{We prove that each totally positive Schubert cell is a real component of the Jacobian of the maximal Mumford ($\mathtt{MM}$) curve 
whose dual graph is the Le-graph 
representing such cell. At this aim we construct a bijection between the points of this positroid cell and Dubrovin-Natanzon Cartier divisors on the $\mathtt{MM}$-curve.}

\textbf{Keywords:} {Weil divisors and Cartier divisors, Tropical curves and Jacobians, Maximal Mumford curves, Dubrovin-Natanzon divisors, Positroid cells, Degeneration of Spectral curves, Kadomtsev-Petviashvili real regular solitons.}

\maketitle

\tableofcontents

\section{Introduction}

Our interest in the connection between totally non-negative Grassmannians and divisors on rational degenerations of ${\mathtt{M}}$-curves was originally motivated by the theory of real and regular multiline soliton solutions for the Kadomtsev-Petviashvili-II (KP-II) equation. Such solutions can be constructed using two different approaches.

The first one is based on the Darboux transform, and real and regular KP-II multiline solitons correspond to data belonging to the totally non-negative Grassmannians $\GTNN$ \cite{Malanyuk1991,KodamaWilliams2013}. In his seminal paper Postnikov \cite{Postnikov2006} constructed a positroid cells decomposition for 
$\GTNN$. The positroid cells are obtained intersecting Gelfand-Serganova strata \cite{GelfandSerganova1987} with $\GTNN$. Each positroid cell is represented by an equivalence class of planar bicolored perfectly oriented graphs in the disc and has a canonical representative, the so--called Le--graph; positroid cells are parameterized by positive edge weights on these graphs \cite{Postnikov2006}. 
Real and regular KP-II multiline solitons with data in the same positroid cell share the same asymptotic behaviour in space--time, and the combinatorics of positroid cells has been used to characterize such asymptotic properties in \cite{ChakravartyKodama2009,KodamaWilliams2014}.

The second approach is based on finite gap theory of soliton equations where the spectral data are a spectral curve $\Gamma$ and a divisor $\mathcal D$. In this approach multiline soliton  solutions can be constructed as degenerations of the finite-gap ones, and this procedure uses proper degenerations of spectral curves. The idea of constructing soliton solutions from the finite-gap ones by degenerating spectral curves was first used in the theory of the Korteweg-de Vries equation by Krichever \cite {Krichever1975} and Mumford \cite{Mumford1984}. 

In the KP-II finite-gap theory, real regular solutions correspond to Dubrovin-Natanzon divisors on ${\mathtt{M}}$-curves, as proven in \cite{Dubrovin1982,DubrovinNatanzon1988}. In our papers \cite{AG1,AG3,AG4} we proved that all real regular KP-II multiline soliton solutions are obtained degenerating real and regular KP-II finite--gap solutions. Indeed, we proved that such multi solitons can be obtained from spectral problems associated with Dubrovin--Natanzon divisors on so-called maximal Mumford (${\mathtt{MM}}$) curves, whose dual graph is a member in the Postnikov class representing the positroid cell for the given soliton data.

More recently, Agostini et al \cite{AgostiniFevolaMandelshtamSturmfels2023}, Ichikawa \cite{Ichikawa2023} have studied the degenerations of quasiperiodic KP solutions to the soliton ones in the framework of tropical curves. In particular, in \cite{FevolaMandelshtam2024}  and \cite{AbendaCelikFevolaMandelshtam2025} the inverse spectral problem for KP soliton solutions associated with tropical spectral curves with one or two components is solved using the combinatorics of metric graphs in terms of Voronoy and Delaunay polytopes.

In this paper we continue the investigation of the relations between the two aforementioned approaches, and we establish a bijection between each totally positive Schubert cell and the Dubrovin-Natanzon Cartier divisors on the ${\mathtt{MM}}$-curve whose dual graph is the corresponding Le-graph. In this way we show that such positroid cells can be treated as real components of the corresponding tropical Jacobians. 

\bigskip

Let us recall that an $\mathtt M$-curve is a curve with a real involution fixing the maximal possible number of real ovals. This number is equal to $g+1$, where $g$ is the genus. The theory of  $\mathtt M$-curves and the procedure for constructing regular curves from the degenerate ones are discussed in details, in particular, in \cite{Viro1990}.

In \cite{DubrovinNatanzon1988} it was shown that a finite-gap KP-II solution associated to a smooth spectral curve $\Gamma$ and a divisor $\mathcal D$ is real and regular if and only if
\begin{itemize}
\item $\Gamma$ is an  $\mathtt M$-curve,
\item the essential singularity of the wave function belongs to one of these ovals,
\item each other oval contains exactly one point of $\mathcal D$. In the following we call such divisors \textbf{Dubrovin-Natanzon} (DN) divisors 
\end{itemize}

${\mathtt{MM}}$-curves are rational degenerations of ${\mathtt{M}}$-curves in the sense of \cite{KummerSturmfelsVlad2025}. Here ${\mathtt{MM}}$ stands for Maximal Mumford. Mumford curve means that it is ``maximally non smooth and reducible'' in the sense of \cite{Tyurin2003}. A Mumford curve is maximal if, in addition, the number of real ovals is $g+1$ where $g$ is the arithmetic genus.

In our construction \cite{AG3,AG4}, the $\mathtt{MM}$-curve has dual graph in the Postnikov class representing the positroid cell for the given soliton data. In particular,  we associate the ovals to the faces of the dual graph, and for any point in the positroid cell we choose a normalization time and construct a unique DN divisor by solving a system of relations on the corresponding network.

It is well-known that for smooth curves Weil divisors and Cartier divisors are equivalent notions, and not equivalent ones for singular curves. In the aforementioned papers we worked with Weil divisors (collections of points), and avoided the case of divisors located at the intersections of rational components by using the freedom of choice of the normalization KP time for the KP wave function.

If all divisor points are located outside the singular points, we call such divisors \textbf{smooth}. Smooth divisors can be equivalently treated as Weil or Cartier ones. In \cite{AG3,AG4}, for any given normalization time, smooth divisors parametrize an open subset of the positroid cell associated with the soliton data. However, to establish a bijection between the positroid cell of the Grassmannian and the DN divisors on a given spectral curve, one has to fix the normalization time and allow divisors to be located at intersection points. If one treats these divisors as Weil, the regularity of the map between soliton data and divisors is lost.

To pass from Weil to Cartier formulation of divisors, we apply blow-up procedures. We show that, if one treats DN divisors as Cartier divisors, the map between soliton data and divisors is regular. 

In the following we use the term \textbf{Dubrovin-Natanzon Cartier divisors} to stress that DN divisors are treated as Cartier ones.

More precisely, in our text we prove that totally positive Schubert cells are parametrized by DN Cartier divisors on  ${\mathtt{MM}}$-curves whose dual graph is Postnikov's Le-graph by the following steps:
\begin{enumerate}
\item  For DN divisors on ${\mathtt{MM}}$-curves whose dual graphs are trivalent Le-graphs of totally positive Schubert cells $S$, only combinations of two basic types of blow-ups may arise.
\item The Krichever's wave function never vanishes identically in the spectral parameter at any rational component. Therefore the divisor of zeroes of the KP wave function is a well-defined DN Cartier divisor for all KP times.
\item For any fixed normalization time, to any point in $S$ we associate a unique DN Cartier divisor on the spectral ${\mathtt{MM}}$-curve.
\item Vice versa, for any fixed normalization time, to any DN Cartier divisor we associate a unique point in $S$ by reconstructing the KP-II soliton data. 
\end{enumerate}

\bigskip

Let us explain which singular divisor configurations may occur at finite KP times. In the KP dynamics the Abel transform has to remain finite because the KP dynamics is linear in the Jacobian. Therefore the allowed singular configurations are those in which at an intersection point of the curve we have a divisor point on each intersecting component, and all these points simultaneously reach the intersection point in such a way that the Abel transform remains finite. In this case the arithmetic genus of the curve is preserved, the solution remains regular, and the soliton data belong to the same positroid cell also in the limit. For such singular configuration, if divisors are treated as  Weil divisors, i.e. formal sums of points with integer coefficients, the Abel map is not continuous any more, and we apply a blow-up procedure to pass from Weil divisors to Cartier ones.

On the contrary, if only one of the divisor points tends to an intersection point, the Abel transform of the divisor goes to infinity. Then the component containing this divisor point separates, and the arithmetic genus drops. It means that we reach the boundary in the moduli space of curves of given genus. At the level of the soliton data it means that we reach the boundary of the positroid cell, and the limiting point belong to a lower-dimensional cell. For equations like Kadomtsev-Petviashvili such situation never happens at finite KP times. This behavior is not shared by all soliton systems. For instance, for Melnikov-type equations the arithmetic genus may drop at finite times, see \cite{GrinevichTaimanov2007,GrinevichTaimanov2008}. 

Let us remark that the perfect orientations of any Postnikov graph are totally cyclic orientations in the sense of \cite{Alexeev2004}, and that the positroid cell parametrizes an open component of the real tropical Jacobian. Moving to the boundary of the positroid cell is equivalent to moving some divisor points to the nodal point in a way not compatible with the KP dynamics. It would be relevant to study the compactification of Jacobians for singular curves in terms of KP divisors. The problem of the  compactification of tropical Jacobians
was discussed by several authors, in particular we would like to cite \cite{Caporaso1994}, \cite{Esteves2001}, \cite{AbreuPacini2020}. In this respect it would be useful to implement the isomorphism between equivalence classes of totally cyclic orientations on subgraphs of metric graphs and equivalence classes of faces of the Voronoi polytopes found in \cite{Amini2011},  in the case when the metric graph is a Le-graph for some orientation.

These questions are naturally connected to the tropical version of the Schottky problem, see for instance \cite{CaporasoViviani2010} and references therein, and integrable systems, see the review \cite{GrushevskyXie2025} for more details. 

Singular Riemann surfaces naturally arise in soliton theory, for a review see \cite{Taimanov2011}. In particular, an interesting application of degenerate curves is associated with the problem of constructing classical orthogonal coordinate systems in the framework of inverse scattering. Applicability of inverse scattering to this geometrical problem was shown in \cite{Zakharov1998}, and an algebro-geometric version of this construction was suggested in \cite{Krichever1997}. The solution of this problem was obtained by Mironov and Taimanov in \cite{MironovTaimanov2006,MironovTaimanov2007}, for further applications see \cite{MironovSenningerTaimanov2023}, \cite{GlukhovMokhov2020}. The spectral curves used in these papers are unions of Riemann spheres with double points.

The theory of divisors on singular curves was started by Knudsen and Mumford in \cite{KnudsenMumford1976,Knudsen1983,Knudsen1983_2}, the divisors and Jacobians of such curves were also discussed by Artamkin \cite{Artamkin2004,Artamkin2004_2,Artamkin2006}. These papers use the Cartier definition of divisors, i.e. the language of line bundles. 

On the contrary, in the majority of soliton papers dealing with singular spectral curves, divisors were treated as collections of points, i.e. as Weil divisors, and the cases where the divisor points are located at singular points were not discussed. Cartier divisors were used in \cite{Parshin2001,KurkeOsipovZheglov2014} for the study of multidimensional generalizations of the Krichever construction of commuting operators, where the spectral varieties are multidimensional and  always singular. 
Positroid cells appear also in connection to dimer models in the disc \cite{PostnikovSpeyerWilliams2009}, \cite{Lam2015}. Such statistical model is naturally linked to the spectral problem for KP-II real and regular multiline solitons since the system of relations used to compute the KP wave function at the nodes is governed by a Kasteleyn matrix \cite{Abenda2021}. More recently, in \cite{Meszarosetal2026} a relevant connection between non-vanishing Pl\"ucker coordinates and cluster monomials was discovered, and it will be interesting to understand if our parametrization of totally positive Schubert cells can be reformulated in terms of cluster algebras \cite{FominZelevinsky2002}. 

Finally, we would like to  mention another open problem connected with tropical geometry and  KP multiline real regular solitons, which we postpone for a future publication. Indeed, in the study of the dynamics of these solitons in $x,y$-space both for fixed $t$ and in the large $t$ limit, it was pointed out in \cite{BiondiniChakravarty2006,KodamaWilliams2014}  that the maxima of these solitons are tropical curves in the sense of \cite{MikhalkinZharkov2008}, see also the book \cite{MikhalkinRau2018}. Since $\mathtt{MM}$-curves are tropical curves in the sense of \cite{KummerSturmfelsVlad2025}, it would be interesting to understand the relation between the tropical geometry associated with the KP dynamics and the tropical geometry of the spectral $\mathtt{MM}$-curves.

\section{Kadomtsev-Petviashvili-II equation, totally non-negative Grassmannians and ${\mathtt{MM}}$-curves}

The Kadomtsev-Petviashvili (KP) equation \cite{KP}
\begin{equation}\label{eq:KPs}
(-4u_t+6uu_x+u_{xxx})_x+ \alpha^2  u_{yy}=0,
\end{equation}
was originally derived as a universal model for wave propagation in weakly dispersive weakly non-linear systems with 2 spatial variables under the additional assumption that the solution slowly depends on the $y$-variable. 

KP equation has two real forms. If $\alpha^2<0$, i.e. $\alpha$ is purely imaginary, the corresponding equation is known as KP-I; if  $\alpha^2>0$, i.e. $\alpha$ is real, the  corresponding equation is known as KP-II. In particular, equation (\ref{eq:KPs}) is used as model for water surface waves \cite{Osborne2010}. If the depth of the water is less than about a centimeter and the surface tension prevails, the wave motion is guided by KP-I; if the water is shallow, but much deeper than a centimeter, KP-II can be used (see \cite{Osborne2010}).
The analytic properties of the solutions of these two KP real forms essentially differ, and the study of these solutions require different mathematical tools.

In our paper we concentrate on the study of the Kadomtsev-Petviashvili-II equation, and for the sake of shortness we will call it KP. Without loss of generality we assume that $\alpha^2=3$, therefore KP has the form
\begin{equation}\label{eq:KP}
(-4u_t+6uu_x+u_{xxx})_x+ 3  u_{yy}=0.
\end{equation}

Formally the KP equation was derived under the assumption that the $y$-dependence of the wave is much slower than the $x$-dependence, and this condition is not true for real surface water waves. Nevertheless the shape of real water waves is well-approximated by multisoliton KP solutions \cite{AblowitzBaldwin2012}.

KP equation is very interesting from the mathematical point of view, including its applications to Riemann-Schottky problem and the theory of the moduli spaces of algebraic curves.

One of the KP key features is that it can be integrated using the inverse scattering transform (IST). The zero-curvature representation for KP was found by Druma \cite{Druma1973}, Zakharov and Shabat  \cite{ZakharovShabat1974}
\begin{equation}\label{eq:KP-Lax}
(\partial_y -B_2)\Psi(P,x,y,t) =0, \qquad (\partial_t -B_3)\Psi(P,x,y,t) =0,
\end{equation}
$$
B_2 = \partial_x^2 + u, \ \ B_3= \partial_x^3 + \frac{3}{4}\bigg(\partial_x\circ u+  u\circ  \partial_x\bigg) + w, 
$$
$$
u=u(x,y,t), \quad w=w(x,y,t),\quad \partial_xw(x,y,t) = \frac{3}{4}  \partial_yu(x,y,t).  
$$
The IST method allows to construct many different classes of KP solutions. In this paper we focus on two such classes:
\begin{enumerate}
\item Multiline solitons.
\item Finite-gap solutions.
\end{enumerate}

\subsection{Multiline solitons}
Multiline solitons can be written explicitly in terms of Wronskian determinants \cite{Matveev1979}. Let ${\mathcal K}$ be a set of phases $\kappa_1$, \ldots, $\kappa_n$ (real or complex), and  $A =( A^i_j )$, $1\le i \le k$, $1\le j \le n$, be a $k\times n$ matrix, $k<n$,  $\mathrm{rank} A=k$. Following \cite{Matveev1979}, let us define  $k$ linearly independent solutions for the zero-curvature ``vacuum'' problem 
(\ref{eq:KP-Lax}) with $u\equiv w \equiv 0$ 
$$
f^{(i)}(x,y,t) = \sum_{j=1}^n A^i_j \exp\left(\kappa_j x +\kappa_j^2 y + \kappa_j^3 t\right), \ \ i=1,\ldots,k.
$$
Then
\begin{equation}\label{eq:KPsol}
  u(x,y,t) = 2\partial_{x}^2 \log(\tau (x,y,t)), \ \ \mbox{where} \ \ \tau (x,y,t) = \mathrm{Wr}_{x}(f^{(1)},\dots, f^{(k)}),
\end{equation}
is a multiline soliton solution to (\ref{eq:KP}). Here $\mathrm{Wr}_{x}(f^{(1)},\dots, f^{(k)})$ denotes the Wronskian with respect to the variable $x$.

From the Cauchy-Binet formula it follows that
\begin{equation}\label{eq:KP-CB}
  \tau (x,y,t)=  \sum\limits_{I} \Delta_I (A)\prod_{\mycom{i_1<i_2}{ i_1,i_2 \in I}} (\kappa_{i_2}-\kappa_{i_1} ) \,
  \exp\left( \sum\limits_{i\in I} (\kappa_i x + \kappa_i^2 y + \kappa_i^3 t ) \right),
\end{equation}
where the sum is over all $k$--element ordered subsets $I$ in $\{1,2,\ldots,n\}$, i.e. $I=\{ 1\le i_1<i_2 < \cdots < i_k \le n\}$ and $\Delta_I (A)$ are the maximal minors of the matrix $A$ with respect to the columns $I$. Since linear operations on the rows of $A$ preserve the KP multisoliton solution $u(x,y,t)$ in (\ref{eq:KPsol}), this solution depends only on the \textbf{point of the Grassmannian} $Gr(k,n)$ represented by the matrix $A$. Let us recall that the minors $\Delta_I (A)$ are exactly the Pl\"ucker (projective) coordinates for this point in $Gr(k,n)$.

Therefore the IST data for multiline soliton solutions are $({\mathcal K},[A])$, where:
\begin{enumerate}
\item ${\mathcal K}$ is a set of phases $\kappa_1$, \ldots, $\kappa_n$ (real or complex). If all  $\kappa_j$ are real, it is assumed that $\kappa_1<\kappa_2< \ldots<\kappa_n$.
\item A point $[A]$ of the Grassmannian $Gr(k,n)$.
\end{enumerate}

\subsection{Finite-gap solutions}
\label{sec:2.2}

The finite-gap method, first applied to KP equation by Krichever \cite{Krichever1976}, generates spatially periodic and quasiperiodic solutions. Let us fix a normalization time $(x_0,y_0,t_0)$.  Following \cite{Krichever1976}, the IST data are:
\begin{enumerate}
\item A finite genus $g$ compact Riemann surface $\Gamma$ with a marked point $P_0$.
\item A local parameter $1/\zeta$ near $P_0$ such that $\zeta(P_0)=\infty$.
\item A non-special divisor $\mathcal D=\gamma_1+\ldots+\gamma_g$, of degree $g$ in $\Gamma$, where $\gamma_j=\gamma_j(x_0,y_0,t_0)$.
\end{enumerate}
For generic IST data there exists a unique function $\hat\psi(\gamma,x,y,t)$, $\gamma\in\Gamma$ with the following properties for fixed $(x,y,t)$:
\begin{enumerate}
\item $\hat\psi(\gamma,x,y,t)$ is meromorphic in $\gamma$ on $\Gamma\backslash P_0$.
\item $\hat\psi(\gamma,x,y,t)$ is holomorphic in $\gamma$ outside the points $P_0$, $\gamma_1$, \ldots, $\gamma_g$ and it may have at most first order poles at the points $\gamma_j$.
\item At the point $P_0$ the function $\hat\psi(\gamma,x,y,t)$ has an essential singularity such that the function
$$
 \Xi(\gamma,x,y,t) = \hat\psi(\gamma,x,y,t)\cdot \exp\left(-\zeta (x-x_0) - \zeta^2 (y-y_0) -  \zeta^3 (t-t_0) \right),
$$
is holomorphic in $\gamma$ at the point $P_0$, and
$$
\Xi(\gamma,x,y,t) = 1 + o(1) \ \ \mbox{as} \ \ \gamma\rightarrow P_0.
$$
\end{enumerate}
Denote by $\chi_j(x,y,t)$ the Taylor coefficients for $\Xi(\gamma,x,y,t)$ near $P_0$:
$$
\Xi(\gamma,x,y,t) = 1+ \frac{\chi_1(x,y,t)}{\zeta} +\frac{\chi_2(x,y,t)}{\zeta^2}+\ldots.
$$    
Then the function
$$
u(x,y,t) = 2\partial_x \chi_1(x,y,t)
$$
satisfies the KP equation. This solution can be written explicitly using the Its-Matveev-Krichever formula:
$$
u(x,y,t) = 2\partial_x^2 \log \theta (x \vec W_1+y \vec W_2+t \vec W_3 + \vec C)+2 \hat \omega_{11},
$$
where 
$$
\theta(\vec z) = \theta(\vec z| B) = \sum\limits_{(n_1,\ldots,n_g) \in {\ZZ}^g} \exp\left[\frac{1}{2}\sum_{k,l=1}^g B_{kl}n_kn_l + \sum_{k=1}^g z_k n_k \right],
$$
is the Riemann theta-function in the normalization of \cite{Fay1973}, $\vec W_1$, $\vec W_2$, $\vec W_3$, $\hat\omega_{11}$, $B_{kl}$ are defined in terms of the spectral curve, $\vec C$  is defined in terms of the spectral curve and the divisor.

Let us recall that for a generic time $(x,y,t)$ the function $ \hat\psi(\gamma,x,y,t)$ considered as a function of $\gamma$ has $g$ zeroes $\gamma_1(x,y,t)$, \ldots, $\gamma_g(x,y,t)$ on $\Gamma$. If $x=x_0, y=y_0, t=t_0$  then $\hat\psi(\gamma,x_0,y_0,t_0)\equiv1$, and its divisor of poles coincides with the divisor of zeroes. For this reason, in the following we use the same notation for both divisors. To replace the normalization time $(x_0,y_0,t_0)$ by a new one $(x_1,y_1,t_1)$ it is sufficient to introduce a new wave function $\hat\Psi(\gamma,x,y,t)$:
\begin{equation*}
\hat\Psi(\gamma,x,y,t) = \frac{\hat\psi(\gamma,x,y,t)}{\hat\psi(\gamma,x_1,y_1,t_1)}.
\end{equation*}
It is clear that the divisor of poles for $\hat\Psi(\gamma,x,y,t)$ is $\tilde{\mathcal D}=\gamma_1(x_1,y_1,t_1)+\ldots+\gamma_g(x_1,y_1,t_1)$. From this observation it follows that:
\begin{lemma}\label{lem:norm_time}
  Let  ${\mathcal D}$ be the divisor of poles for the Krichever KP wave function with normalization time $(x_0,y_0,t_0)$, and let the divisor $\tilde{\mathcal D}$ be the divisor of zeroes for this function calculated at the point  $(x_1,y_1,t_1)$.

Then the two pairs of data: $\big( {\mathcal D}, (x_0,y_0,t_0)\big)$ and $\big( \tilde{\mathcal D}, (x_1,y_1,t_1)\big)$      generate the same KP solution. Equivalently, they correspond to the same point of the Sato-Segal-Wilson Grassmannian.
\end{lemma}

As it was pointed out in \cite{Krichever1975}, soliton solutions can be obtained from the finite-gap ones by degenerating the spectral curves, but obtaining a finite-gap solution such that its limit is the given multiline one is a highly  non-trivial task. 

\subsection{Real regular solutions}\label{RealRegular}

If the KP equation is used as a model of wave propagation in weakly dispersive and weakly non-linear media, it is necessary to select \textbf{real and regular} solutions. 

As it was pointed out in \cite{Malanyuk1991}, to obtain real and regular multiline solitons it is sufficient to start with the following IST data:
\begin{enumerate}
\item A set of $n$ real phases ordered in the ascendant way:  $\kappa_1<\kappa_2<\cdots < \kappa_n$. 
\item A point of a real totally non-negative Grassmannian $\GTNN$ (see Definition~\ref{def:tp} in the Appendix). 
\end{enumerate}

Indeed if these conditions are fulfilled, (\ref{eq:KP-CB}) is a sum of non-negative terms of which at least one is strictly positive. Therefore for real $x,y,t$, the function $\tau(x,y,t)$ is smooth and strictly positive, and $u(x,y,t)$ is real, smooth and bounded.

In \cite{KodamaWilliams2013} Kodama and Williams proved that these conditions are also necessary, and this proof is highly non-trivial. 

The multiline solitons are the second logarithmic derivative of a sum of exponents. If there is only one dominant exponent in some space-time region, then the solution is close to zero. At the boundary between two regions there are two dominant exponents approximately of the same order, and the solution is approximately a line soliton. These line solitons may form very nontrivial patterns, see \cite{BiondiniChakravarty2006}, and, following \cite{KodamaWilliams2014}, one can treat the maxima of the solutions as tropical curves in the sense of the book \cite{MikhalkinRau2018}.

For the finite-gap KP solutions the sufficient conditions for being real and regular were found by Dubrovin \cite{Dubrovin1982}; in \cite{DubrovinNatanzon1988} Dubrovin and Natanzon proved that for regular curves these conditions are also necessary.

To construct real regular KP solutions associated to smooth curves it is necessary and sufficient that:
\begin{enumerate}
\item The curve $\Gamma$ is an ${\mathtt M}$-curve, i.e. 
\begin{enumerate}
\item $\Gamma$ is real, that is it admits an antiholomorphic involution $\sigma:\Gamma\rightarrow\Gamma$.
\item The involution $\sigma$ has the maximal possible number $g+1$ of fixed ovals. The letter ${\mathtt M}$ means ``maximal''.
\end{enumerate}     
\item The marked point $P_0$ is real, i.e. $\sigma P_0= P_0$. We denote the oval containing $P_0$ by $\Omega_0$, and the other ones by $\Omega_1$, \ldots, $\Omega_g$. We call the oval $\Omega_0$ \textbf{infinite} and all other ovals  \textbf{finite}.
\item The local parameter $1/\zeta$ satisfies $\sigma\zeta=\bar\zeta$.
\item Each finite ovals $\Omega_j$, $j=1,\ldots,g$ contains exactly one divisor point, which we denote by $\gamma_j$.
\end{enumerate}

We call the divisors satisfying these conditions \textbf{Dubrovin-Natanzon (DN) divisors}. For such real and regular solutions the divisor of zeroes ${\mathcal D}(x,y,t)=\gamma_1(x,y,t)+\ldots+\gamma_g(x,y,t)$ is also a DN-divisor for all real $(x,y,t)$, see \cite{Dubrovin1982}, \cite{DubrovinNatanzon1988}.

\subsection{Real and regular KP solitons and DN divisors on tropical limits of ${\mathtt M}$-curves}
\label{sec:2.4}

The relevance of obtaining real regular soliton solutions from real regular finite-gap solutions was pointed out by S.P. Novikov. In our papers \cite{AG1,AG2,AG3,AG4} we proved that every \textbf{real and regular} KP soliton solution can be obtained from  \textbf{real and regular} KP finite-gap solutions by a proper degeneration of the spectral ${\mathtt M}$-curve into an  ${\mathtt{MM}}$-curve. These degenerated spectral curves can be treated as tropical curves in the sense of the papers \cite{AgostiniCelikStruweSturmfels2021,AgostiniFevolaMandelshtamSturmfels2023}, see also \cite{Ichikawa2023} and  \cite{KummerSturmfelsVlad2025}.

Following Section~\ref{RealRegular}, let $[A]$ be a point in a totally non-negative Grassmannian. In \cite{Postnikov2006} it is proven that the totally non-negative Grassmannian admits a natural cell decomposition into positroid cells (see also the Appendix). The positroid cell containing $[A]$ is represented by an equivalence class of perfectly oriented bicolored graphs in the disk. This class contains the canonical Le-graph which has the property of having total number of faces equal to $1+d$ where $d$ is the dimension of the positroid cell. If the positroid cell belongs to $\GTNN$, then the Le-graph has $d-k$ trivalent white vertices.

In \cite{AG4} we proved that to any graph in Postnikov class representing the positroid cell containing $[A]$ we can associate an ${\mathtt{MM}}$-curve, which is dual to this graph. Such ${\mathtt{MM}}$-curve is a spectral curve for the soliton data $\{[A], \kappa_1<\ldots<\kappa_n\}$. For a generic normalization time, the divisor corresponding to the soliton data is \textbf{smooth}, i.e. it lies outside the singular points of the curve, and is explicitly obtained in \cite{AG3,AG4} by solving the linear system of relations studied in \cite{AG5,AG6}.

Let us recall the main construction from \cite{AG4} for the case of interest in this paper.

In the special case of the Le-graph the arithmetic genus $g$ of the spectral curve coincides with the dimension $d$ of the correspondent positroid cell, and the number of trivalent white vertices (respectively of trivalent black vertices) is $d-k$ (respectively $d-n+k$). The corresponding  ${\mathtt{MM}}$-curve $\Gamma$ is the union of $2d-n+1$ copies of $\mathbb{CP}^1$. One of these copies $\Gamma_0$ corresponds to the boundary of the disk and the remaining $2d-n$ copies correspond to internal vertices. We denote $\Gamma_j$ the copy of $\mathbb{CP}^1$ corresponding to the white vertex $W_j$, and denote $\Sigma_l$ the one corresponding to the black vertex $B_l$ respectively. The component $\Gamma_0$ contains the marked point $P_0$ and the points $\kappa_1$, \ldots, $\kappa_n$ corresponding to the boundary vertices of the graph.

\begin{definition} In the following, we call the $\mathbb{CP}^1$ components corresponding to the white vertices \textbf{white}, and the $\mathbb{CP}^1$ components corresponding to the black vertices \textbf{black}. We also call $\Gamma_0$ the \textbf{boundary component}.
\end{definition}

\begin{figure}[H]
  \centering\includegraphics[width=0.9\textwidth]{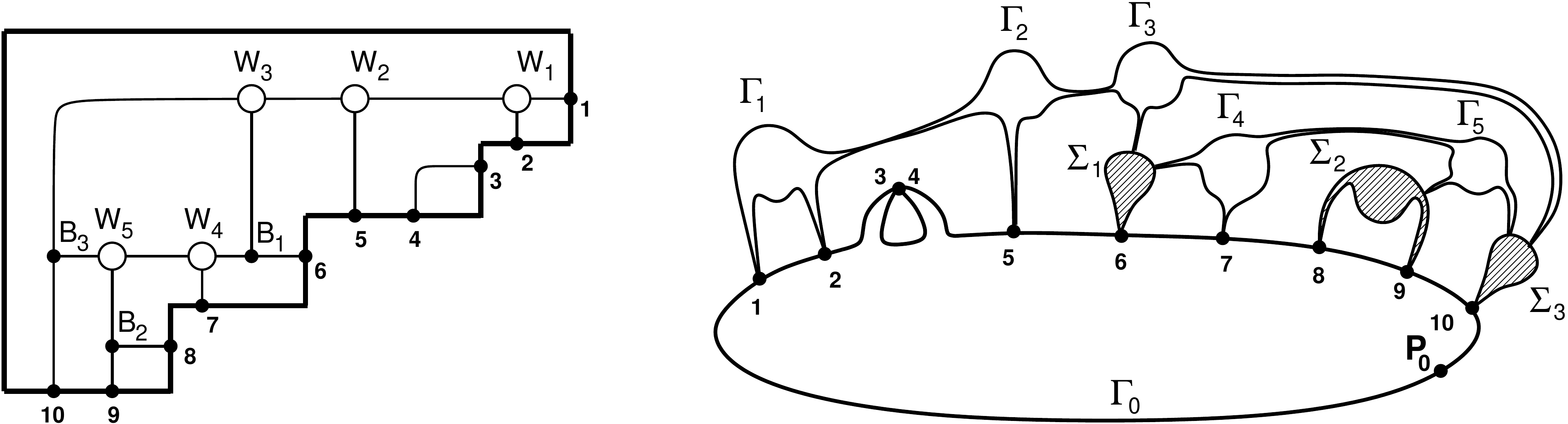}
    \caption{\small{\sl On the left: the Le-diagram discussed in the Appendix (see Figure~\ref{fig:appendix_7}). On the right: the real part of the corresponding curve $\Gamma$. $\Gamma$ is the union of 9 copies of $\mathbb{CP}^1$ ($\Gamma_0$, $\Gamma_1$, $\Gamma_2$, $\Gamma_3$, $\Gamma_4$, $\Gamma_5$, $\Sigma_1$, $\Sigma_2$, $\Sigma_3$). The component $\Gamma_0$ corresponds to the boundary of the disk, the other components correspond to the vertices of the graph. The boundary vertices 3 and 4 correspond to a nodal point. All other boundary vertices correspond to intersections of $\Gamma_0$ with other components.}}
\label{fig:dual}
\end{figure}

The edges of the graph correspond either to intersections between distinct components or to nodal points (see \cite{AG3,AG4} for more details). For an example see Figure~\ref{fig:dual}. We assume that on $\mathbb{CP}^1$ components associated with 3-valent vertices the coordinates of the double points are $0$, $1$ and $\infty$ (see \cite{AG3,AG4} for details). 
 
\begin{figure}[H]
  \centering\includegraphics[width=0.9\textwidth]{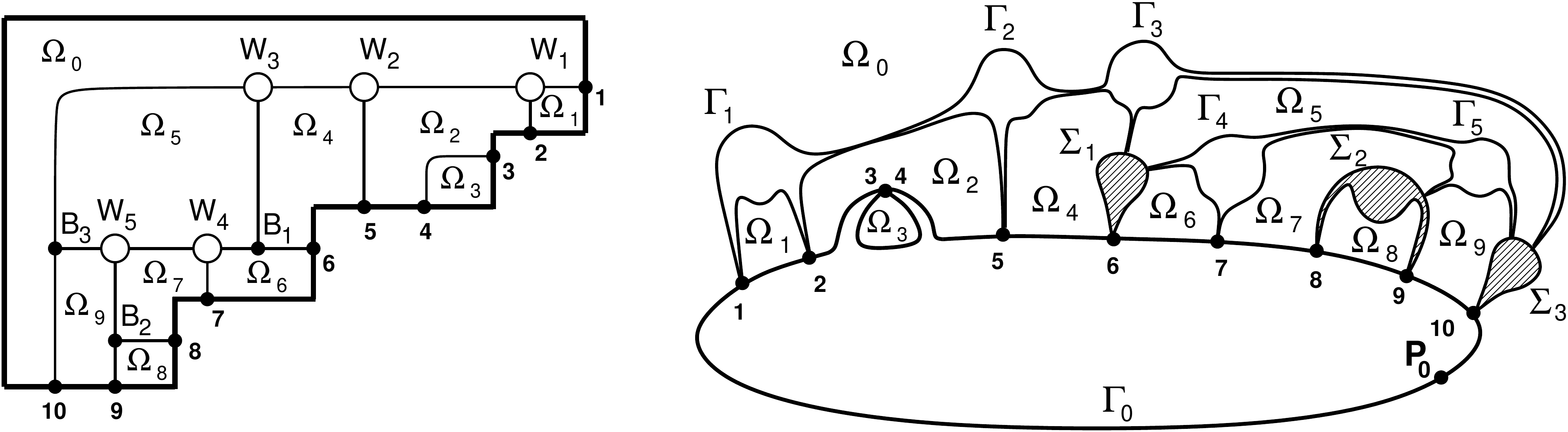}
  \caption{\small{\sl Top: we mark the faces of the Le-diagram of Figure~\ref{fig:dual}.  Bottom: the real ovals in $\Gamma$ are the boundaries of the marked areas with the same labels. The infinite oval $\Omega_0$ is the boundary of the infinite area and contains the points $1$, $10$ and $P_0$.}}
  \label{fig:dual_2}
\end{figure}

In Section~4.2 in \cite{AG3}, see also Proposition~4.4 in \cite{AG4}, we proved that $\Gamma$ is the rational degeneration of a genus $g$ ${\mathtt M}$-curve, where the $g+1$ ovals correspond to the faces of the graph (see Figure~\ref{fig:dual_2} for an example), i.e. $\Gamma$ is an  ${\mathtt{MM}}$-curve in the sense of \cite{KummerSturmfelsVlad2025}.

In \cite{AG3,AG4} we also proved that the KP wave function is extended from $\Gamma_0$ to $\Gamma$ and has the following properties:
\begin{enumerate}
\item It is a rational function of degree 1 in the spectral parameter on each white component, and is constant on each black one.
\item The value of the KP wave function at the nodes and intersection points is computed by solving linear relations on the corresponding edges of the graph. These systems of relations and their solutions are treated in details in \cite{AG5,AG6}.
\item For a generic time $(x_0,y_0,t_0)$ the divisor of poles of the KP wave function is smooth, and each finite oval contains exactly one divisor point. Therefore the divisor is a Dubrovin-Natanzon (DN) divisor, and the solution can be obtained as the limit of a real and regular finite-gap solution, associated to a non-singular ${\mathtt M}$-curve. 
\end{enumerate}
  
In this paper we are interested in parametrizing the full positroid cell by DN divisors; therefore it is necessary to fix $(x_0,y_0,t_0)$  and use both smooth and singular divisors. We treat this problem in the following Sections.

\section{Singular  Dubrovin-Natanzon divisors on ${\mathtt{MM}}$-curves}

In this section we extend the notion of DN divisors to the singular ones. In particular, we explain how to pass from the the language of Weil divisors to that of Cartier divisors in order to obtain a proper description of the singular component.

To define the Cartier divisor of zeros of the wave function, it is necessary to assume that this function does not vanish identically at any connected component of the curve. Therefore we have to modify the definition of the spectral curve by removing all black components. 

Next we recall the definition of the Abel transform for nodal curves, and we study for which singular divisor configurations it may remain finite. We show that we have two basic configurations, and if we pass from Weil divisors to the Cartier ones, the Abel transform became regular near these configurations. 

Finally, we introduce Dubrovin-Natanzon Cartier divisors, and we explain how to modify the Krichever construction for singular divisors.

\subsection{Reduced spectral curve}\label{sec:reduced_curves}

To treat singular DN divisors as Cartier ones, we have to use spectral curves obtained from those constructed in the previous section by applying the following procedure.

\begin{figure}[H]
  \centering\includegraphics[width=0.8\textwidth]{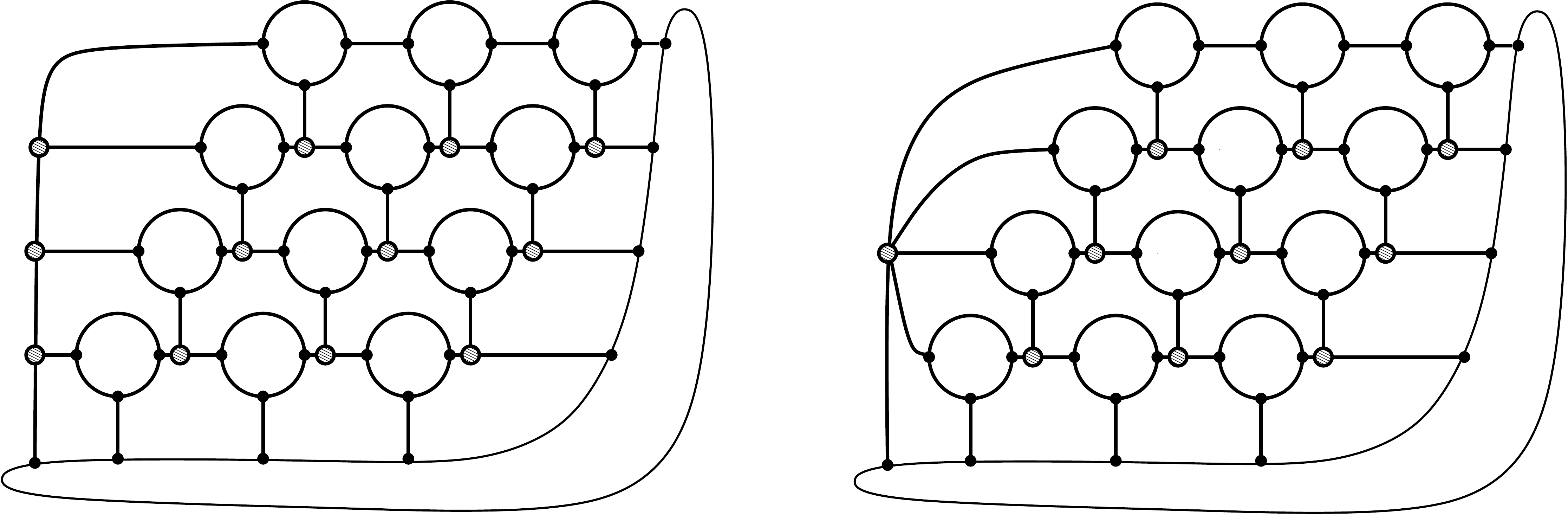}
    \caption{\small{\sl On the left: the original graph. On the right: the graph obtained by merging a chain of 3 trivalent black vertices into one black vertex of valency 5.}}
\label{fig:reduct1}
\end{figure}

\begin{figure}[H]
  \centering\includegraphics[width=0.5\textwidth]{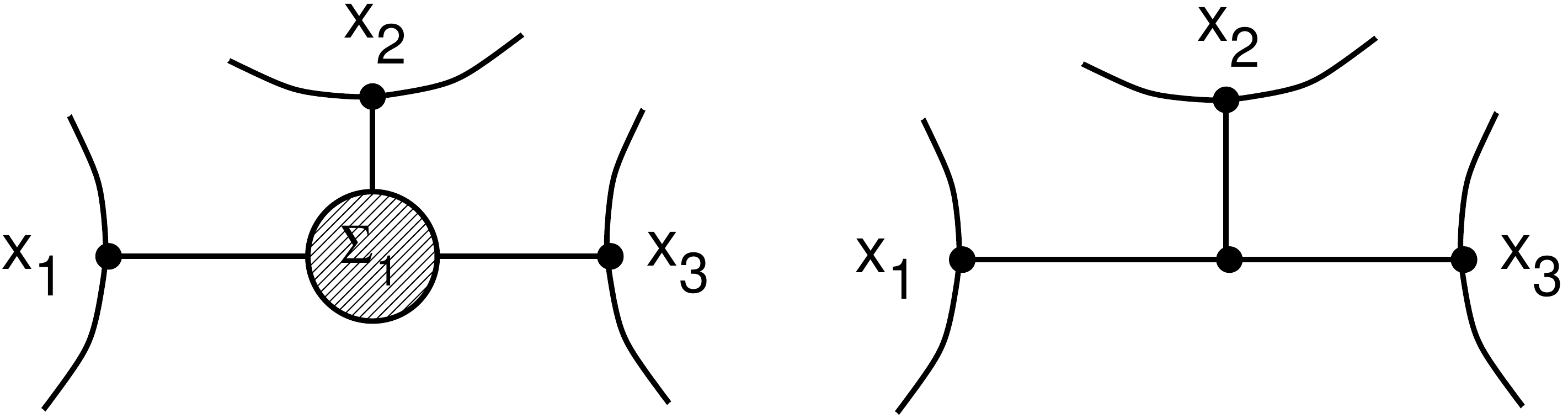}
    \caption{\small{\sl On the left: the original curve. On the right: instead of the component $\Sigma_1$ we have a triple point, i.e. for regular functions we have the relation $f(x_1)=f(x_2)=f(x_3).$ }}
\label{fig:reduct2}
\end{figure}

\begin{enumerate}
\item If we have a chain of $l$ trivalent black vertices we use Postnikov move and transform the chain to a unique black vertex of valency $l+2$ (see Figure~\ref{fig:reduct1}). The wave function constructed in the previous Section remains unchanged because it is constant on black components. 
\item We eliminate each black component and replace it by a point of multiplicity equal to the valency of the corresponding black vertex in the reduced graph. Multiple point means that, if the black component was connected to the points $x_1$, \ldots, $x_{l+2}$, on the new curve a function $f(x)$ is regular near this point, if it is regular near each $x_j$ and satisfies the algebraic identity:
\begin{equation}\label{eq:multiple}
f(x_1)=f(x_2)=\ldots=f(x_{l+2}),
\end{equation}
see Figure~\ref{fig:reduct2}. 

Let us remark that the basic algebraic model of such multiple point is the intersection of $l+2$ generic lines at one point in $\CC^{l+2}$, not in $\CC^2$. 

\end{enumerate}  

It is easy to check that, for smooth divisors, the wave function defined in the previous section remains the same on the reduced spectral curve. Therefore it is regular in the spectral parameter outside the divisor of poles and the essential singularity at $P_0$. 

\begin{figure}[H]
  \centering\includegraphics[width=\textwidth]{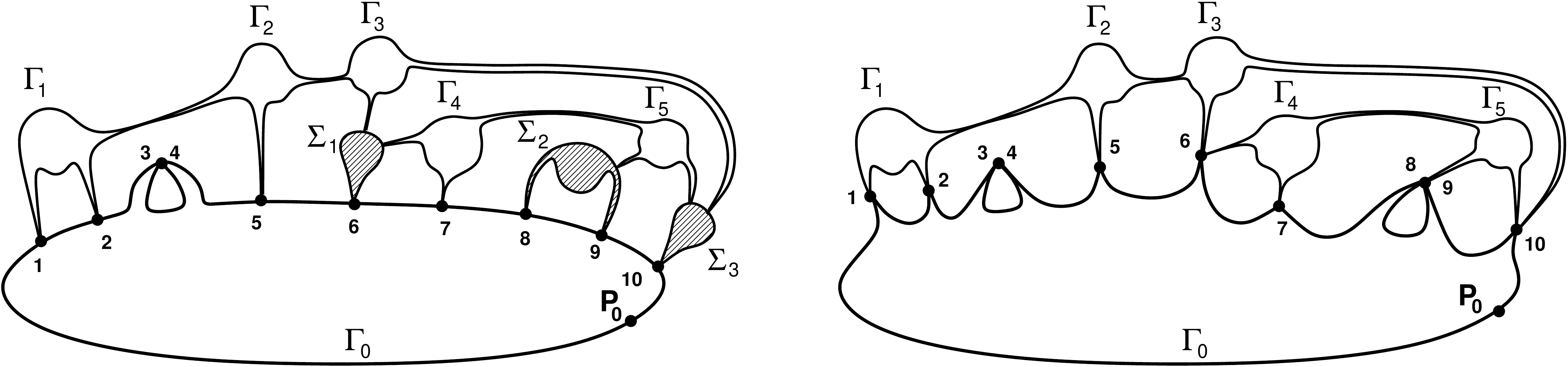}
    \caption{\small{\sl On the left: the original spectral curve corresponding to the Le-diagram of Figure~\ref{fig:dual} left. On the right: the corresponding reduced spectral curve.}}
\label{fig:dual_new}
\end{figure}

In Figure~\ref{fig:dual_new} we apply the reduction procedure to the curve of Figure~\ref{fig:dual}.
  
\begin{remark}
If we apply the reduction procedure to an ${\mathtt{MM}}$-curve, we again obtain an  ${\mathtt{MM}}$-curve.
\end{remark}

\subsection{The Abel transform of Dubrovin-Natanzon divisors on ${\mathtt{MM}}$-curves}

Let us start by defining the Abel transform for singular curves following \cite{Tyurin2003,Artamkin2004}.
 
\begin{definition}\label{def:holomorhic}
A differential $\omega$ is holomorphic on the reduced ${\mathtt{MM}}$-curve if:
\begin{enumerate}
\item Its restriction to each component is meromorphic.
\item At each component this restriction is regular outside the multiple points and has at most first order poles at the multiple points.
\item If the points $x_1$, \ldots, $x_{l+2}$ are glued together, the sum of residues at these points is equal to 0.
\end{enumerate}
\end{definition}

It is convenient to choose $b$-cycles corresponding to the finite faces of the plabic graph, and to orient them clockwise.  To simplify notations, we enumerate them from  $1$ to $g$, where $g$ denotes the arithmetic genus of $\Gamma$. The $c$-cycles go around the edges. The $a$-cycles are integer linear combinations of $c$-cycles.

\begin{figure}[H]
  \centering\includegraphics[width=0.6\textwidth]{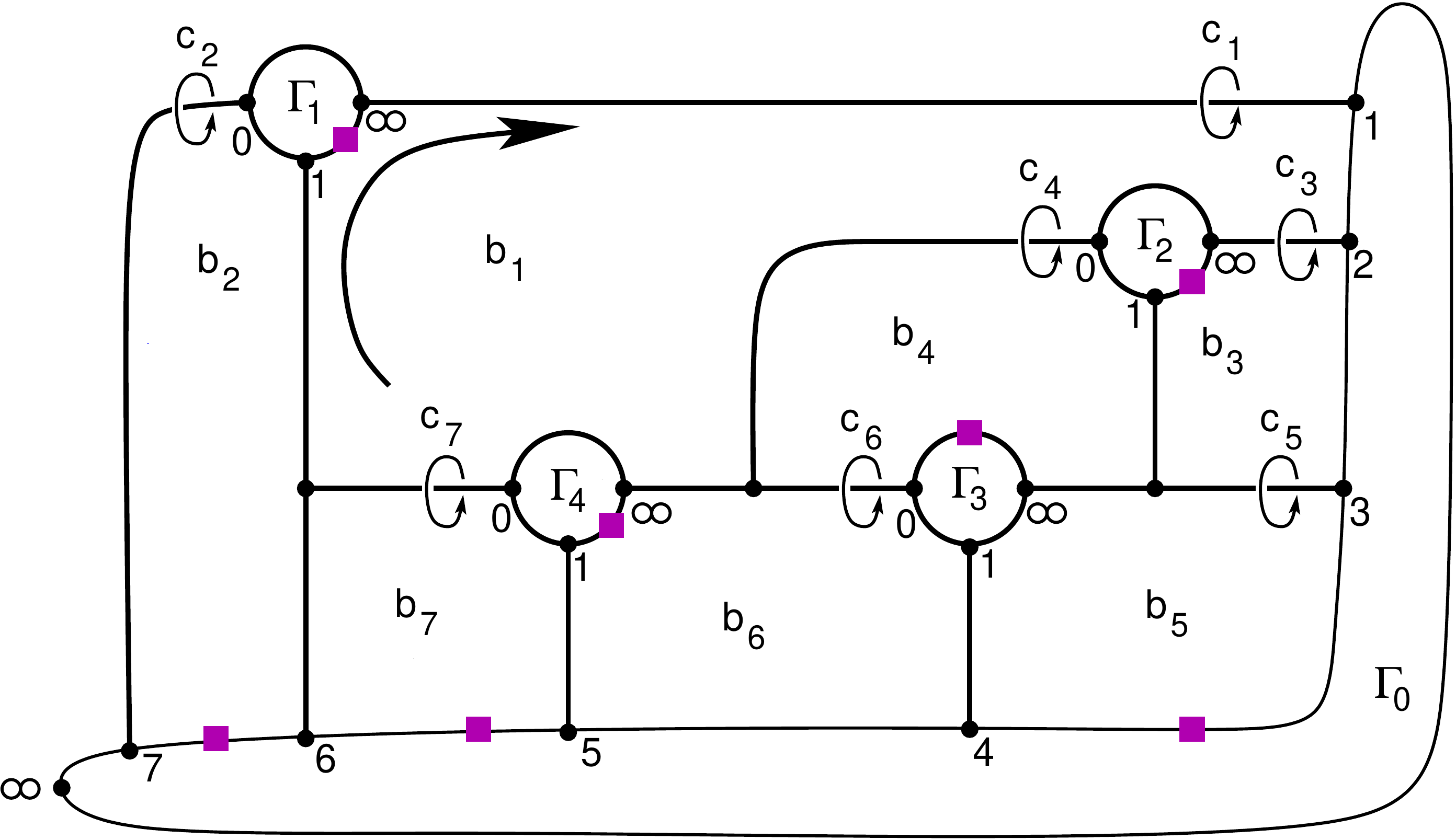}
      \caption{\small{\sl A basis of $b$-cycles and $c$-cycles on $\Gamma$}}
			\label{fig:basis2}
\end{figure}

\begin{example}
For the graph in Figure~\ref{fig:basis2} we have the following expressions:
\begin{align*}
  c_1 = a_1, \ \ c_2 = a_2,  \ \ c_3 = a_3 - a_1 \ \ c_4 = a_4 - a_1 \ \ c_5 = a_5 - a_3, \\
  c_6 = a_6 - a_4, \ \  c_7 = a_7 - a_1,
\end{align*}  
therefore
\begin{align*}
  a_1 = c_1, \ \ a_2 = c_2,  \ \ a_3 = c_3 + c_1 \ \ a_4 = c_4 + c_1 \\
  a_5 = c_5 + c_3 + c_1, \ \  a_6 = c_6 + c_4 + c_1 , \ \  a_7 = c_7 + c_1,
\end{align*} 
\end{example}

In our text we use the following normalization for the canonical differentials $\omega_k$:
\begin{equation*}
  \oint_{a_j}\omega_k = 2\pi i \delta_{jk},
\end{equation*}
as in \cite{Fay1973}. The differentials $\omega_k$ are defined in the following way (see \cite{Artamkin2004}):
\begin{enumerate}
\item If the cycle $b_k$ does not pass through a white vertex, then the restriction of $\omega_k$ to the corresponding copy of $\mathbb{CP}^1$ is equal to zero.
\item If the cycle $b_k$ does not pass through any boundary vertex, then the restriction of $\omega_k$ to $\Gamma_0$ is equal to zero.  
\item If the cycle $b_k$ passes through a white vertex, the differential  $\omega_k$ has exactly two first-order poles on the corresponding component, with residue $+1$ (respectively $-1$) at the point associated to the outgoing (respectively incoming) edge.
\item If the cycle  $b_k$ passes through a boundary vertex, $\omega_k$ has a first-order pole at the corresponding point of $\Gamma_0$ with residue $1$  (respectively $-1$) if $b_k$ leaves (respectively enters) the boundary of the disk at such vertex. 
\end{enumerate}  

\begin{example}
For the graph in Figure~\ref{fig:basis2} the differential $\omega_1$ has the following properties:
\begin{enumerate}
\item It is equal to zero at $\Gamma_3$.
\item At all other components it has two simple poles with residues $+1$ and $-1$. Let us compute $\omega_1$ restricted to $\Gamma_1$. The intersection of the cycle $b_1$ with the component $\Gamma_1$ is the line connecting the points $1$ and $\infty$. Therefore the restriction of the holomorphic form $\omega_1$ to the component $\Gamma_1$ is:
  $$
  \omega_1\bigg|_{\Gamma_1}=-\frac{dz}{z-1}.
  $$
Proceeding in the same way on the other components one concludes that:
\begin{align*} 
  \omega_1&= \frac{dz}{z}& \ \ &\mbox{on} \ \ \Gamma_4,\\
  \omega_1&= \frac{dz}{z}& \ \ &\mbox{on} \ \ \Gamma_2,\\
  \omega_1&= \frac{dz}{z-\kappa_2}- \frac{dz}{z-\kappa_1} & \ \ &\mbox{on} \ \ \Gamma_0.\\
\end{align*}  
\end{enumerate}

\begin{remark}
Symbolically we draw the $c$-cycles around the edges. Since every edge is connected with a boundary vertex or a white vertex, the cycle shall be drawn on the corresponding component, for an example see Figure~\ref{fig:cycleb4}. 
\end{remark}  

\begin{figure}[H]
  \centering\includegraphics[width=0.6\textwidth]{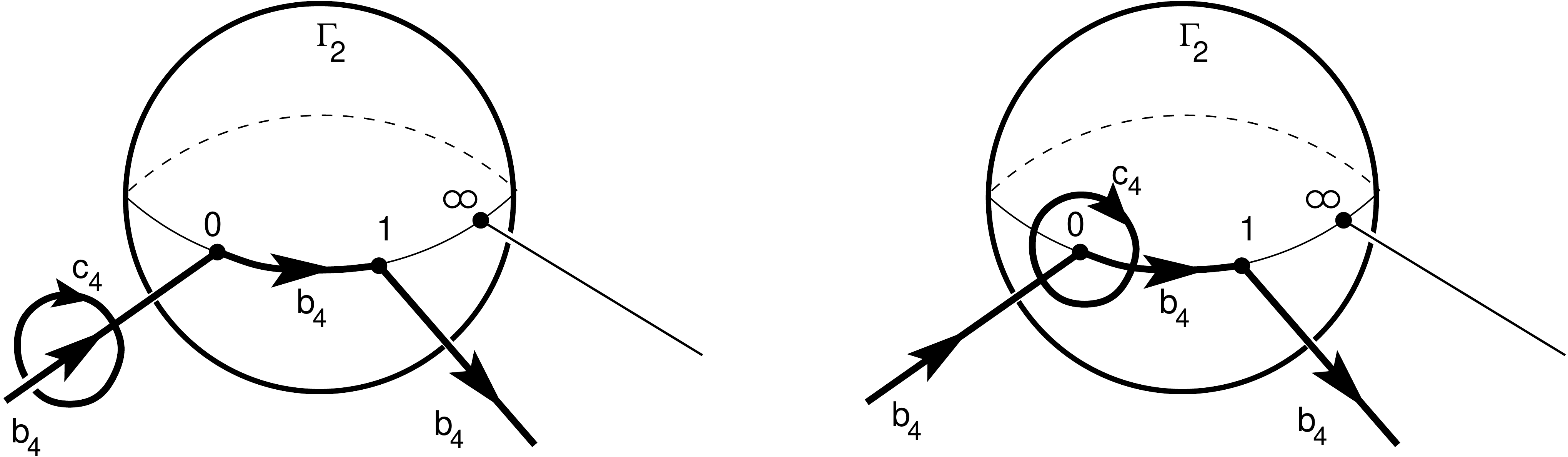}
  \caption{\small{\sl On the left: the symbolic representation of the cycle $c_4$. On the right: the cycle $c_4$ is drawn on $\Gamma_2$.}}
	\label{fig:cycleb4}
\end{figure}

\end{example}

It is convenient to define the Abel transform of the divisor in a slightly non-standard way. 
\begin{definition}\label{def:1}
  Let ${\mathcal D} = \gamma_1+\ldots +\gamma_g$, be a smooth divisor where $g$ is the arithmetic genus of $\Gamma$. Let us recall that each component $\Gamma_j$ with  $j>0$  contains exactly one divisor point. We define the $l$-th component of the Abel transform of ${\mathcal D}$ by: 
  \begin{equation}
    A_l({\mathcal D}) = \sum_{\gamma_s\in\Gamma_0}\,\int_{\infty}^{\gamma_s}\omega_l\big|_{\Gamma_0}+  \sum_{j=1}^{g-k}\, \int_{P_{0,j}}^{\gamma_{j}}\omega_l\big|_{\Gamma_{j}},
  \end{equation}
where the starting points of integration $P_{0,j}$ are fixed using the following rule. If the cycle $b_l$ passes through the component $\Gamma_{j}$, $P_{0,j}$  denotes the unique marked point on $\Gamma_{j}$ through which it does not pass. If the cycle $b_l$ does not pass through the component $\Gamma_{j}$, the restriction of $\omega_l$ to this component it identically zero, therefore this term does not contribute to the sum, and the choice of $P_{0,j}$ is irrelevant.
\end{definition}
A divisor $\mathcal D$ satisfies the DN condition if and only if
$$
\Im A_l({\mathcal D}) = \pi \mod(2 \pi) \ \ \mbox{for all} \ \  l=1,\ldots, g.
$$

\subsection{Resolution of singularities in the space of Weil divisors}\label{sec:res_sing}

If $\Gamma$ is a regular $M$-curve, the Abel map and its inverse are both regular on the DN component, and Weil divisors are equivalent to Cartier divisors. On the contrary, in the case of rational curves such as the ${\mathtt{MM}}$-curves considered here, the Abel map is discontinuous at singular Weil divisors, therefore singular configurations need additional treatment. 

\begin{figure}[H]
  \centering \includegraphics[width=0.35\textwidth]{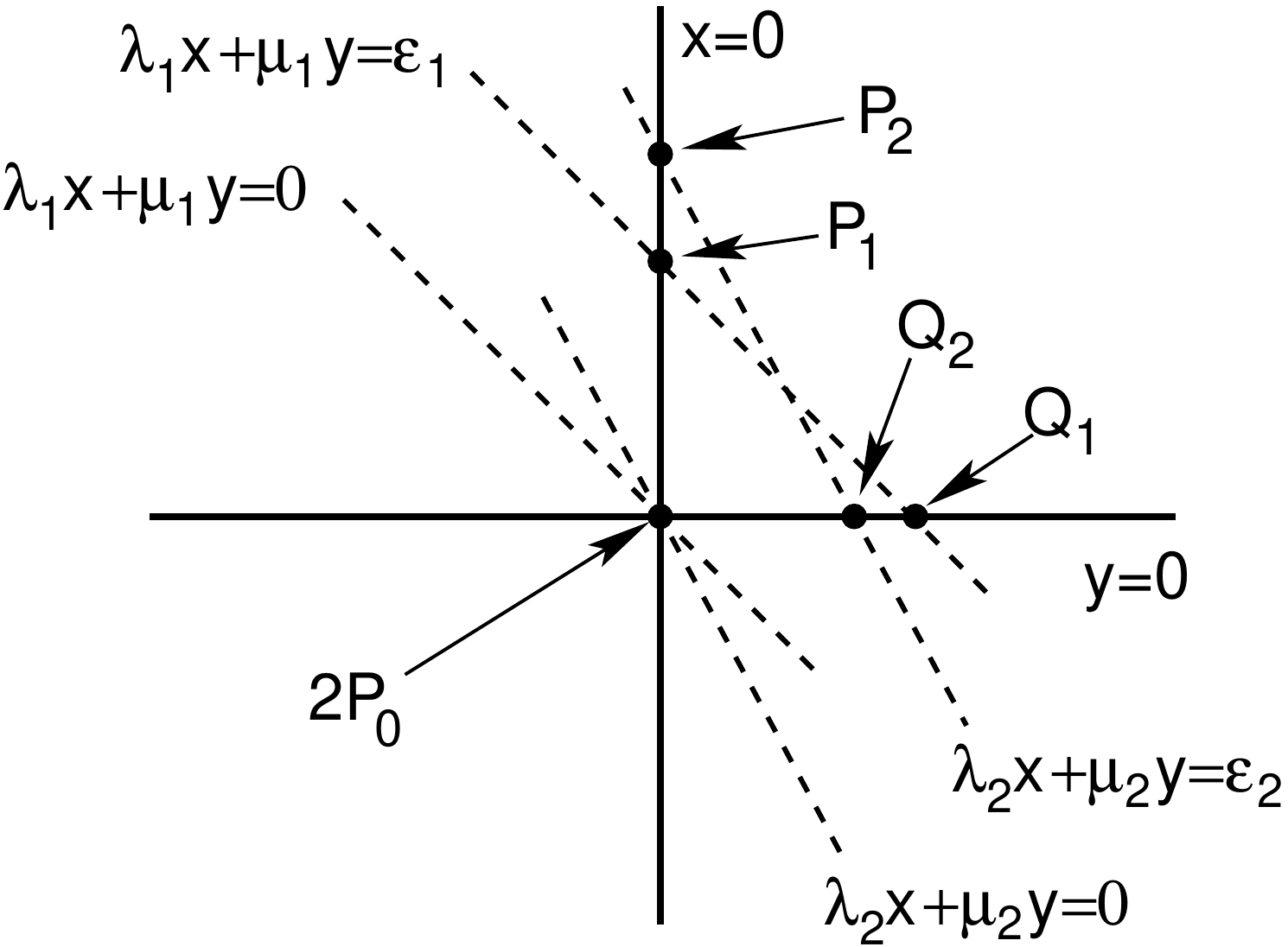}
      \caption{\small{\sl We illustrate Example~\ref{ex:3}. For $\varepsilon_j\ne 0$  Cartier divisors are completely defined by Weil ones. For $\varepsilon_j= 0$ Cartier divisors corresponding to the Weil divisor $2P_0$ form a one-parameter family parameterized by the slope of the line $f_j(x,y)=0$.}\label{fig:divisors}}
\end{figure}

Let us start illustrating the difference between Weil and Cartier divisors in the simplest examples of either a rational curve with two irreducible components, or a rational curve with three irreducible components intersecting in one point.
\begin{example}\label{ex:3}
  Let $\Gamma$ be the rational curve $xy=0$, and assume that the divisor ${\cal D}= P+Q$ is a pair of points, $P$ belonging to the component $x=0$ and $Q$ to $y=0$.

 If both points lie outside the point $x=y=0$, we may describe ${\cal D}$ as a Cartier divisor by covering $\Gamma$ with two Zariski-open sets: $U_o=\Gamma\backslash\infty$ and $U=\Gamma\backslash\{P,Q\}$, and introducing two functions: $1$ on $U$ and $f(x,y)$ on $U_0$ where 
$$ 
  f(x,y) = \lambda x + \mu y - \varepsilon, \ \ \mbox{such that} \ \ f(P)=f(Q)=0.   
$$
Equivalently, this Cartier divisor is defined by a straight line passing through the points $P$ and $Q$ (see Figure~\ref{fig:divisors}). 

If we consider the limit where $\lambda\ne0$, $\mu\ne0$, $\varepsilon\rightarrow 0$, the Weil divisor tends to $2P_0$,
where $P_0=(0,0)$, and this limit does not depend on the ratio $\lambda/\mu$. On the contrary, the limiting Cartier divisor depends on this ratio.  

Indeed, consider a pair of linear functions
 $$ 
  f_j(x,y) = \lambda_j x + \mu_j y - \varepsilon_j, \ \ j=1,2,   
 $$
  where all $\lambda_j$ and $\mu_j$ are non-zero. If $\lambda_1/\mu_1\ne\lambda_2/\mu_2$, in the limit $\varepsilon_j\rightarrow 0$, $j=1,2$,  the corresponding Cartier divisors are different, because the function $f_1(x,y)/f_2(x,y)$  
is discontinuous on $\Gamma$, since its values at the double point $x=y=0$ are different at the components $x=0$ and $y=0$ (see Figure~\ref{fig:divisors}). 

In this case, to pass from Weil divisor to Cartier divisor, one has to apply the standard blow-up procedure to the variety of Weil divisors near the point $P_j=Q_j=0$. We describe this procedure for DN Weil divisors later in this Section, see Equation (\ref{eq:cartier2}).
\end{example}

\begin{figure}[H]
  \centering \includegraphics[width=0.35\textwidth]{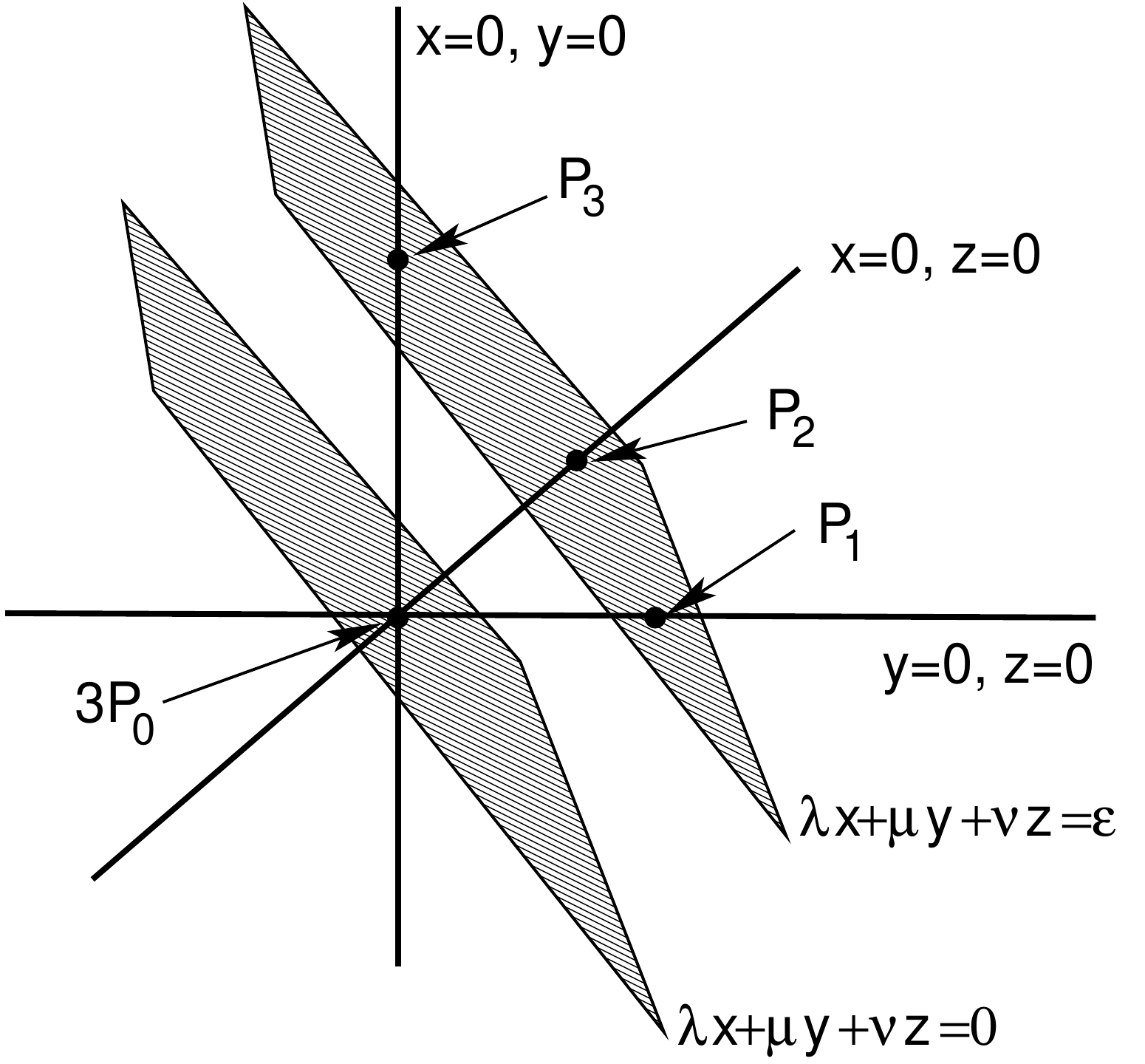}
      \caption{\small{\sl We illustrate Example~\ref{ex:4}. For $\varepsilon\ne 0$ the Cartier divisors are completely defined by the Weil ones. But for $\varepsilon= 0$ Cartier divisors corresponding to the Weil divisor $3P_0$ form a two-parameter family, parameterized by the normal vectors to the plane $f(x,y,z)=0$.}\label{fig:divisors2}}
\end{figure}

\begin{example}\label{ex:4}
  Let $\Gamma$ by the rational curve in $\CC^3$ defined by the following equations
  $$
  \left\{
\begin{split}
  xy&=0, \\
  xz&=0, \\
  yz&=0. 
\end{split}
\right.
 $$
$\Gamma$ is the union of three lines intersecting at the point $(0,0,0)$ (see Figure~\ref{fig:divisors2}). Let 
${\cal D} = P_1+P_2+P_3$ be a generic Weil divisor such that $P_1$ belongs to the component $y=0,z=0$,  $P_2$ to $x=0,z=0$,  $P_3$ to $x=0,y=0$. As a Cartier divisors $\cal D$ is locally defined by a linear function 
 $$ 
  f(x,y,z) = \lambda x + \mu y + \nu z - \varepsilon, \ \   
  $$
  where $\lambda\mu\nu\ne 0$ and $\varepsilon\ne 0$.

Instead, if $\varepsilon=0$, two functions  $f_1(x,y,z) = \lambda_1 x + \mu_1 y + \nu_1 z$ and $f_2(x,y,z) = \lambda_2 x + \mu_2 y + \nu_2 z$ define the same  Cartier divisor if and only if
$$ 
\frac{\lambda_1}{\lambda_2}=\frac{\mu_1}{\mu_2}=\frac{\nu_1}{\nu_2}.
$$
Therefore we obtain a two-parameter family of  Cartier divisors corresponding to the same Weil divisor ${\cal D} = 3 P_0$. Again, we describe the blow-up procedure for DN Weil divisors later in this Section, see Equation (\ref{eq:cartier3}).

\end{example}

Let us focus on the study of Dubrovin-Natanzon divisors on ${\mathtt{MM}}$-curves. For rational curves, the Jacobian has more than one connected component. As shown in \cite{AG3,AG4}, for a generic normalization KP time $(x_0,y_0,t_0)$ the divisor ${\mathcal D} = \gamma_1(x_0,y_0,t_0)+\ldots+\gamma_g(x_0,y_0,t_0)$ corresponding to a real and regular solution lies in the finite ovals outside intersection and nodal points. Moreover, there is exactly one divisor point in each finite oval, i.e. the divisor satisfies DN condition. The Abel transform for this divisor is finite and remains finite for all finite $(x,y,t)$. Therefore, if the Abel transform of a point lies in the open interior of one of these Jacobian components, for all finite $(x,y,t)$ the Abel transform of its KP trajectory remains inside the open interior of this component.

\begin{figure}[H]
  \centering\includegraphics[width=0.55\textwidth]{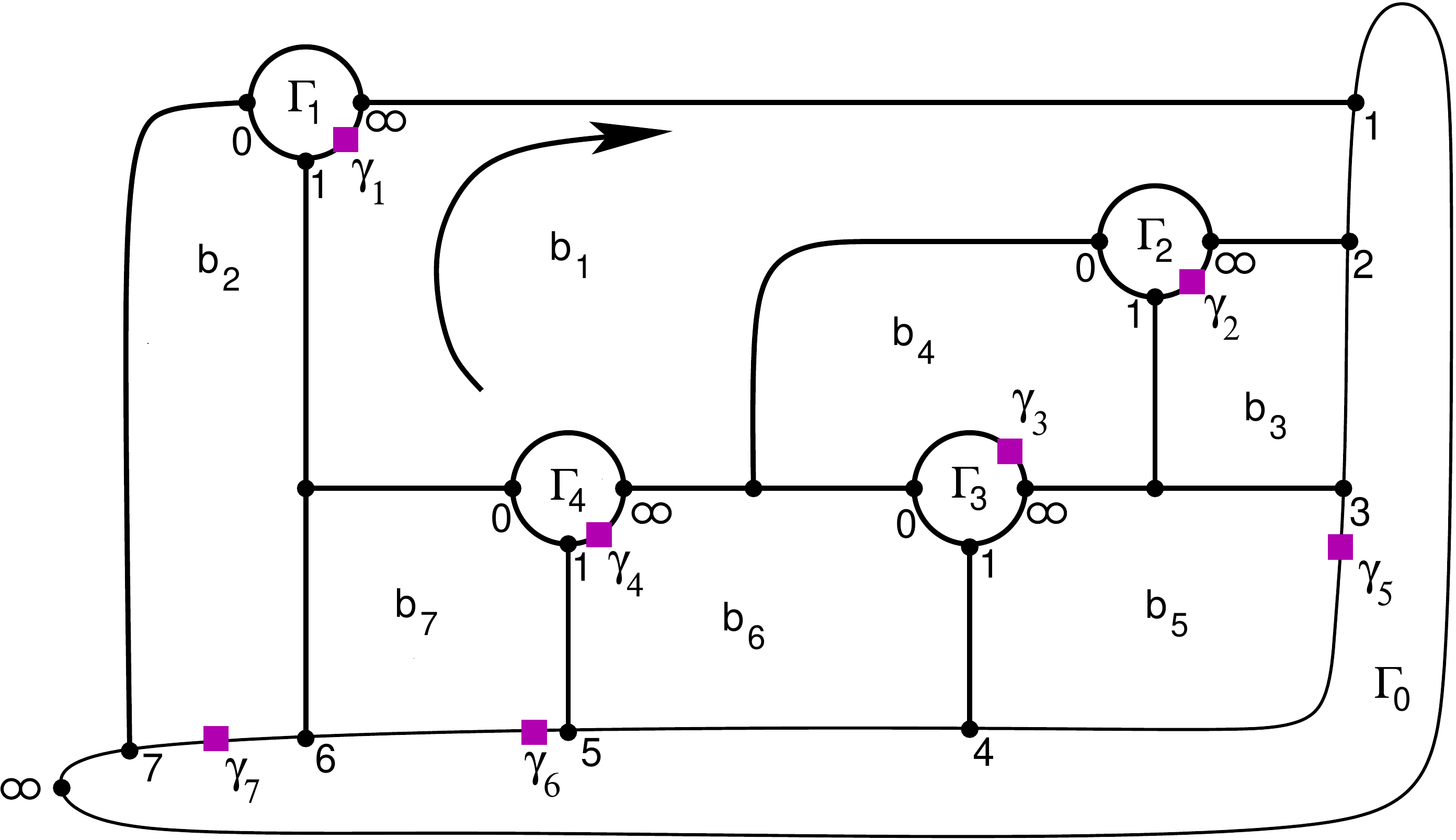}\hspace{0.5cm}
    \includegraphics[width=0.25\textwidth,trim= 0 -3cm 0 0 ]{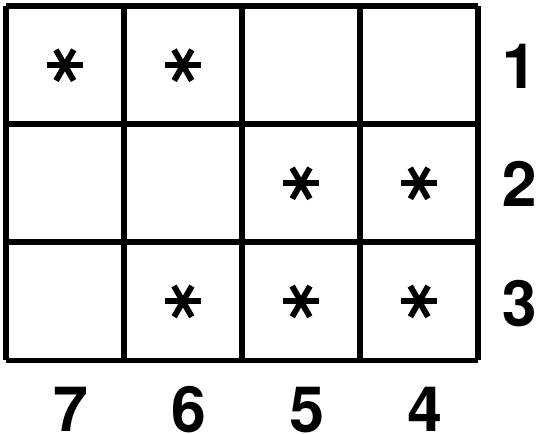}
      \caption{\small{\sl A divisor on $\Gamma$ satisfying Dubrovin-Natanzon conditions.}}\label{fig:basis3}
\end{figure}

In the following example we show that divisor points may reach nodes of the rational curve in finite time, keeping the Abel transform finite. Consider the curve in Figure~\ref{fig:basis3} (left), describing the KP-2 soliton solutions in the 7-dimensional positroid cell of $Gr^{\mbox{\tiny TNN}} (3,7)$, represented by the Le-diagram in Figure~\ref{fig:basis3} (right).

The edge connecting $\Gamma_0$ and $\Gamma_4$ represents the node between the two corresponding components of the ${\mathtt{MM}}$-curve which bound the homological cycles $b_6$ and $b_7$.  Therefore the differentials $\omega_6$, $\omega_7$  corresponding to the components $A_6({\mathcal D})$ and  $A_7({\mathcal D})$  of the Abel transform have first-order poles at this nodal point, while all other differentials in the basis are regular near this node.

We choose the local coordinate on $\Gamma_4$ to take the value 1 at the node, and  the local coordinate on $\Gamma_0$ to take the value $\kappa_5$. Therefore:

\begin{align}
&\omega_6 = -\frac{d\gamma}{\gamma-1} + O(1), \ \omega_7= \frac{d\gamma}{\gamma-1} + O(1), \ \mbox{at} \ \Gamma_4 \  \mbox{near the point} \ \gamma=1,\\
&\omega_6 = \frac{d\gamma}{\gamma-\kappa_5} + O(1), \ \omega_7= -\frac{d\gamma}{\gamma-\kappa_5} + O(1), \ \mbox{at} \ \Gamma_0 \ \mbox{near the point} \ \gamma=\kappa_5.
\end{align}  
When $\gamma_4$ and $\gamma_6$ lie near this nodal point, we have
$$
A_6({\mathcal D}) =\log\left(\frac{\gamma_6-\kappa_5}{\gamma_4-1}\right) + O(1), \ \
A_7({\mathcal D}) =\log\left(\frac{\gamma_4-1}{\gamma_6-\kappa_5}\right) + O(1),
$$
and all other components of the Abel transform are regular in $\gamma_4$, $\gamma_6$.

\begin{figure}[H]
  \centering\includegraphics[height=3.5cm]{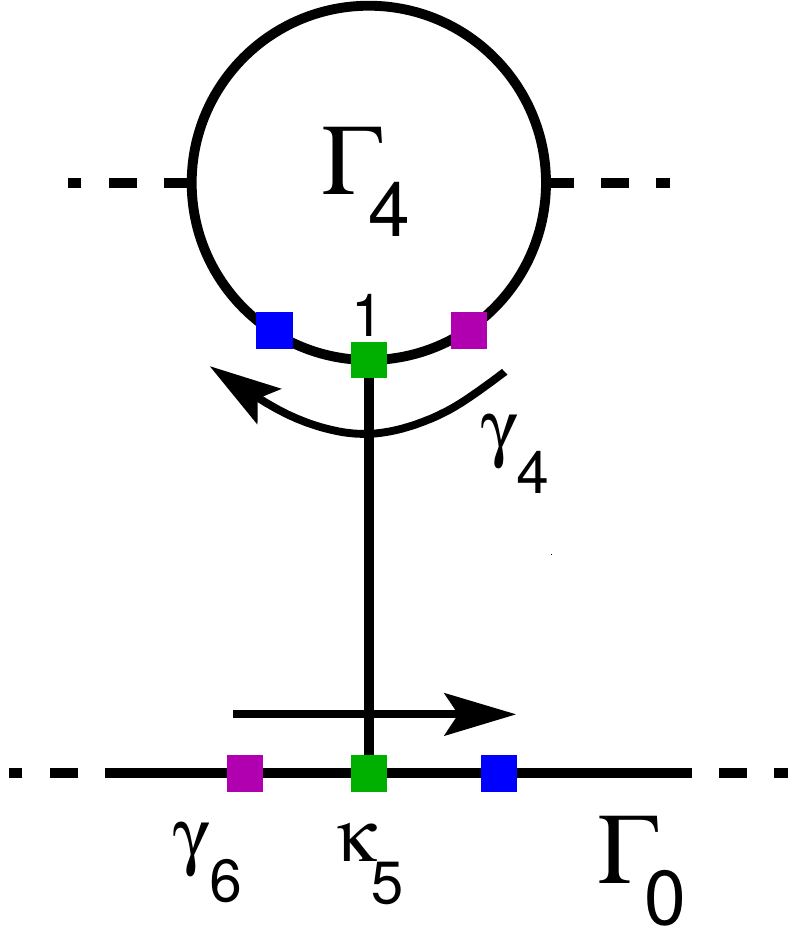}\hspace{0.5cm}\includegraphics[height=5cm]{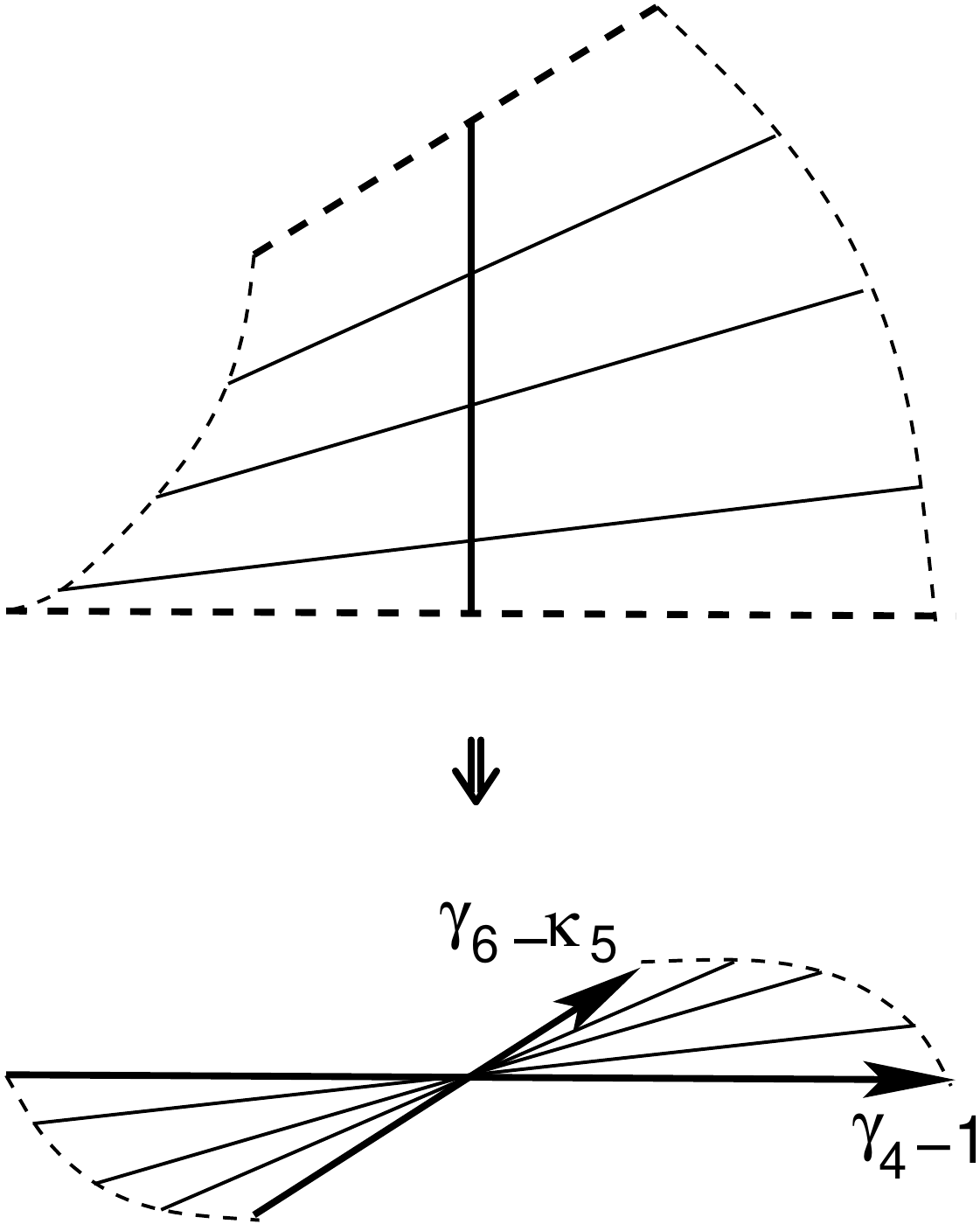}\hspace{0.5cm}\includegraphics[height=5cm]{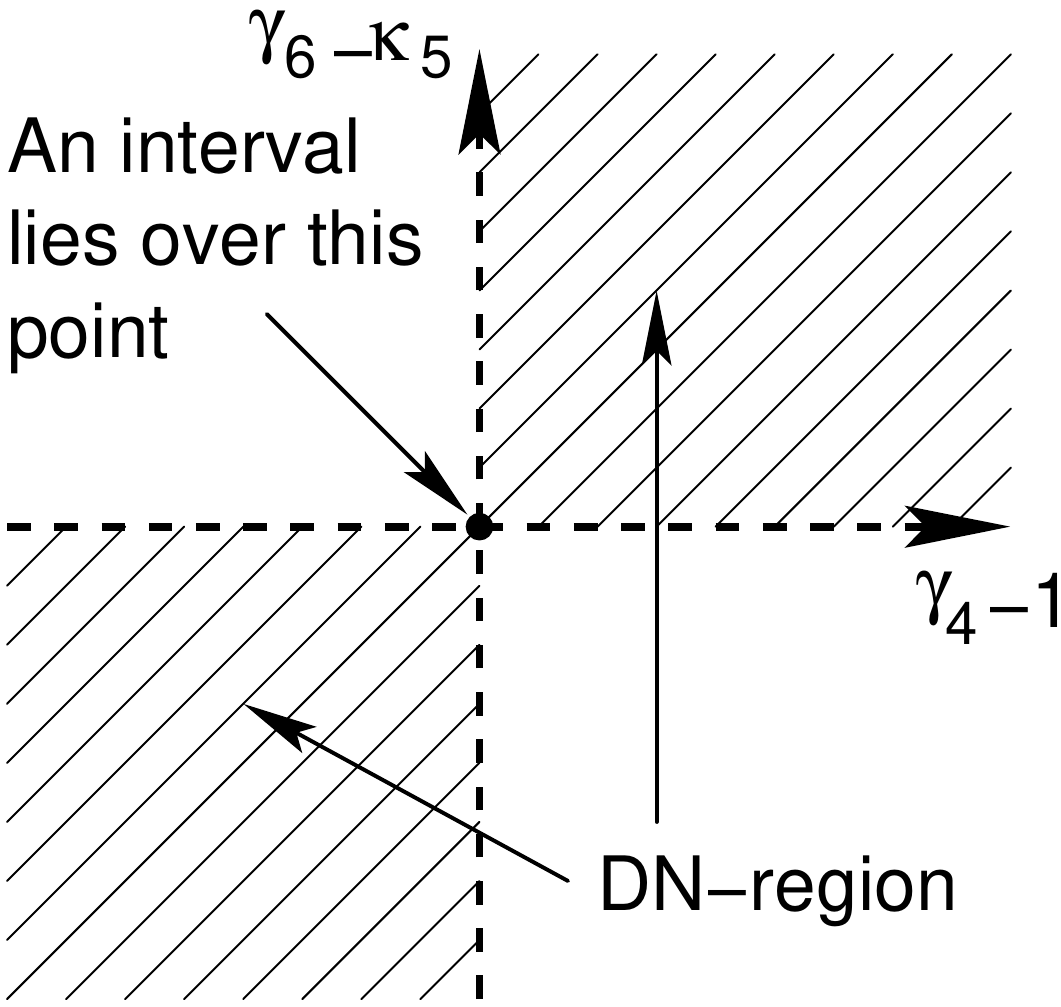}
    \caption{\small{\sl On the left: the divisor points $\gamma_4$ and $\gamma_6$ simultaneously pass through the double point; in the middle: the resolution of singularity for this pair of divisor points; on the right: the positions of divisor points $\gamma_4$, $\gamma_6$ satisfying the DN reality regularity condition.}}\label{fig:sing2}
  \end{figure}

Therefore, the Abel transform remains finite if and only if the two divisor points $\gamma_4$ and $\gamma_6$ simultaneously go to point $1$ in $\Gamma_4$ and point $\kappa_5$ in $\Gamma_0$ respectively, and the ratio
$$
\frac{\gamma_4-1}{\gamma_6-\kappa_5}
$$
has finite non-zero limit. In such case we have to apply the standard blow-up procedure near this point (see Figure~\ref{fig:sing2}):
\begin{equation}\label{eq:cartier2}
\RR^2 \rightarrow \RR^2\times \RR P^1, \ \ \big(\gamma_4,\gamma_6\big)\rightarrow \big(\gamma_4,\gamma_6,(\gamma_4-1):(\gamma_6-\kappa_5) \big).
\end{equation}

This situation exactly coincides with the situation discussed in Example~\ref{ex:3}. Therefore the Abel transform remains regular in terms of Cartier divisors.

\begin{definition}\label{def:case1}
We call \textbf{elementary non-smooth divisor configuration of the first type} the situation in which two divisor points of a DN divisor simultaneously approach a node of multiplicity two of the given $\mathtt{MM}$-curve in such a way that the ratio of the distances from this nodal point has a finite non-zero limit, and all other divisor points remain at finite distance from all other nodes of the curve.
\end{definition}

\begin{figure}[H]
  \centering\includegraphics[width=0.65\textwidth]{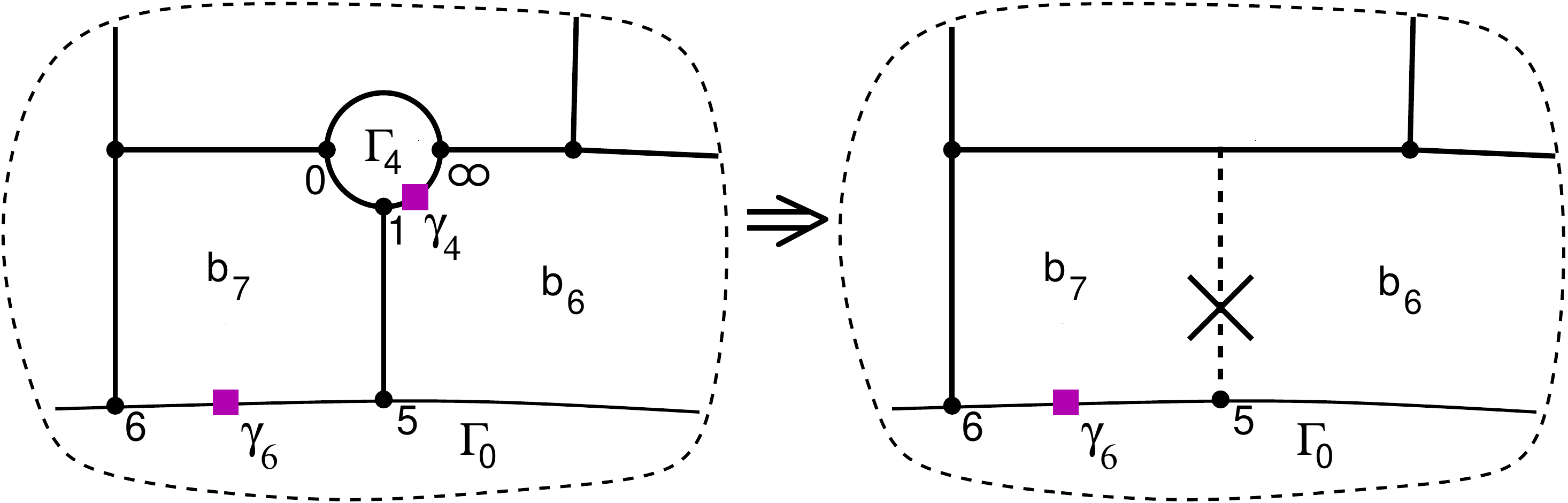}
      \caption{\small{\sl The ``ungluing'' of the double point connecting $\Gamma_4$ and $\Gamma_0$ when $\gamma_4\rightarrow 1$. After ungluing it, i.e. removing the edge connecting $\Gamma_0$ to $\Gamma_4$, $\Gamma_4$ becomes a bi-valent component without divisor, therefore it is removed.}}\label{fig:degen1}
  \end{figure}

\begin{remark}
  If all divisor points except $\gamma_4$ lie outside double points and $\gamma_4\rightarrow 1$, then $A_6({\mathcal D})\rightarrow\infty$ and $A_7({\mathcal D})\rightarrow\infty$, all other components of the Abel transform remain finite. The limit $\gamma_4=1$ means that we reached the boundary of the corresponding Jacobian component, and the boundary of the positroid cell which is of lower dimension. In terms of the spectral curve it means that we ``unglue'' the double point (see Figure~\ref{fig:degen1}). Such situation never happens in the KP dynamics for finite KP times.

  In contrast to the KP case, if one considers Melnikov-type equations associated with conformal transformations of surfaces in $\RR^4$, the ungluing of double points may happen in finite time, see \cite{GrinevichTaimanov2007,GrinevichTaimanov2008}.
\end{remark} 

\begin{figure}[H]
  \centering\includegraphics[height=4.cm]{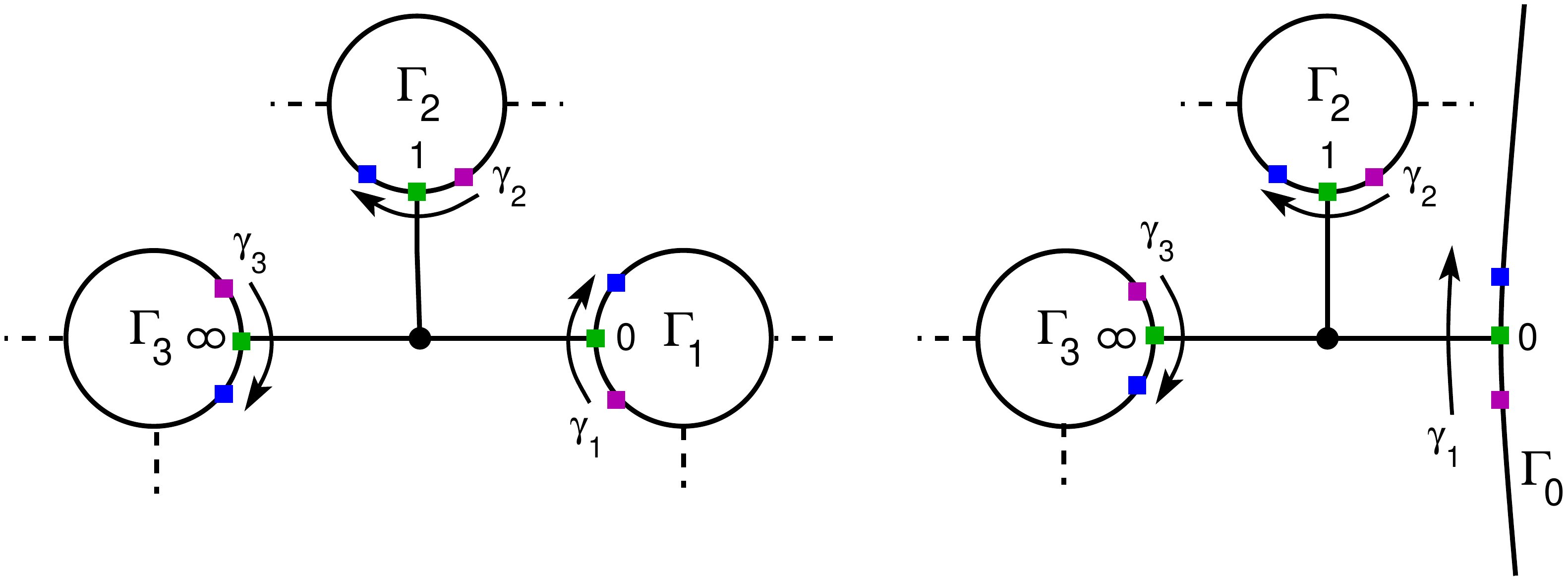}
    \caption{\small{\sl An elementary non-smooth divisor configuration of the second type: 3 divisor points  $\gamma_1$, $\gamma_2$ and $\gamma_3$ simultaneously pass through a triple point. On the left: all white components correspond to internal vertices. On the right: one of the components corresponds to the boundary of the disk.}}\label{fig:sing3}
  \end{figure}

Another elementary non-smooth divisor configuration corresponds to the situation illustrated in  Figure~\ref{fig:sing3}. Assume that we have a trivalent black vertex connected with 3 distinct trivalent white vertices, or with 2 distinct trivalent white vertices and the boundary of the disk. On the curve this configuration corresponds to a triple point. If these 3 divisor points simultaneously approach this triple point, and the ratios of distances
  $$
  \frac{\gamma_2-1}{\gamma_1}, \ \  \frac{1/\gamma_3}{\gamma_1},`
  $$
have finite non-zero limits, the Abel transform remain finite, and we have to apply the standard blow-up procedure near this point:
\begin{equation}\label{eq:cartier3}
\RR^3 \rightarrow \RR^3\times \RR P^2, \ \ \big(\gamma_2,1/\gamma_3,\gamma_1\big)\rightarrow \big(\gamma_2,1/\gamma_3,\gamma_1,(\gamma_2-1):(1/\gamma_3):\gamma_1 \big).
\end{equation}

For the configuration presented at Figure~\ref{fig:basis3} it may happen at the intersection point between the components $\Gamma_0$, $\Gamma_2$, $\Gamma_3$.

This case coincides with the one discussed in Example~\ref{ex:4}. Again, the Abel transform remains regular in terms of the Cartier divisor.

\begin{definition}\label{def:case2}
Consider a triple point belonging to the boundary of three ovals. Let three  points of a DN divisor on these ovals simultaneously approach this triple point in such a way that the ratio of the distances to this triple point has a finite non-zero limit, and all other divisor points remain at finite distance from all other nodes of the curve. We call such situation \textbf{elementary non-smooth divisor configuration of the second type.} 
\end{definition}  

\begin{remark}
In the example of Figure~\ref{fig:basis3} the two non-smooth divisor configurations described in Definitions~\ref{def:case1} and \ref{def:case2} may happen simultaneously. In general, it is natural to expect that several degenerations of this type may occur simultaneously.  
\end{remark}

\subsection{Dubrovin-Natanzon Cartier divisors}

It is natural to ask whether more complicated singular configurations of divisors may take place. It is clear that, if we do not impose the reality and regularity conditions, more complicated degenerations can happen. Before entering in the classification of singular Weil divisors configurations compatible with the Dubrovin-Natanzon condition, let us reformulate Dubrovin-Natanzon constraints in terms of Cartier divisors.

\begin{definition}\label{def:DNC}
Let $\Gamma$ be a reduced $\mathtt{MM}$-curve whose dual graph is the Le-graph of the given positroid cell $S\subset\GTNN$.  We call a Cartier divisor Dubrovin-Natanzon, if it satisfies the following constraints:
\begin{enumerate}
\item The degree of the divisor is equal to the arithmetic genus $g$ of the curve and coincides with the dimension of $S$.
\item The divisor is effective, i.e. all functions defining the divisors are holomorphic on the open sets where they are defined.
\item The functions defining the divisor do not vanish identically at any component, and at each white component the divisor has degree 1.
\item The divisor is real, i.e. the functions defining the divisors are real on the real part of the curve.
\item The divisor is odd at each finite oval (see Definition~\ref{def:parity}).  
\end{enumerate}
\end{definition}  

It is easy to check that such divisor has degree $k$ on the boundary component, since the number of white components is $g-k$. From DN condition it also follows that all zeroes of functions defining the divisor at $\Gamma_0$ are simple.

\begin{figure}[H]
  \centering\includegraphics[width=8cm]{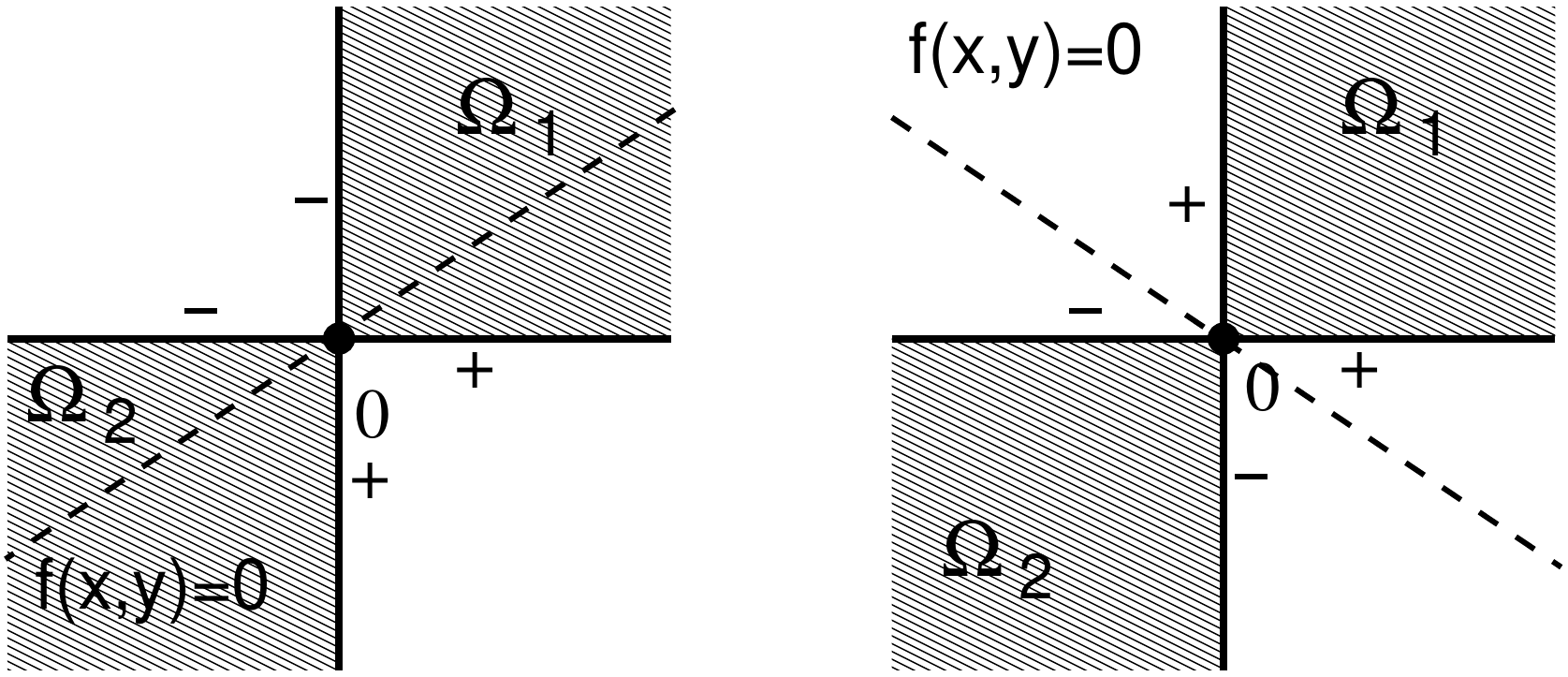}
  \caption{\small{\sl The parity of Cartier divisor point $2 P_0$, $P_0=0$ on ovals. We assume that the Weil divisor $2P_0$ corresponds to the Cartier divisor defined by a linear function $f(x,y)$.  On the left: the divisor point $2 P_0$ is odd on both ovals $\Omega_1$ and $\Omega_2$. On the right:  the divisor point $2 P_0$ is even on both ovals.}}\label{fig:parity}
  \end{figure}

\begin{figure}[H]
  \centering\includegraphics[width=11cm]{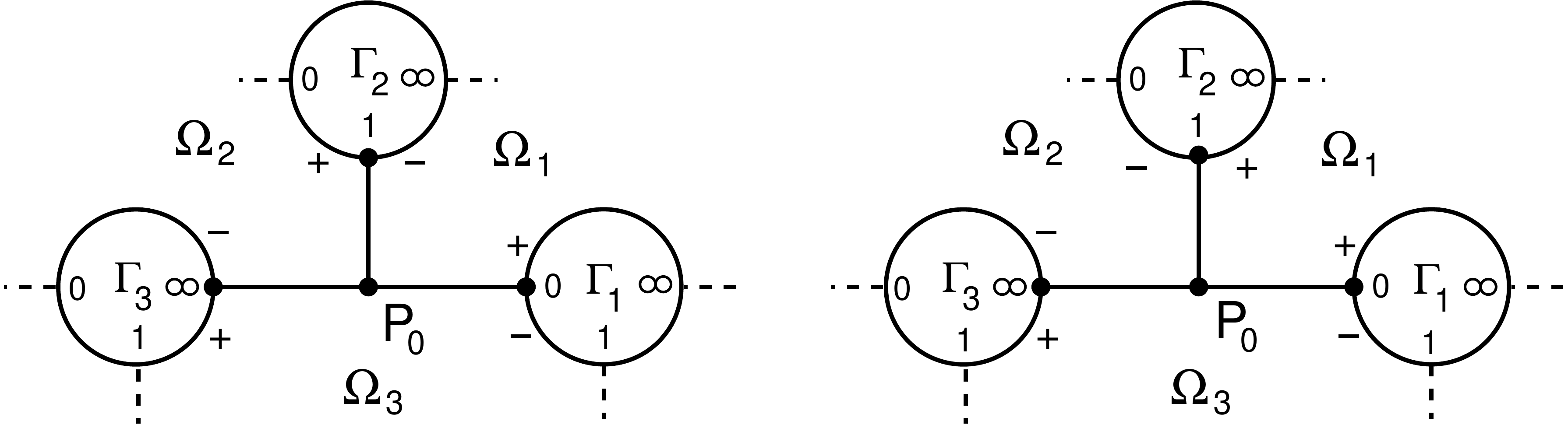}
  \caption{\small{\sl The parity on ovals of Cartier divisor at a triple point $P_0$. The functions defining the divisor near $P_0$ are $f_1=\lambda_1 z$ on $\Gamma_1$,  $f_2=\lambda_2 (z-1)$ on $\Gamma_2$,  $f_3=\lambda_3 /z$ on $\Gamma_3$.  
      On the left: $\lambda_1<0$, $\lambda_2<0$, $\lambda_3>0$, the divisor point $3 P_0$ is odd on all ovals $\Omega_j$.
      On the right: $\lambda_1<0$, $\lambda_2>0$, $\lambda_3>0$, the divisor point $3 P_0$ is even on $\Omega_1$, $\Omega_2$ and odd on  $\Omega_3$.}}\label{fig:parity3}
  \end{figure}

\begin{definition}\label{def:parity}
Let us define the parity of a real divisor on an oval. The parity of the divisor on an oval is the sum of the parities of all zeroes of the functions defining the divisor at this oval.

\begin{enumerate}
\item If a zero lies outside a nodal point, the parity of this zero coincides with the parity of its order.
\item If a zero lies in a multiple point, consider the part of the oval near this point. In a sufficiently small neighborhood of the multiple point the functions defining the divisor are different from zero outside this point. We say that the zero is even if the signs coincide on both sides of the interval and odd otherwise, see Figures~\ref{fig:parity} and \ref{fig:parity3}. 
\end{enumerate}
\end{definition}

Let us remark that if a DN Weil divisor is smooth, i.e. has empty intersection with the set of multiple points of $\Gamma$, then it is a DN Cartier divisor in the sense of Definition~\ref{def:DNC}.

\subsection{Krichever KP wave function from Cartier divisors}
\label{sec:3.5}

Assume now that $\Gamma$ is the reduced $\mathtt{MM}$-curve from Section~\ref{sec:reduced_curves},  ${\cal D}_C$ is a Dubrovin-Natanzon Cartier divisor on $\Gamma$, ${\cal D}_W$ is the corresponding Weil divisor. In order to include the case of singular divisors, we modify the definition of KP wave function from Sections~\ref{sec:2.2} and \ref{sec:2.4} in the following way:

The function $\hat\psi(\gamma,x,y,t)$, $\gamma\in\Gamma$ has the following properties for fixed $(x,y,t)$:
\begin{enumerate}
\item $\hat\psi(\gamma,x,y,t)$ is meromorphic in $\gamma$ on $\Gamma_0\backslash P_0$ and all $\Gamma_j$, $j>0$. 
\item $\hat\psi(\gamma,x,y,t)$ is holomorphic in $\gamma$ outside the point $P_0$ and the support of ${\cal D}_W$.
\item If $P_j$ is a multiple point of $\Gamma$ obtained by gluing the points $Q_1$, \ldots, $Q_{l+2}$, and $P_j$ contains no points of ${\cal D}_W$, than $\hat\psi(\gamma,x,y,t)$ is holomorphic at these points and
  $$
 \hat\psi(Q_1,x,y,t)= \hat\psi(Q_2,x,y,t) = \ldots =  \hat\psi(Q_{l+2},x,y,t).
  $$
\item If $\gamma_j$ belongs to ${\cal D}_W$, and $\gamma_j$ lies outside multiple points,  $\hat\psi(\gamma,x,y,t)$ may have at most first-order pole at  $\gamma_j$.\label{def:item:4}
\item If $P_j$ is a multiple point of $\Gamma$ obtained by gluing the points $Q_1$, \ldots, $Q_{l+2}$, and $P_j$ belongs the support of ${\cal D}_C$, due to the DN conditions the multiplicity of the divisor point $P_j$ is equal to $l+2$. This divisor point is locally defined by a collection of functions
  $$
f_j(\gamma) = \lambda_j (\gamma-Q_j)+ O((\gamma-Q_j)^2), \ \ j=1,\ldots,l+2.
$$
Then  $\hat\psi(\gamma,x,y,t)$ may have first-order poles at $Q_1$, \ldots, $Q_{l+2}$, and 
 \begin{equation}\label{eq:mult_div}
  f_1(\gamma)\hat\psi(\gamma,x,y,t)\bigg|_{\gamma=Q_1}= f_2(\gamma)\hat\psi(\gamma,x,y,t)\bigg|_{\gamma=Q_2}=\ldots =f_{l+2}(\gamma)\hat\psi(\gamma,x,y,t)\bigg|_{\gamma=Q_{l+2}}.
 \end{equation}\label{def:item:5}
 Equivalently, this collection $\{f_1(\gamma),\ldots,f_{l+2}(\gamma) \}$ is a regular function $f(\gamma)$ defined near the point $P_j$ and condition~(\ref{eq:mult_div}) means that the product $f(\gamma)\hat\psi(\gamma,x,y,t)$ is regular at $P_j$.

\item \label{item:6} At the point $P_0$ the function $\hat\psi(\gamma,x,y,t)$ has an essential singularity such that
$$
 \Xi(\gamma,x,y,t) = \hat\psi(\gamma,x,y,t)\cdot \exp\left(-\zeta (x-x_0) - \zeta^2 (y-y_0) -  \zeta^3 (t-t_0) \right),
$$
is holomorphic in $\gamma$ at the point $P_0$, and
$$
\Xi(\gamma,x,y,t) = 1 + o(1) \ \ \mbox{as} \ \ \gamma\rightarrow P_0.
$$
\end{enumerate}

Let us remark that the condition~\ref{item:6} assumes that that $P_0$ lies outside the support of ${\cal D}_C$, and for DN divisors it is always true. 

  In the next section we prove that, if the soliton data belong to a totally positive irreducible Schubert cell and the reduced ${\mathtt{MM}}$-curve is dual to its Le-graph, then the divisors of zeroes of the KP wave function are well-defined DN Cartier divisors for all times. In particular, the only singular divisor configurations which may occur are the ones introduced in Definitions~\ref{def:case1} and \ref{def:case2}. Finally we show that DN Cartier divisors fully parametrize this positroid cell.

\section{Main results} 

In this section we fix a totally-positive Schubert cell $S$ in $\GTNN$ (see Appendix for definitions) and denote by $\Gamma$ the reduced ${\mathtt{MM}}$-curve dual to the Le-graph associated to $S$ as in Section~\ref{sec:reduced_curves}. We also use the notation: $\vec t_0=(x_0,y_0,t_0)$,  $\vec t_1=(x_1,y_1,t_1)$.

Our aim is to show that $S$ can be interpreted as a real component of the Jacobian of $\Gamma$.

In \cite{AG4} we proved that for any point  $[A]\in S$ one can choose a normalization time such that the corresponding Dubrovin-Natanzon divisor of poles is smooth on $\Gamma$ (in \cite{AG4} the divisors were treated in Weil sense). In the next Theorem  we prove that the divisor of zeroes of the wave function is well-defined in the Cartier sense because the wave function does not vanish identically at any component. 

\begin{theorem}\label{thm:1}

Let $[A]$ be point in $S$.  Consider the following algebro-geometrical Krichever data associated to this point:
\begin{itemize}
\item  The reduced ${\mathtt{MM}}$-curve $\Gamma$.
\item A generic normalization point  $\vec t_0$ such that the corresponding divisor of poles $\mathcal D$ is smooth and satisfies the  Dubrovin-Natanzon (DN) condition formulated in Section~\ref{sec:2.4}.
\end{itemize}
Then for any collection of times $\vec t$
\begin{enumerate}
\item  The Krichever wave function does not vanish identically in the spectral parameter at any $\mathbb{CP}^1$ component of $\Gamma$.
\item  The divisor of zeroes is well-defined for all $\vec t$ as a Cartier divisor and fulfills Definition~\ref{def:DNC}.
\item  During the evolution of the KP hierarchy flows, the Dubrovin-Natanzon Weil divisor of zeroes has singularities only of the two types described in Definitions~\ref{def:case1} and \ref{def:case2}.
\end{enumerate}
\end{theorem}

Theorem~\ref{thm:1} is proven in the next section. 

\begin{remark}
  The proof of this Theorem relies heavily on total non-negativity. For divisors without DN condition it is very likely that the Krichever wave function may vanish identically at white components for some times and for these times the Cartier divisor of zeroes is not defined any more.
\end{remark}

\begin{theorem}\label{thm:2}
Let us fix the reference time $\vec t_0$. Then, any given point $[A]\in S$ is associated to a Dubrovin-Natanzon Cartier divisor on the curve $\Gamma$. Moreover, the only possible singular configurations of divisors are those described in Theorem~\ref{thm:1}.
\end{theorem}  

\begin{proof}
Indeed, thanks to Theorem~\ref{thm:1}, to construct the DN Cartier divisor corresponding to a given point $[A]\in S$, we can use the following procedure. From $\cite{AG4}$ it follows that the exists a normalization time $\vec t_1$ such that the corresponding divisor ${\cal D}(\vec t_1)$ is a smooth Dubrovin-Natanzon Weil divisor, and we can construct the corresponding Krichever wave function. By Theorem~\ref{thm:1} the divisor of zeroes ${\cal D}(\vec t_0)$ calculated at time $\vec t_0$ is a well-defined Dubrovin-Natanzon Cartier divisor. By Lemma~\ref{lem:norm_time} divisor ${\cal D}(\vec t_0)$  corresponds to the point $[A]\in S$ for normalization time  $\vec t_0$, and does not depend on the choice of $\vec t_1$. 

Therefore, we associate a unique  Dubrovin-Natanzon Cartier divisor to each point$[A]\in S$, for generic points these divisors are smooth Weil divisors constructed in $\cite{AG4}$, and the only possible singularities are those described in Theorem~\ref{thm:1}.
\end{proof}

\begin{theorem}\label{thm:3}
  Let the following spectral data be given:
\begin{enumerate}
\item The spectral curve $\Gamma$.  
\item A marked point $P_0\in\Gamma_0$  in the infinite oval of $\Gamma$,
\item The coordinates $\kappa_1<\kappa_2<\cdots < \kappa_n$ of the double points of $\Gamma_0$ (all of them are real and correspond to the chosen normalization of the boundary vertices of the Le-graph).
\item A Dubrovin-Natanzon Cartier divisor ${\cal D}_C$ on $\Gamma$ as in Definition~\ref{def:DNC}.
\end{enumerate}
  
Then
\begin{enumerate}
  
\item There exists a unique Krichever wave function as defined in Section~\ref{sec:3.5}.
\item The corresponding soliton data $[A]$ belongs to $S$ and its representative matrix may be computed from the reduced canonically oriented Le-network with weights as in Lemmas~\ref{lem:weights} and \ref{lem:pos_weights1} for smooth divisors and as in Lemma~\ref{lem:pos_weights2} if the divisor is singular.
\end{enumerate}
\end{theorem}  

We prove Theorem~\ref{thm:3} in Section~\ref{sec:6}. 

From the proof of Theorem~\ref{thm:3}, it easily follows that:

\begin{corollary}
The map from the Dubrovin-Natanzon Cartier divisors to the positroid cell defined in Theorem~\ref{thm:3} is continuous.
\end{corollary}

Next Theorem is the main result of our paper.

\begin{theorem}\label{thm:4}
The Dubrovin-Natanzon Cartier divisors on the curve $\Gamma$ fully parametrize the positroid cell $S$. Therefore $S$ can be naturally interpreted as a real component of the Jacobian of $\Gamma$.
\end{theorem}

\begin{proof}
First, let us fix a normalization time $\vec t_0$, and the phases $\kappa_1$, \ldots, $\kappa_n$ on $\Gamma_0$. In Theorem~\ref{thm:2} we associate a Dubrovin-Natanzon Cartier divisor to any point in $S$, and the map is regular.

Vice versa, by Theorem~\ref{thm:3},  to any Dubrovin-Natanzon Cartier divisor we associate a unique point in $S$, and the map is regular. This completes the proof.
\end{proof}

\section{Proof of Theorem~\ref{thm:1}}

The idea of the proof of the first part of this Theorem is the following. Let us assume that for some $\vec t_1$ the Krichever wave function $\hat\psi(\gamma,\vec t_1)$ vanishes identically in $\gamma$ at some $\Gamma_j$, $j>0$. Then we prove a sequence of lemmas which provide constraints on the possible configurations of these components. Finally we show that these constraints are incompatible with the DN conditions on the divisor. 

Let us introduce some useful terminology. If the wave function $\hat\psi(\gamma,\vec t_1)$ vanishes at some $\mathbb{CP}^1$ components, then we mark all the corresponding white vertices and all edges at which the wave function vanishes. In our Figures we use blue color to mark such vertices and edges.

Then the collections of marked vertices and edges have the following properties:

\begin{lemma}\label{lem:01}
\begin{enumerate}
\item If an edge at a black vertex is marked, then all edges connected to this vertex are also marked, and this black vertex is marked.
\item Let $\Gamma_j$ be a white vertex such that at least two edges at $\Gamma_j$ are marked. Then $\Gamma_j$ is also marked.
\end{enumerate}  
\end{lemma}

The proof is omitted because it directly follows from the properties of the wave function. Indeed, at each $\mathbb{CP}^1$ component corresponding to a white vertex, the wave function $\hat\psi(\gamma,\vec t)$ is a degree 1 meromorphic function of $\gamma$, and it either has exactly one zero, or is constant. Therefore, if it has zeroes at two distinct points of a white component, it necessarily vanishes identically at this component. 

\begin{figure}[H]
  \centering\includegraphics[width=0.9\textwidth]{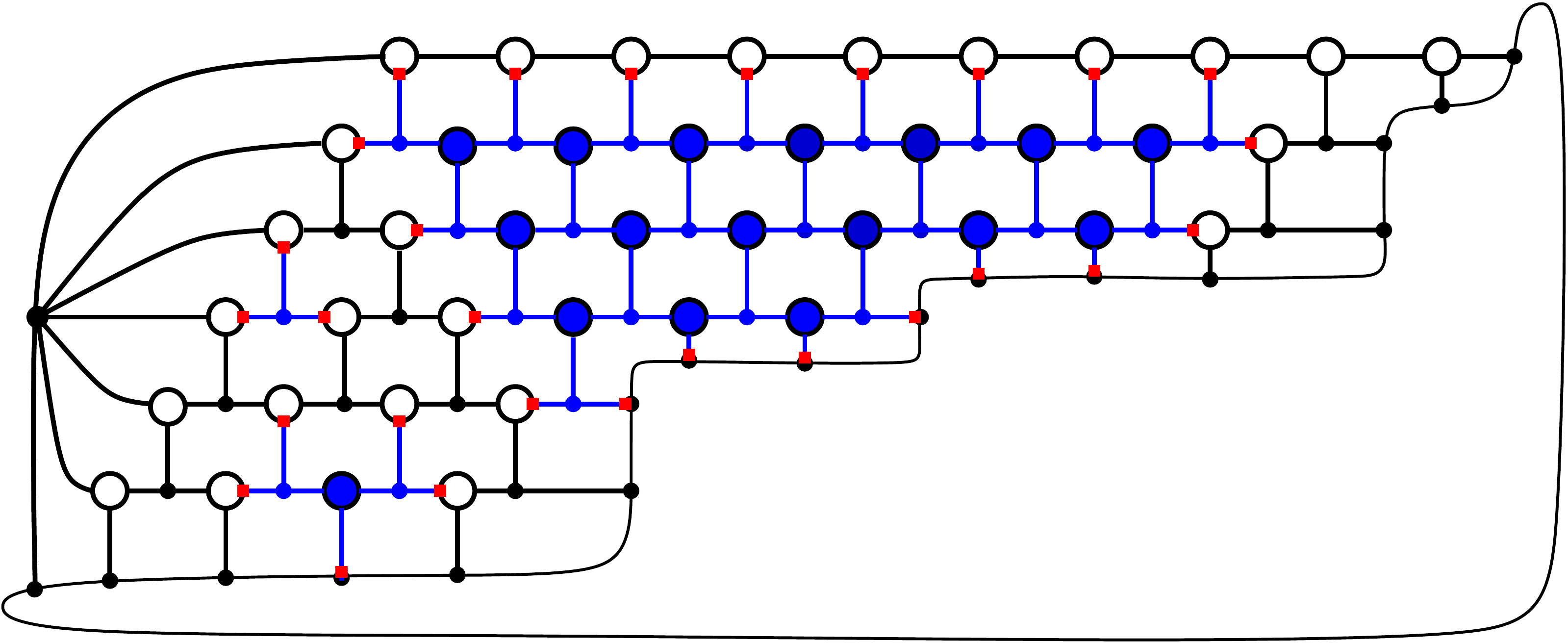}
      \caption{\small{\sl Example of the reduced Le-graph of a totally positive Schubert cell with marked vertices and edges. In the current Figure marked vertices and edges are colored blue. If an edge connects a marked component to a non-marked one, we have a divisor point (marked red) at the non-marked end. }}\label{fig:marked1}
			
\end{figure}

\begin{corollary}
If an edge connects a marked component to a non-marked one, then the non-marked component either corresponds to a white vertex, or to $\Gamma_0$. In both cases, there is a divisor point at the non-marked end which we color red in the Figures.  
\end{corollary}
An example is presented in Figure~\ref{fig:marked1}.

\begin{figure}[H]
  \centering\includegraphics[width=0.85\textwidth]{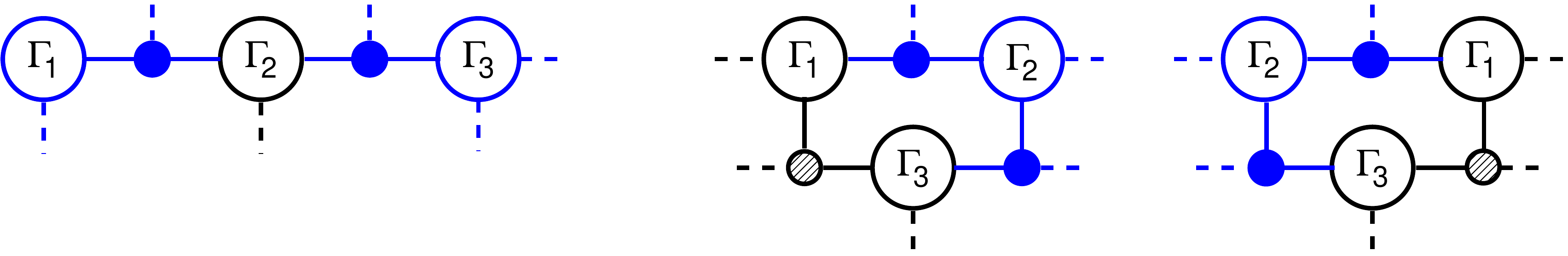}
      \caption{\small{\sl This Figure illustrates Lemma~\ref{lem:2}.}}\label{fig:lemma2}
  \end{figure}

\begin{lemma}\label{lem:2}
 \begin{enumerate} 
 \item  Let the white vertices  $\Gamma_1$,  $\Gamma_2$ and  $\Gamma_3$ lie in the same row and be consecutive as in  Figure~\ref{fig:lemma2} left. If both $\Gamma_1$ and  $\Gamma_3$ are marked, then $\Gamma_2$ is also marked.
 \item Let $\Gamma_2$ be a marked white vertex, $\Gamma_1$ lie in the same row and be the closest to $\Gamma_2$ either at its left or at its right, and let $\Gamma_3$ lie in the next row downwards between $\Gamma_1$ and $\Gamma_2$ (see Figure~\ref{fig:lemma2} middle and right). Then either both $\Gamma_1$ and $\Gamma_3$ are marked, or both of them are not marked.  
\end{enumerate}

\end{lemma}
  
The proof immediately follows from Lemma~\ref{lem:01}.

\begin{lemma}\label{lem:3}
The white vertices in the upper row of the diagram and the white vertices located at the left end of each row are never marked.
\end{lemma}
\textbf{Proof of Lemma~\ref{lem:3}.}

If the divisor satisfies the DN conditions, then the wave function in Krichever normalization $\hat\psi(\gamma,\vec t)$ is real and strictly positive at the infinite oval. The infinite oval passes through all components  corresponding to the white vertices in the upper row; therefore these vertices cannot be marked. If a white vertex is located at the left end of a given row, then the corresponding component is connected to a black vertex intersecting the infinite oval, and this white vertex cannot be marked. \qed
\medskip

\begin{figure}[H]
  \centering\includegraphics[width=0.75\textwidth]{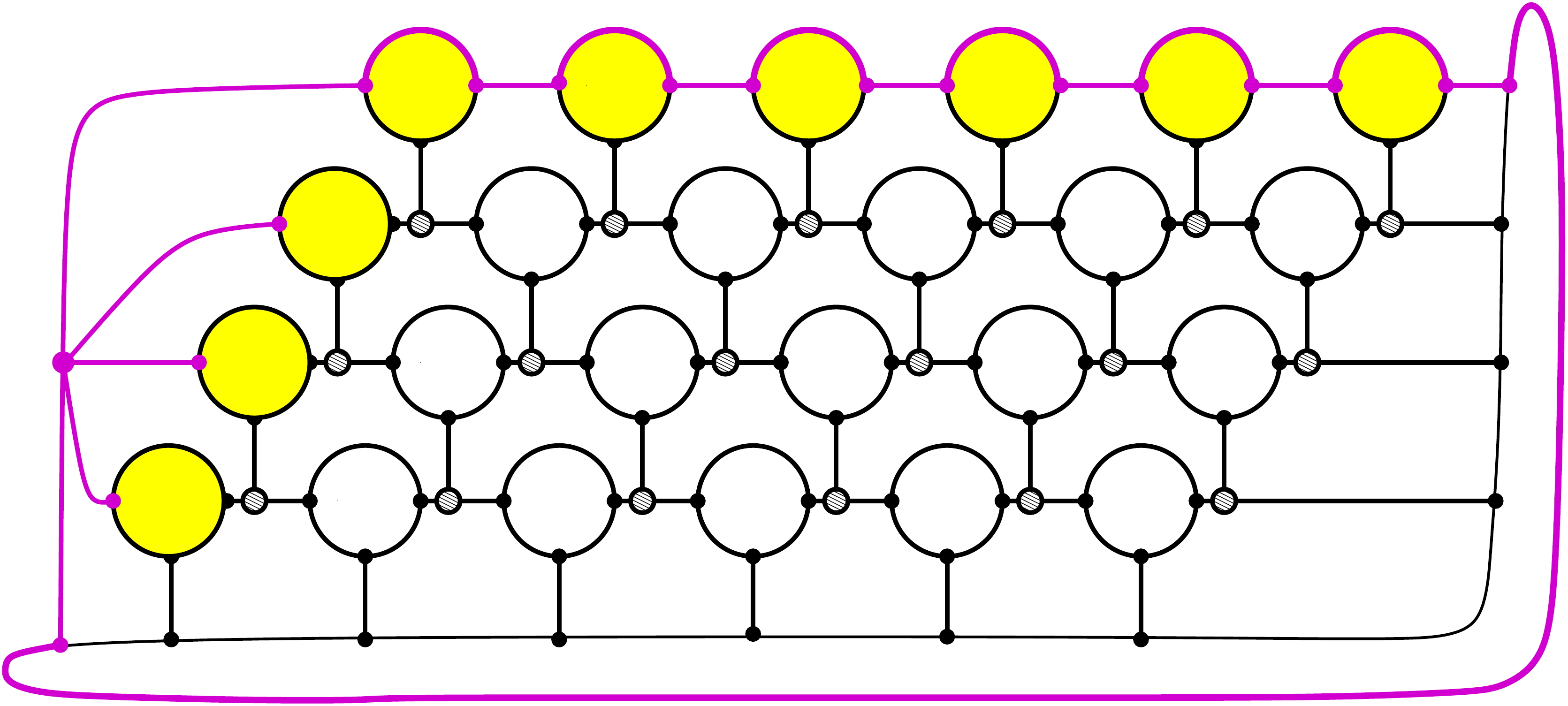}
      \caption{\small{\sl The infinite oval is painted magenta. At the points marked magenta the wave function is strictly positive for any $\vec t$. At the components marked yellow the wave function never vanishes identically in the spectral parameter.}}\label{fig:lemma3.2}
\end{figure}

We illustrate Lemma~\ref{lem:3} in Figure~\ref{fig:lemma3.2}.

\begin{definition}
Let us call two vertices of the graph \textbf{directly connected} if the exists an edge connecting them. Similarly, a vertex is directly connected to the boundary if there exists an edge connecting such vertex to the boundary.
\end{definition}

\begin{figure}[H]
  \centering\includegraphics[width=0.75\textwidth]{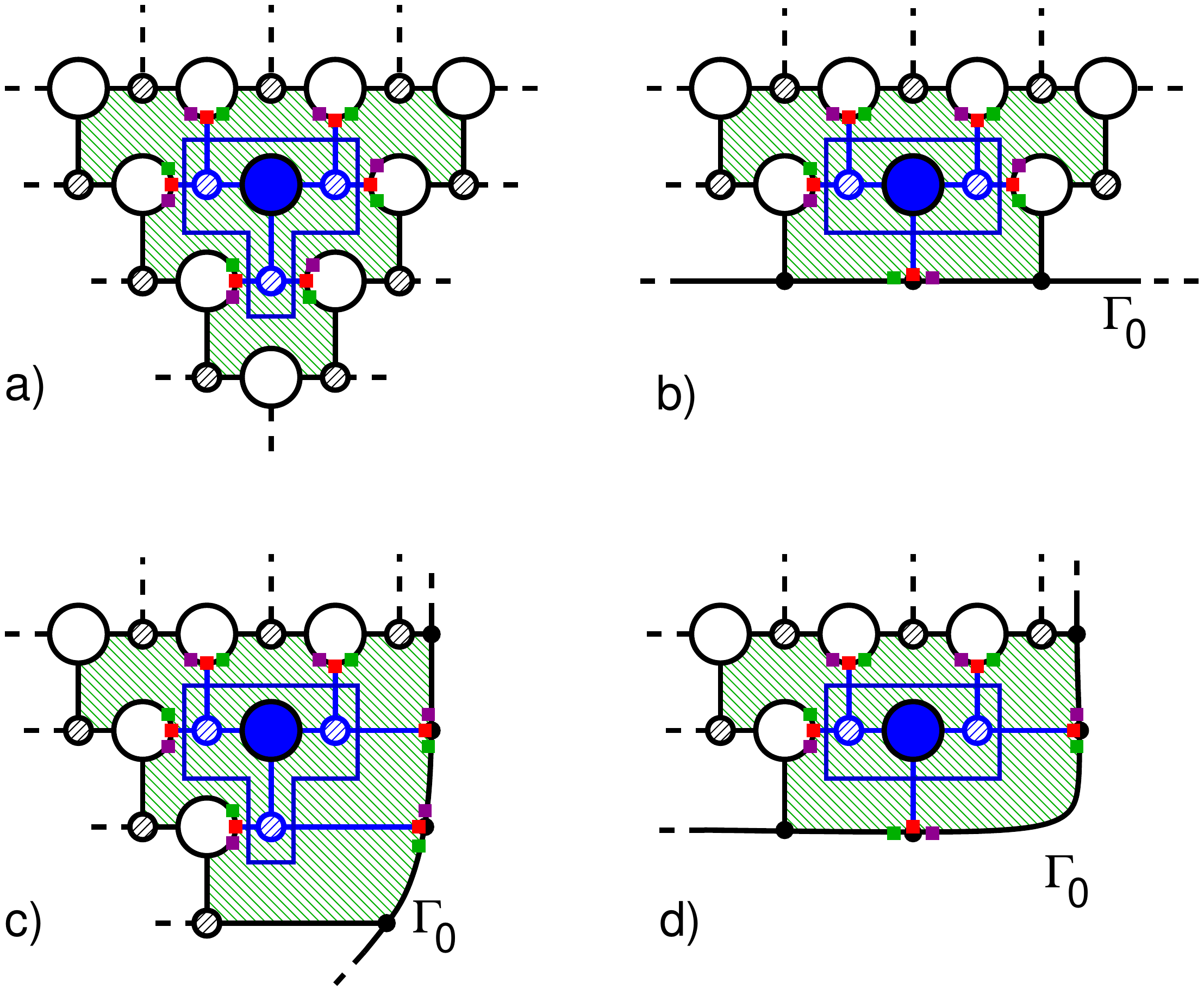}
    \caption{\small{\sl This Figure illustrate the case of one isolated marked white vertex. This vertex and the black vertices adjacent to this white one are marked blue and enclosed by a blue polygon. The faces surrounding this marked area are colored salad green.}}\label{fig:marked0}
\end{figure}

\textbf{Proof of Theorem~\ref{thm:1}.} Let us start from the simplest case of one isolated marked white vertex,
see Figure~\ref{fig:marked0}. By Lemma~\ref{lem:01} all black vertices directly connected to this white one are also marked. If an edge connects a marked black vertex either to a non-marked white one or to a boundary vertex, then this edge is also marked, and a divisor point is located at the end of the edge at the non-marked component. These divisor points are marked red in Figure~\ref{fig:marked0}. The faces surrounding the marked area are colored salad green. Each red divisor point is located at the edge between two such faces; therefore the total number of green faces is equal to the total number of red divisor points. 

For generic time no divisor point lies at the double points of $\Gamma$. A small generic perturbation of $\vec t_1$ will move these red divisor points inside green areas, and each green area contains exactly one red divisor point. Therefore only two possible shifts of divisor points are compatible with the DN conditions; one is colored green in Figure~\ref{fig:marked0}, the other one is colored magenta. As a consequence, no position of the divisor point at the marked white vertex is compatible with the  DN conditions. We come to a contradiction, and this ends the proof in this case.

Let us consider the general situation when more than one white vertex is marked in some connected region. We proceed inside the marked area from top to bottom, and we move within a row from left to right.
We also use the following notation. If  $\Gamma^{i}_{j-1}$, $\Gamma^{i}_j$ are two subsequent white vertices on row $i$, and there exists a white vertex in the next row $i+1$ lying between  $\Gamma^{i}_{j-1}$ and $\Gamma^{i}_j$, then we denote it by  $\Gamma^{i+1}_j$ (see, for example, Figure~\ref{fig:lemma4}).

Let us prove the following:

\begin{lemma}\label{lem:4}
  Let $\Gamma^{i}_1$,\ldots, $\Gamma^{i}_r$, be subsequent white vertices in the row $i$ such that:
\begin{enumerate}  
\item All of them are marked;
\item The white vertex $\Gamma^{i}_0$ preceding $\Gamma^{i}_1$ in this row is not marked;
\item Either the white vertex $\Gamma^{i}_{r+1}$ next to  $\Gamma^{i}_r$ in this row is not marked, or $\Gamma^{i}_r$ is the last white vertex in the row $i$.   
\end{enumerate}
Then:
\begin{enumerate}
\item If none of  $\Gamma^{i}_1$,\ldots, $\Gamma^{i}_r$ are directly connected to the boundary, then
\begin{enumerate} 
\item The white vertex  $\Gamma^{i+1}_1$ exists and is not marked.
\item The white vertices $\Gamma^{i+1}_2$, \ldots, $\Gamma^{i+1}_r$ exist and are marked.
\item If a white vertex $\Gamma^{i+1}_{r+1}$ lying next to $\Gamma^{i+1}_r$ exists, it is not marked.
\end{enumerate}
We illustrate this case in Figure~\ref{fig:lemma4}.
\begin{figure}[H]
  \centering\includegraphics[width=0.8\textwidth]{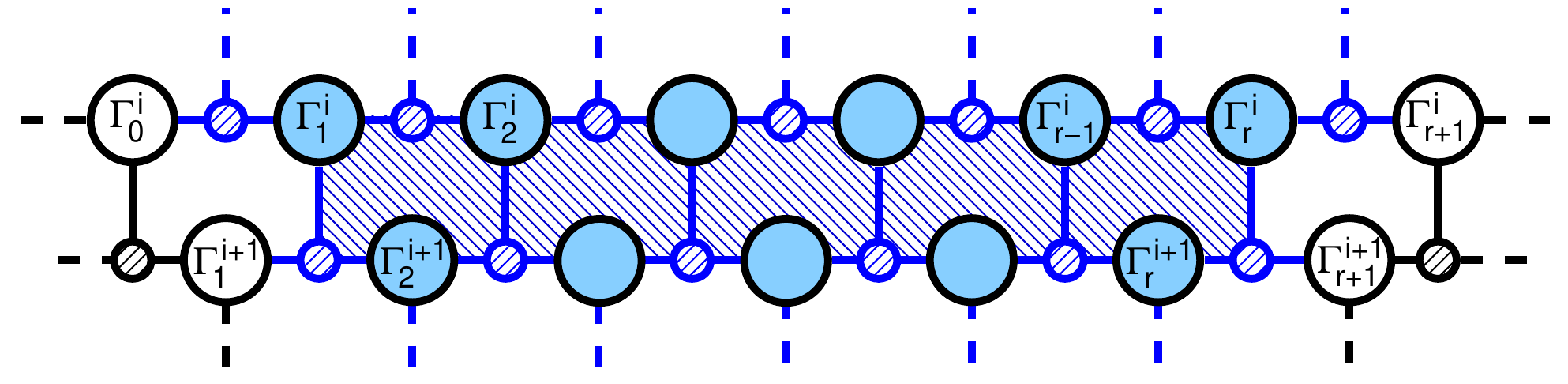}

  \vspace{5mm}
  
  \centering\includegraphics[width=0.8\textwidth]{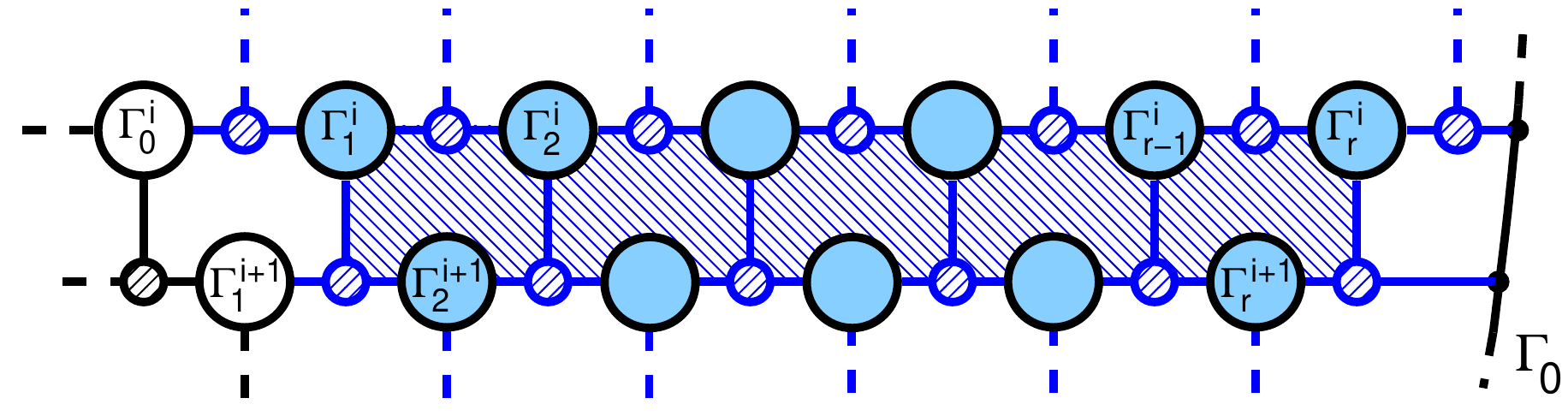}
  \caption{\small{\sl This Figure illustrates Lemma~\ref{lem:4} case 1. Internal faces of the marked area between the rows $i$ and $i+1$ are marked blue. Top: $\Gamma^{i}_r$ is not the last white vertex in the row $i$. Bottom: $\Gamma^{i}_r$ is the last white vertex in the row $i$.}} \label{fig:lemma4}
\end{figure}

\item Let $\Gamma^{i}_s$,  $s>2$, be the first vertex directly connected to the boundary, then:
\begin{enumerate}
\item The white vertex  $\Gamma^{i+1}_1$ exists and is not marked.   
\item All white vertices $\Gamma^{i}_j$, $j>s$ are also directly connected to the boundary.
\item The white vertices $\Gamma^{i+1}_2$, \ldots, $\Gamma^{i+1}_{s-1}$ exist and are marked.
\item The white vertex $\Gamma^{i+1}_{s-1}$ is the last one in the row $i+1$.
\end{enumerate}
We illustrate this case in Figure~\ref{fig:lemma4.3}.
\begin{figure}[H]
  \centering\includegraphics[width=0.9\textwidth]{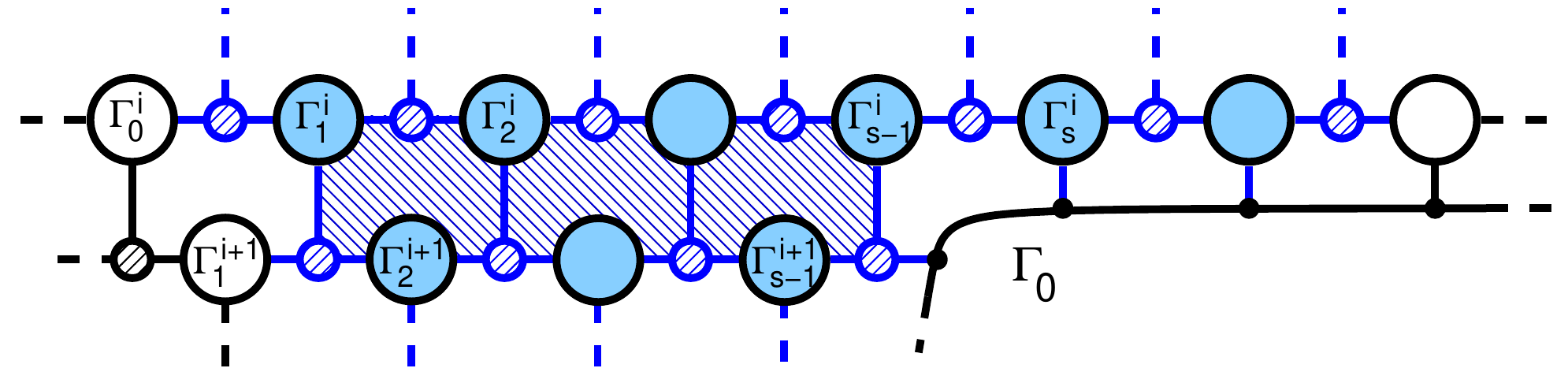}
    \caption{\small{\sl This Figure illustrates Lemma~\ref{lem:4} case 2. Internal faces of the marked area between the rows $i$ and $i+1$ are marked blue. $\Gamma^{i}_s$ is the first white vertex in the row $i$ directly connected to the boundary.}} \label{fig:lemma4.3}
\end{figure}  
\item  Let $\Gamma^{i}_2$ be the first vertex directly connected to the boundary, then:
\begin{enumerate}
\item The white vertex  $\Gamma^{i+1}_1$ exists and is not marked.   
\item  The white vertex  $\Gamma^{i+1}_1$ is the last one in the row $i+1$.
\end{enumerate}
We illustrate this case in Figure~\ref{fig:lemma4.5}.
\begin{figure}[H]
  \centering\includegraphics[width=0.9\textwidth]{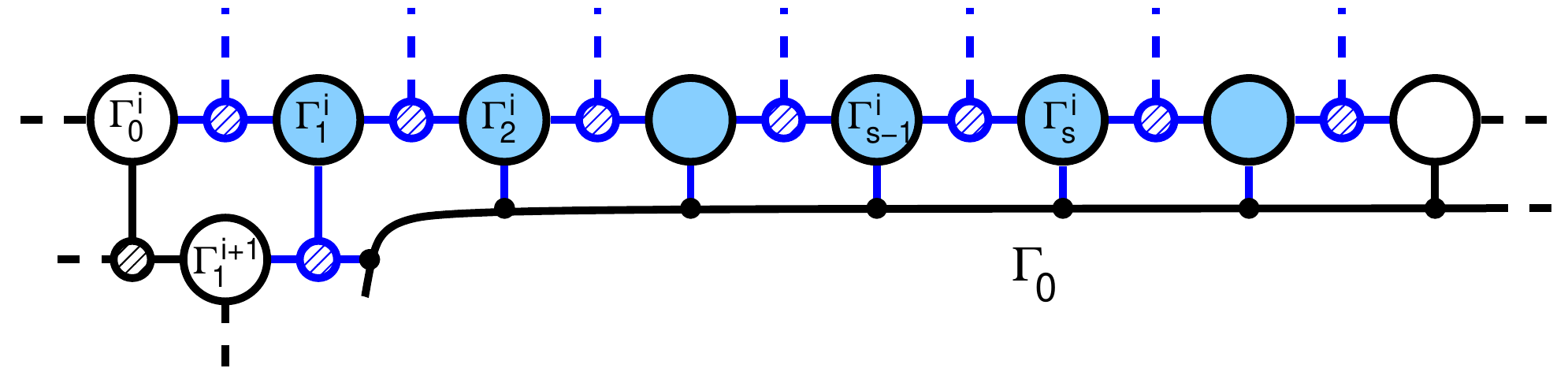} 
    \caption{\small{\sl This Figure illustrates Lemma~\ref{lem:4} case 3.  $\Gamma^{i}_2$ is the first white vertex in the row $i$ directly connected to the boundary. No faces between the rows $i$ and $i+1$ are marked.}}\label{fig:lemma4.5}
\end{figure}  
\item If $\Gamma^{i}_1$ is directly connected to boundary, there are no marked components in the row $i+1$ between $\Gamma^{i}_0$ and $\Gamma^i_r$. (see Figure~\ref{fig:lemma4.4}).
\begin{figure}[H]
  \centering\includegraphics[width=0.9\textwidth]{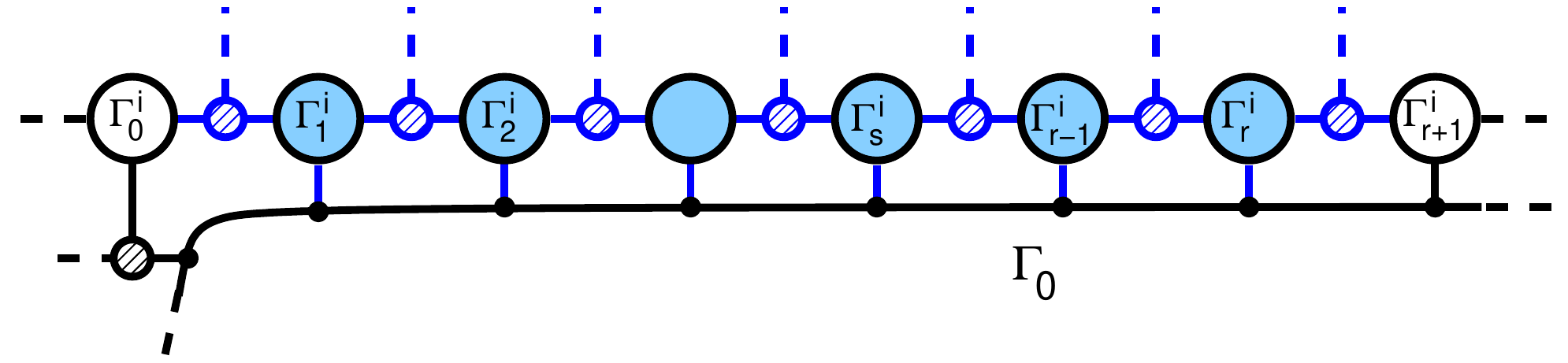}
    \caption{\small{\sl This Figure illustrates Lemma~\ref{lem:4} case 4.}}\label{fig:lemma4.4}
\end{figure}
\end{enumerate}
\end{lemma}

The proof of Lemma~\ref{lem:4} follows immediately from Lemmas~\ref{lem:01}-\ref{lem:3}.

From Lemma~\ref{lem:4} it immediately follows:

\begin{corollary}\label{cor:1}
\begin{enumerate}  
\item  The number of internal faces in a connected marked area between lines $i$ and $i+1$ is equal to the number of marked white vertices in row $i+1$ belonging to this marked area.
\item  The number of marked internal faces between the lines $i$ and $i+1$ is less than the number of marked white vertices in row $i$.
\item The number of internal faces of the marked areas is smaller than the number of marked white vertices.  
\end{enumerate}  
\end{corollary}  

Let us return to the proof of  Theorem~\ref{thm:1}. Consider a connected marked area. Again, the faces surrounding the marked area are colored salad green. Let $i_0$ be the first row of the blue area, and let a connected component of the marked part of this row contain white vertices $\Gamma^{i_0}_1$, $\Gamma^{i_0}_r$ enumerated from left to right. Then, using Lemma~\ref{lem:4} we can move down row by row, and at each step we obtain a connected component, which is located below the component from the previous row, and is shorter. Therefore each connected component of the marked set has a triangular shape. The number of salad green faces is equal to the number of edges, connecting the blue component either to $\Gamma_0$, or to non-marked white components. If an edge connects the marked area to boundary $\Gamma_0$, we have a divisor point at  $\Gamma_0$ at the corresponding boundary point; analogously, if an edge connects the marked area to a non-marked white component, we have a divisor point at this non-marked component at the corresponding point.

\begin{figure}[H]
  \centering\includegraphics[width=0.9\textwidth]{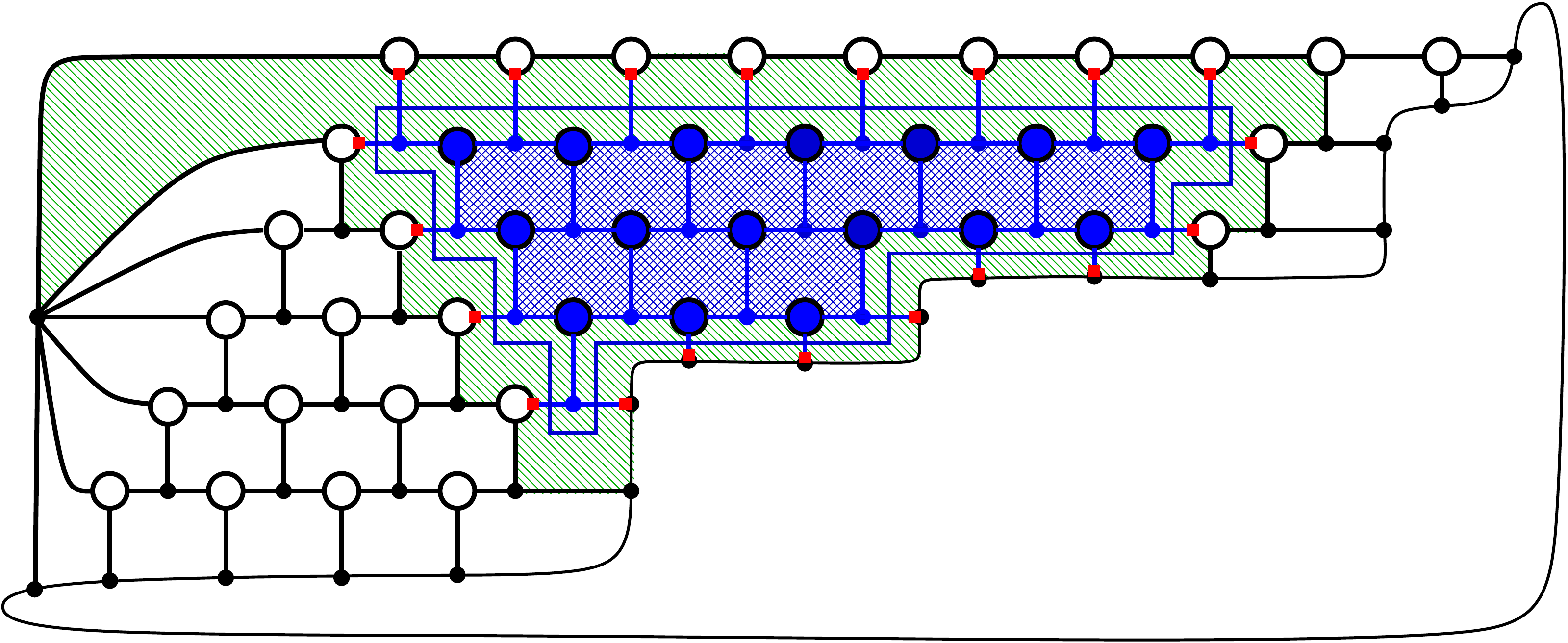}
    \caption{\small{\sl One connected component of the marked region is enclosed by a blue polygon. The faces surrounding this marked area are colored salad green. Internal faces of the marked area are colored blue.}} \label{fig:marked2}
\end{figure}

Applying a small shift of time, we obtain a generic configuration such that all divisor points corresponding to the boundary faces are located in the salad green area. By DN condition we have exactly one divisor point at each face, therefore all divisor points at the white vertices from the marked area shall correspond to the internal faces of the marked area. By Corollary~\ref{cor:1}, the number of these internal faces from this component is smaller then the number of marked white vertices from this component. We obtained a contradiction, therefore the first part is proved (see Figure~\ref{fig:marked2} for an example).

We proved that for all $\vec t$ the Krichever wave function never vanishes identically at any component of the reduced spectral curve, therefore we have a well-defined Cartier divisor $(\hat\psi(\gamma,\vec t))$ associated to it. By definition, the divisor of zeroes ${\cal D}(\vec t)$ can be written as: ${\cal D}(\vec t) = (\hat\psi(\gamma,\vec t))+{\cal D}(\vec t_0)$, where ${\cal D}(\vec t_0)$ is the divisor of poles, $\vec t_o$ is the normalization time. The divisor  $(\hat\psi(\gamma,\vec t))$ is even at all ovals as a divisor of a function, ${\cal D}(\vec t_0)$ is odd at all finite ovals due to the DN property, therefore ${\cal D}(\vec t)$ is also odd at all finite ovals. It ends the proof the second statement. 

To prove the last statement of Theorem~\ref{thm:1}, let us remark that singular divisors occur only if the KP wave function vanishes at a multiple point of $\Gamma$. Since  $S$ is a totally-positive Schubert cell, we have at most one multiple point of multiplicity greater than 3, and it intersects the infinite oval $\Omega_0$, therefore the wave function at this point never vanishes. All other multiple points have either degree 2 or 3.

If the KP wave function vanishes at a double point, obtained by gluing of a pair $x_1$, $x_2$, it has first-order poles in spectral parameter at $x_1$ and at $x_2$, and we have the first singular configuration. Equivalently, if  KP wave function vanishes at a triple point, obtained by gluing of  $x_1$, $x_2$ and $x_3$ it has first-order poles in spectral parameter at $x_1$, $x_2$ and $x_2$, and we have the second singular configuration.
This ends the proof of Theorem~\ref{thm:1}.

\section{Proof of Theorem~\ref{thm:3}}\label{sec:6}

To construct the Krichever wave function we use the following procedure.
\begin{enumerate}
\item To a given Dubrovin-Natanzon Cartier divisor we associate a collection of weights on the canonically oriented Le-graph such that the boundary measurement map matrix is totally non-negative.
\item Using the linear relations on the black and white vertices and at the boundary of the graph, we obtain the wave function as the solution of a full-rank linear system. 
\end{enumerate}  

Let us start by  associating a collection of Krichever's compatibility conditions on the components and on the multiple points to a given  Dunbrovin-Natanzon Cartier divisor.

To start with, let us assume that the divisor is smooth, i.e. this Cartier divisor can be treated as a Weil divisor.

The Krichever wave function on rational curves is rational at each irreducible rational components, associated to a white vertex, is a product of exponent to a rational function at the boundary component and is continuous at multiple points. Dubrovin-Natanzon conditions implies that at each white component  $\Gamma_{i,j}$ the rational function has degree 1 and has a unique real pole $\gamma_{i,j}$:
\begin{equation}
\label{eq:restriction1}
\hat\psi_{i,j}(\zeta, \vec t) =  c_{i,j}^{(0)}(\vec t) + \frac{c_{i,j}^{(1)}(\vec t)}{\zeta-\gamma_{i,j}}.
\end{equation}

\begin{figure}[H]
  \centering\includegraphics[width=0.2\textwidth]{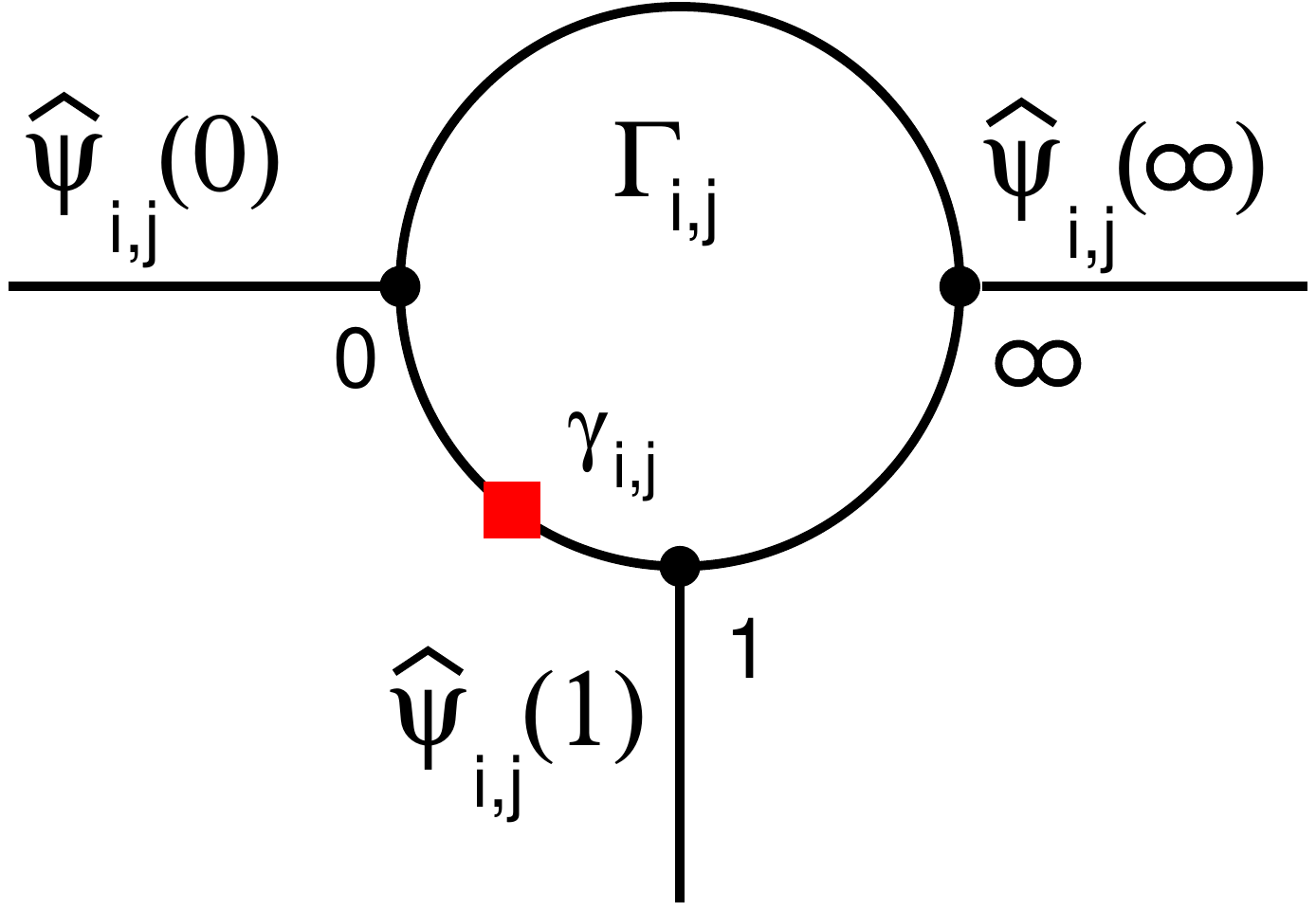}
      \caption{\small{\sl The restriction of the Krichever wave function to the component $\Gamma_{i,j}$.}}\label{fig:one_white}
\end{figure}

Therefore the wave function satisfies the following algebraic relations:

\begin{enumerate}
\item At each component $\Gamma_{i,j}$ associated to a white vertex we have:
\begin{equation}
\label{eq:relation_white}
\hat\psi_{i,j}(\infty, \vec t) =  \gamma_{i,j} \, \hat\psi_{i,j}(0, \vec t) + (1- \gamma_{i,j}) \hat\psi_{i,j}(1, \vec t).
\end{equation}
(This relation directly follows from (\ref{eq:restriction1})).
\item For each black vertex the values of the wave function at the corresponding points on the components connected to this vertex are the same.
\item If an edge connects either two white components or a white component to the boundary or two points of the boundary component, then  the values of the wave function at the corresponding points are the same.   
\end{enumerate}  

Let us illustrate conditions 2) and 3) by some examples:

\begin{enumerate}
\item If a black point is a triple point, the values of the wave function at the corresponding points of the white components are equal, see Figure~\ref{fig:one_black}:
\begin{equation}    
\label{eq:relation_black1}
\hat\psi_{i_{l},j_{s}}(\infty, \vec t) =  \hat\psi_{i_{l-1},j_{s}}(1, \vec t) = \hat\psi_{i_{l},j_{s-1}}(0, \vec t).
\end{equation}

\begin{figure}[H]
  \centering\includegraphics[width=0.5\textwidth]{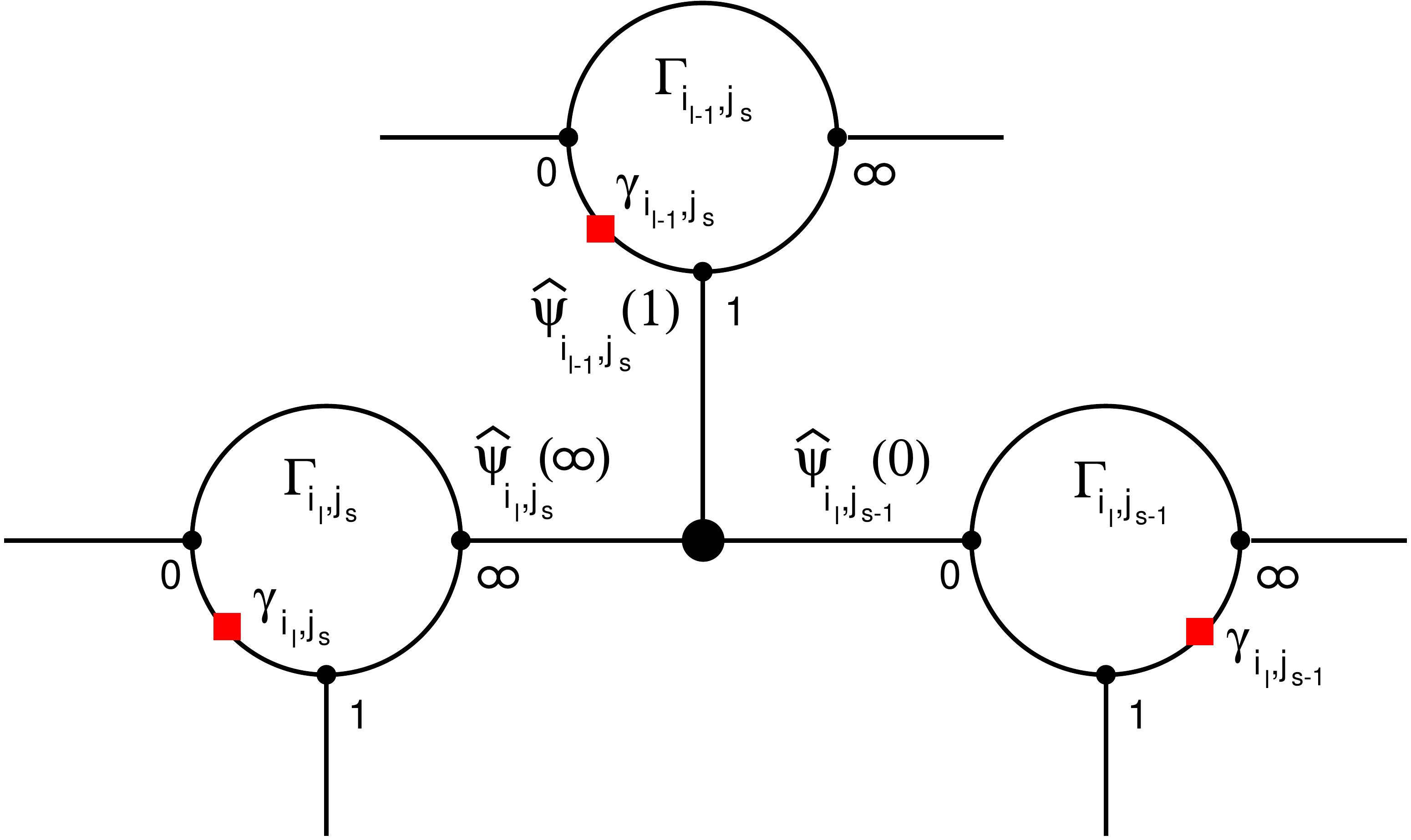}
      \caption{\small{\sl Intersection of three white components at one point. The values of the wave functions at the intersection points are equal, see Equation~\ref{eq:relation_black1}.}}\label{fig:one_black}
    \end{figure}

\item If two white vertices in the first row are connected by an edge, the values of the wave function at the corresponding points of the white components are equal, see Figure~\ref{fig:one_black2}:
\begin{equation}    
\label{eq:relation_black2}
\hat\psi_{1,j_{s}}(\infty, \vec t) =  \hat\psi_{1,j_{s-1}}(0, \vec t).
\end{equation}
\begin{figure}[H]
  \centering\includegraphics[width=0.5\textwidth]{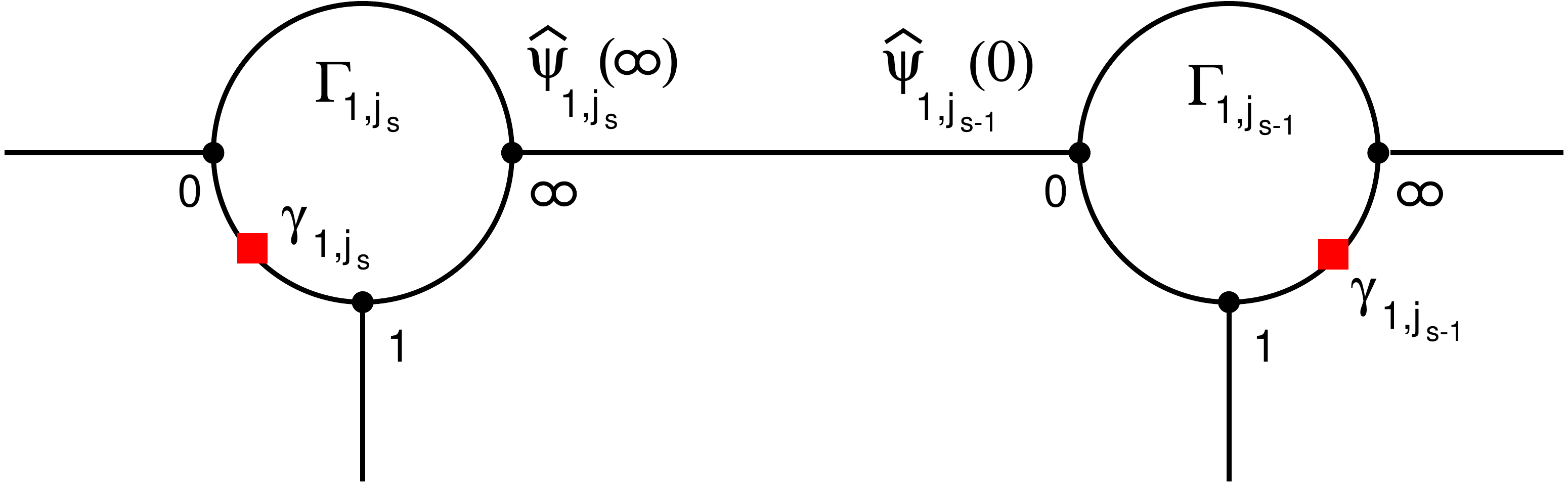}
      \caption{\small{\sl Intersection of two white components from the first row at one point. The values of the wave functions at the intersection points are equal, see Equation~\ref{eq:relation_black2}.}}\label{fig:one_black2}
    \end{figure}

\item For a sink $\kappa_{j_s}$, $1\le s \le n-k-1$, the value of the wave function at this sink coincides with the value of the wave function at the point of the white vertex connected to it, see Figure~\ref{fig:one_sink} right:
\begin{equation}    
\label{eq:relation_sink1}
\hat\psi(\kappa_{j_s}, \vec t) =  \hat\psi_{i_{l},j_{s}}(1, \vec t).
\end{equation}

\item For the sink $\kappa_n$ the value of the wave function at this sink coincides with the values of the wave function at all  points of the white vertices connected to it, see Figure~\ref{fig:one_sink} left :
\begin{equation}    
\label{eq:relation_sink2}
\hat\psi(\kappa_{n}, \vec t) =  \hat\psi_{1,n-1}(0, \vec t)= \ldots = \hat\psi_{k,n-1}(0, \vec t).
\end{equation}
\begin{figure}[H]
  \centering\includegraphics[width=0.6\textwidth]{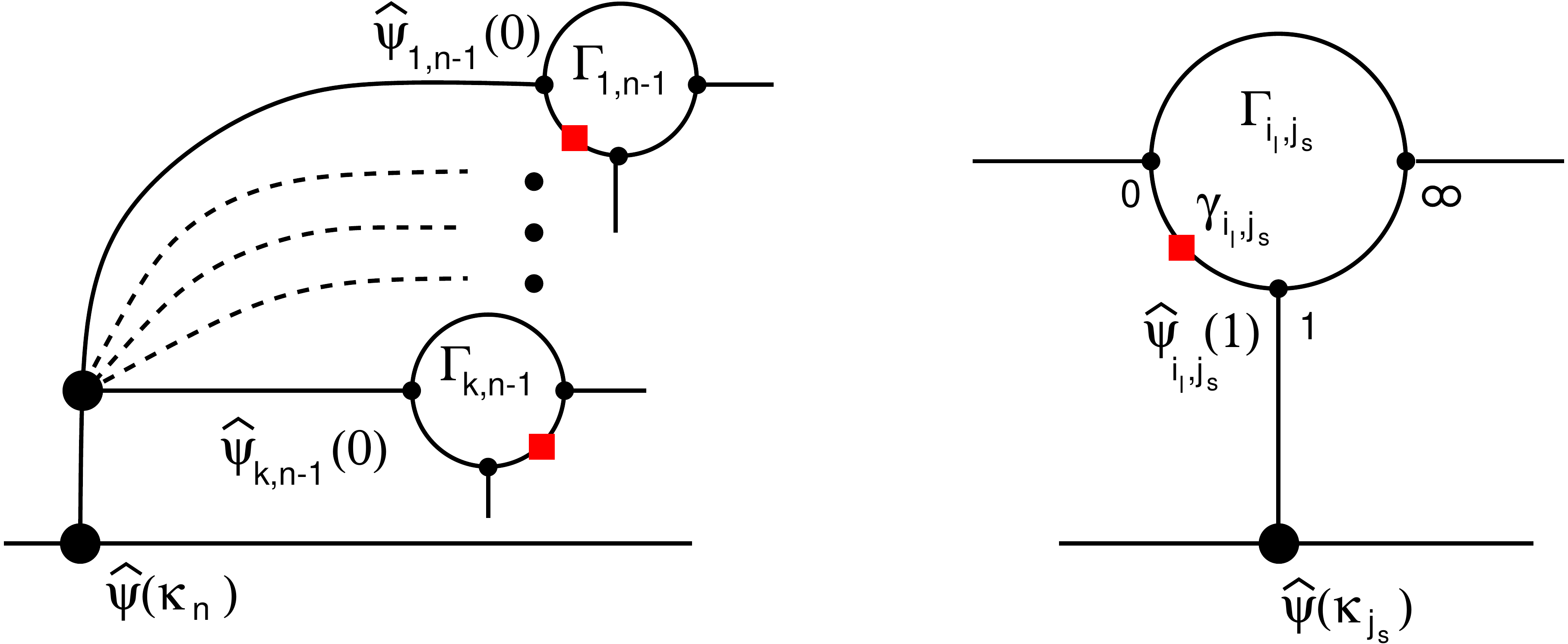}
      \caption{\small{\sl On the left: relations at the sing $\kappa_n$, see Equation~\ref{eq:relation_sink2}. On the right: relations at the sing $\kappa_{j_{s}}$, for $j_s<n$ see Equation~\ref{eq:relation_sink1}.}}\label{fig:one_sink}
    \end{figure}

\end{enumerate}

These conditions provide a full-rank linear system of relations in $3(g-k) + n$ variables. We have
\begin{enumerate}
\item  $(g-k)$ equations associated to white vertices,
\item $k$ equations associated to the sink $\kappa_n$,
\item  $2(g-n+1)$ associated to the trivalent black vertices,
\item  $n-k-1$ conditions associated to the sinks different from $\kappa_n$,
\item  $n-k-1$ conditions associated to the edges in the first line.
\end{enumerate}

Therefore the total number of equations is $3(g-k) + k$. 

Starting from the last line of the diagram and moving left to right and from bottom to top, we may express all the $3(g-k)$ variables associated to the white vertices in function of the values of the wave function at the boundary sinks. The remaining $k$ equations then provide the values of the wave functions at the boundary sources in function of the boundary sinks.

The resulting equations take the form
\begin{equation}
\label{eq:bm1}  
\hat\psi(\kappa_{i_l},\vec t)+ \sum_{s=1}^{n-k}b_{l,s} \hat\psi(\kappa_{j_s},\vec t)=0, 
\end{equation}
where $l=1,\ldots,k$, $i_l$ enumerate the boundary sources while  $s=1,\ldots,n-k$, $j_s$ enumerate the boundary sinks,   and the coefficients $b_{l,s}$ depend on the positions of the divisor points $\gamma_{i,j}$, associated to the white vertices.

In case of singular divisors the matrix $b_{l,s}$ may become singular; moreover it does not depend on the position of the divisor at the components $\Gamma_0$. Therefore this matrix does not contain the full information about the boundary measurement map. Let us observe that on $\Gamma_0$ the function $\hat\psi(\zeta,\vec t)$ can be written as
\begin{equation}    
\label{eq:psi_on_gamma}
\hat\psi(\zeta, \vec t) = \frac{\psi(\zeta, \vec t)}{P(\zeta)},
\end{equation}
where
\begin{equation}    
\label{eq:psi_on_gamma2}
\psi(\zeta, \vec t)=\left(\zeta^k - \sum_{j=0}^{k-1}w(\vec t) \zeta^j \right)\exp\left(\zeta x + \zeta^2 y + \zeta^3 t\right), 
\end{equation}
\begin{equation}    
\label{eq:funct_P}
P(\zeta)= \psi(\zeta, \vec t_0)= \prod_{j=1}^k (\zeta-\gamma_j)\exp\left(\zeta x_0 + \zeta^2 y_0 + \zeta^3 t_0 \right).
\end{equation}
where $\gamma_j=\gamma_j(\vec t_0)$ are the coordinates of the divisor points belonging to $\Gamma_0$ for the normalization time $\vec t_0$.

\begin{figure}[H]
  \centering\includegraphics[width=0.8\textwidth]{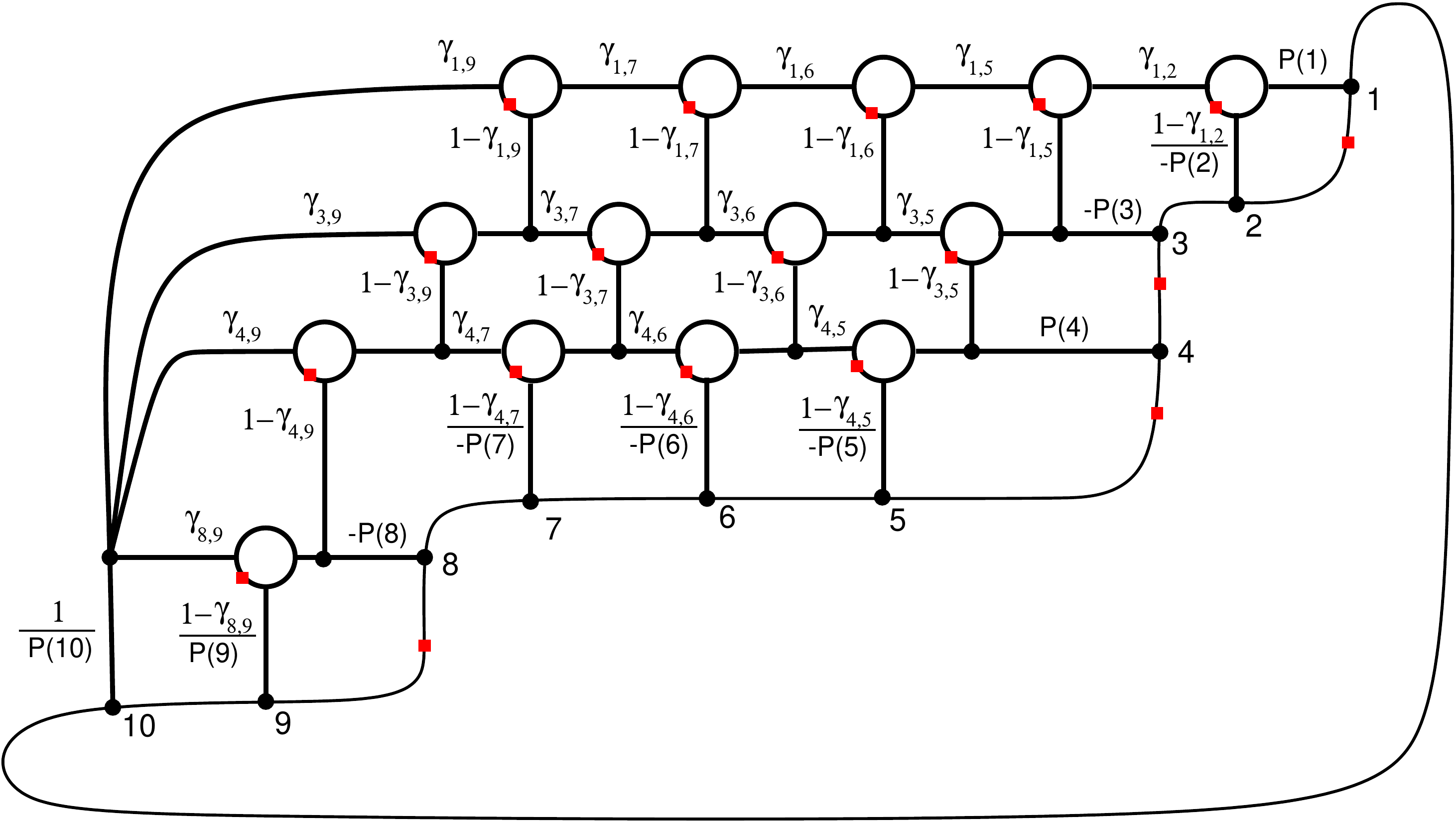 }
    \caption{\small{\sl This Figure illustrates Lemma~\ref{lem:weights}.}}\label{fig:network1}
  \end{figure}

\begin{lemma}\label{lem:weights}
Equation (\ref{eq:bm1}) is equivalent to  
\begin{equation}
\label{eq:bm2}  
\sum_{j=1}^{n}a_{l,j} \psi(\kappa_{j},\vec t)=0, \ \ l=1,\ldots,k
\end{equation}
where $A=[a_{ij}]$ is the matrix in the reduced row echelon form generated by the boundary measurement map of the canonically oriented reduced Le-network of the cell $S$ (canonical means that all horizontal edges are oriented right to left, and all vertical edges are oriented up to down) with weights:
\begin{enumerate}
\item The horizontal edge starting at the white vertex $\Gamma_{i,j}$, carries the weight $\gamma_{i,j}$.
\item The horizontal edge starting at a black vertex carries the weight $1$.
\item The horizontal edge starting at the boundary source $\kappa_{i_l}$ carries the weight $(-1)^{k+l-1} P(\kappa_{i_l})$.
\item The vertical edge starting at the white vertex $\Gamma_{i,j}$ and ending at a black vertex carries the weight $1-\gamma_{i,j}$.   
\item The vertical edge starting at the white vertex $\Gamma_{i_s,j}$ and ending at the boundary sink $\kappa_{j}$ carries the weight $(-1)^{k+s}(1-\gamma_{i_s,j})/P(\kappa_j)$.
\item The vertical edge starting at a black vertex and ending at the boundary sink $\kappa_{n}$ carries the weight $1/P(\kappa_n)$.  
\end{enumerate}
\end{lemma}
For an example see Figure~\ref{fig:network1}.

 \begin{remark}\label{rem:6}
In the special case $S=Gr^+(1,2)$ we add a bivalent black vertex to the edge connecting the boundary point $\kappa_1$, $\kappa_2$ to make the rules of Lemma~\ref{lem:weights} applicable to this case.
 \end{remark}   
\begin{figure}[H]
  \centering\includegraphics[width=0.2\textwidth]{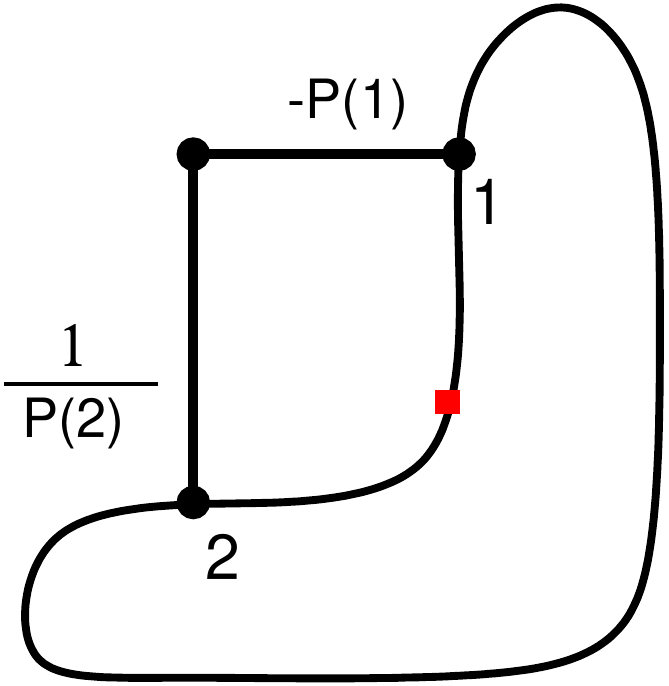 }
    \caption{\small{\sl This Figure illustrates Remark~\ref{rem:6}.}}\label{fig:network2}
  \end{figure}
  The proof is straightforward, since it is sufficient to substitute the right-hand side of (\ref{eq:psi_on_gamma}) into the system of equations fulfilled by the KP wave function to check (\ref{eq:bm2}). Indeed, the boundary measurement map provides the value of the maximal minors of the matrix $A$ associated to this network. In particular the element
  $$
  a_{l,j}=(-1)^{\sigma_{l,j}} \Delta_{lj},
  $$
  where $\sigma_{lj}$is the number of boundary sources in the open interval $]i_l,j[$, and $\Delta_{lj}$ is the minor of the columns $\{i_1,\ldots,i_k\}\backslash \{i_l\} \cup \{j\}$. Following Postnikov \cite{Postnikov2006}, this minor is computed via discrete integration summing the weights of all directed paths starting at the source $i_l$ and ending at the sink $j$. The weight of a path is the product of the weights of all edges belonging to such path. 
  
\begin{lemma}\label{lem:darboux}
The Darboux transformation generated by the matrix $A$, computed in Lemma~\ref{lem:weights}, provides the wave function satisfying (\ref{eq:bm2}). 
\end{lemma}
The proof follows directly from the Darboux transformations theory, see, for example, Equation~6.5 in \cite{Abenda2021}

\begin{lemma}\label{lem:positivity}
The matrix $A$ defined in Lemma~\ref{lem:weights} is totally non-negative.
\end{lemma}
\begin{proof}
  The weights defined in Lemma~\ref{lem:weights} may be either positive or negative, and the signs depend on the position of the divisors. Let us show that for divisors satisfying the Dubrovin-Natanzon conditions the weights of all paths connecting sources and sinks are positive.

 Consider a path $\cal P$ connecting the source $\kappa_{i_l}$ to the sink $\kappa_j$. Its weight $w({\cal P})$ satisfies
\begin{equation}\label{eq:parity1}
w({\cal P)} = \frac{w({\cal P}_{ij})\, P (\kappa_{i_l}) \, (-1)^{k+l+1}}{P(\kappa_j)\,  (-1)^{k+s}}
\end{equation}
where $w({\cal P}_{ij})$ is the product of the weights of type $\gamma_{ij}$, $(1-\gamma_{ij})$ on the path, and $l$, $s$ are as in Items 3 and 5 in Lemma~\ref{lem:weights}.

Let us point out that
\begin{equation}\label{eq:parity2}
\frac{ (-1)^{k+l+1}}{(-1)^{k+s}} = (-1)^{\sigma(i_l,j)+1}.
\end{equation}

Let us denote:
\begin{enumerate}
 \item OA = \# of ovals in the area to the right of $\cal P$,
 \item WA = \# of white vertices in in the area to the right of $\cal P$,
 \item D0 = \# of divisor points on $\Gamma_0$ belonging to the area to the right of $\cal P$,
 \item DP = \# of divisor points on $\cal P$. 
\end{enumerate}
First of all, 
\begin{equation}\label{eq:parity3}
 \mathrm{OA} = \mathrm{WA} + \sigma(i_l,j)+1,
\end{equation}
Since the divisor satisfies Dubrovin-Natanzon conditions, the following identity holds
\begin{equation}\label{eq:parity4}
\mathrm{OA -WA = D0 + DP}\quad   (\mathrm{mod}\ 2).
\end{equation}

Next
\begin{equation}\label{eq:parity5}
\sign (P (\kappa_{i_l}) )\, \sign (P(\kappa_j))  = (-1)^{\mathrm{D0}},
\end{equation}
\begin{equation}\label{eq:parity6}
\sign(w({\cal P}_{ij})) = (-1)^{\mathrm{DP}}.
\end{equation}
therefore, inserting equations (\ref{eq:parity2})--(\ref{eq:parity6}) into (\ref{eq:parity1}), we obtain that 
\begin{equation}
\begin{split}  
  \sign(w({\cal P})) &= \sign(w({\cal P}_{ij}))\, \sign (P (\kappa_{i_l}) )\, \sign (P(\kappa_j))\, (-1)^{\sigma(i_l,j)+1} =\\
                      &= (-1)^{\mathrm{D0+DP}} \, (-1)^{\sigma(i_l,j)+1}=\\
                      & =(-1)^{\mathrm{OA-WA}}\,(-1)^{\sigma(i_l,j)+1} = 1.
\end{split}  
\end{equation}
\end{proof}

\begin{figure}[H]
  \centering\includegraphics[width=0.8\textwidth]{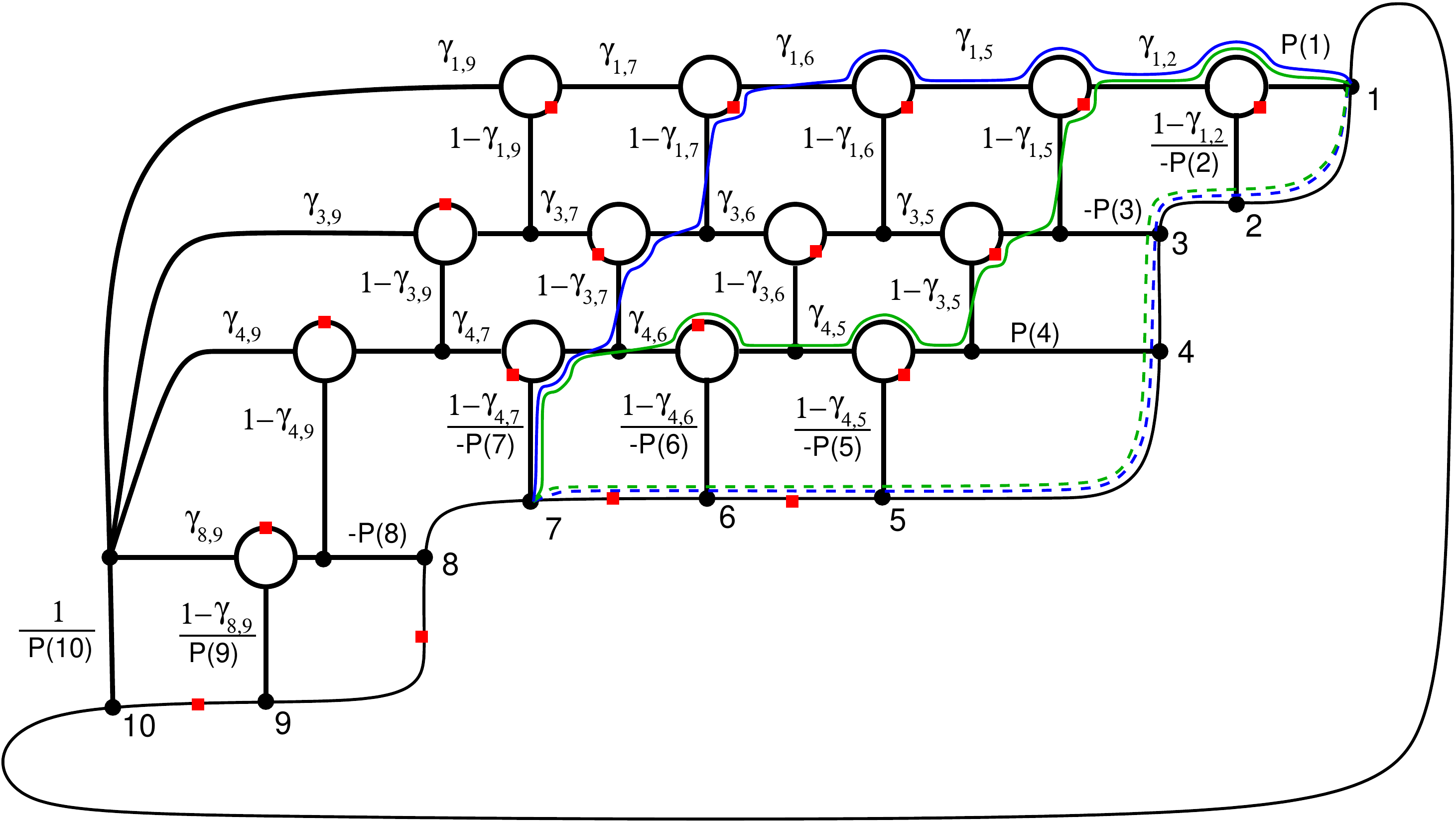 }
    \caption{\small{\sl Along each directed path the product of weights is positive.}}\label{fig:network3}
  \end{figure}

We invite the reader to check these relations for the blue and green paths in Figure~\ref{fig:network3}.

Since the boundary measurement map associates to each smooth Dubrovin-Natanzon configuration of divisor points a unique point of the positroid cell $S$, it is always possible to gauge the edge weights in such a way that they are all positive. 

\begin{figure}[H]
  \centering\includegraphics[width=0.8\textwidth]{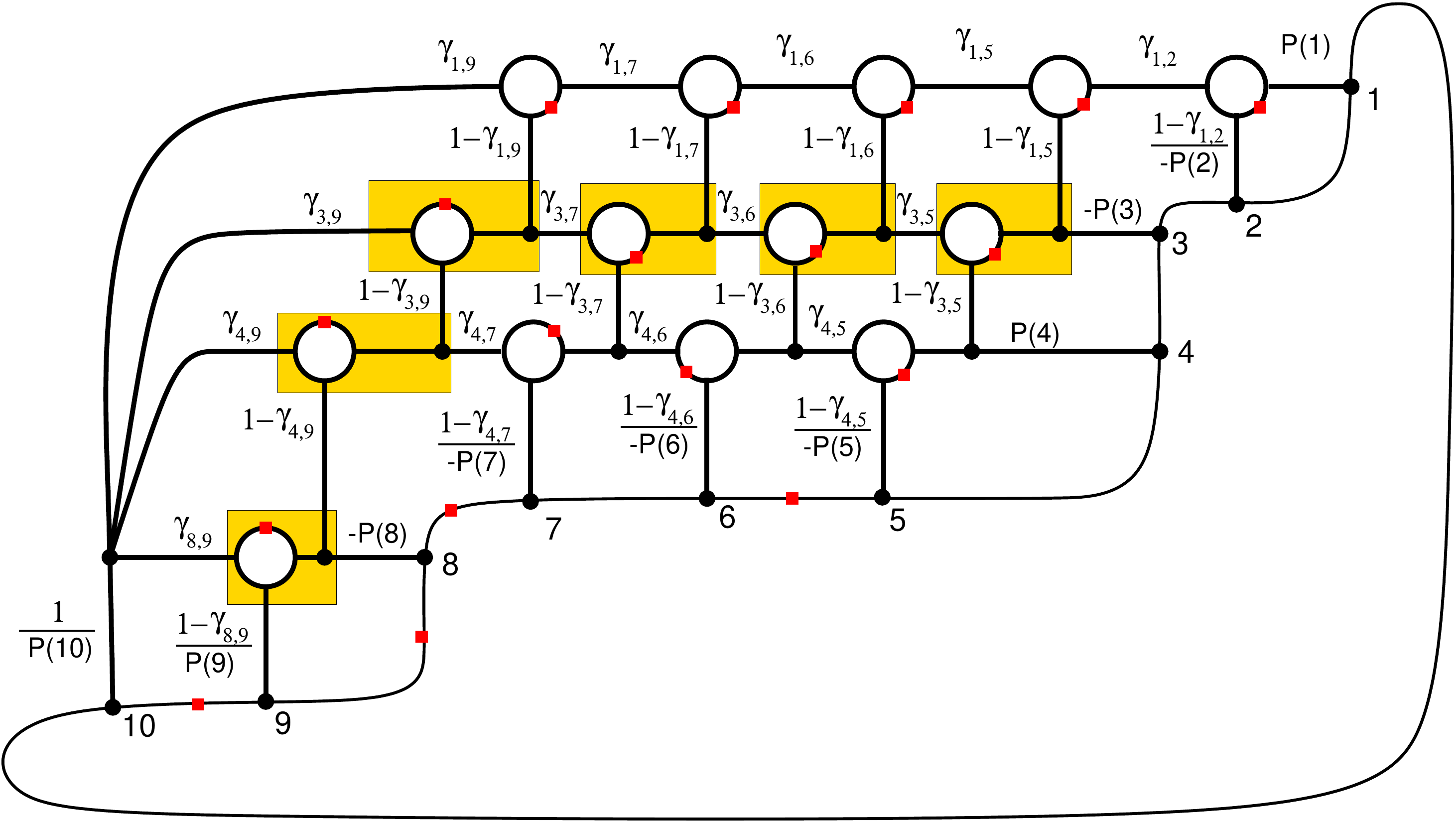 }
    \caption{\small{\sl We illustrate Lemma~\ref{lem:pos_weights1} marking in yellow the pairs of black and white vertices where we apply the weight gauge.}}\label{fig:network4}
  \end{figure}

\begin{figure}[H]
  \centering\includegraphics[width=0.8\textwidth]{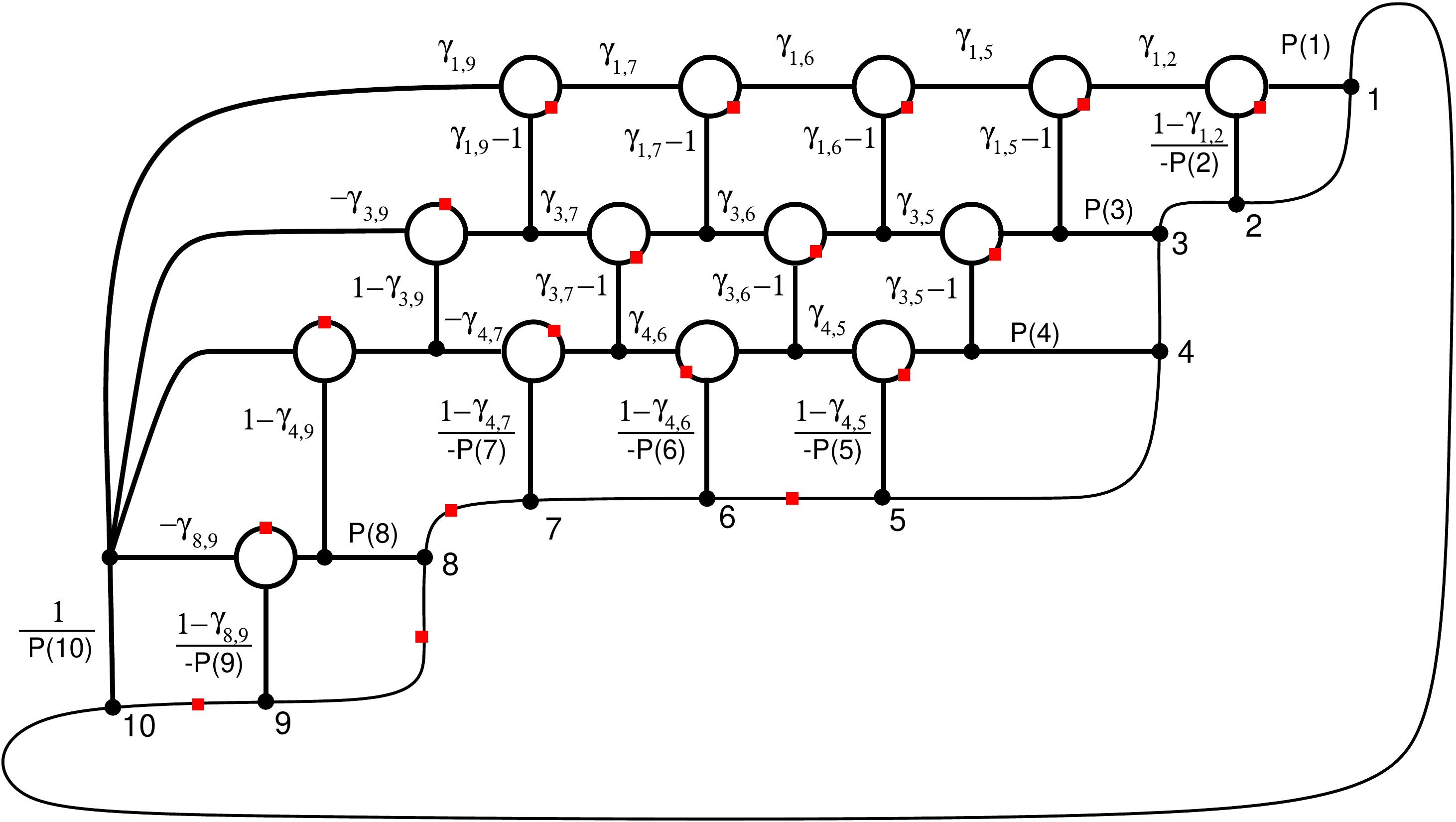 }
    \caption{\small{\sl We illustrate Lemma~\ref{lem:pos_weights1} by showing that after applying the weight gauge to the pairs of vertices marked at the Figure~\ref{fig:network4} all   edges carry positive weights.}}\label{fig:network5}
  \end{figure}

\begin{lemma}\label{lem:pos_weights1}
There exists a weight gauge transformation of the network defined in Lemma~\ref{lem:weights} such that all edge weights
are positive. In particular, if all divisor points on the white components are located in the interval $]0,1[$ and the divisor points on $\Gamma_0$ are located in the intervals $]i_s,i_s+1[$, $s=1,\ldots, k$, the weights defined in Lemma~\ref{lem:weights} are all positive (for an example see Figure~\ref{fig:network1}).
\end{lemma}
\begin{proof}
  The construction of the  weight gauge transformation is explicit and it is done traveling along the network. One starts from the first row and  moves from right to left inside the row. Then one moves to the next row downwards, and repeats the procedure until one exhausts all rows.

All weights in the first row have the form $\gamma_{1,j}$ and are positive due to the Dubrovin-Natanzon conditions. The weight of the edge starting at the boundary source 1 is always positive for the same reason. The weight of each vertical edge starting at a white vertex $\Gamma_{1,j}$ and reaching the boundary vertex $j$ is also positive. In this way one exhausts  the first row and moves to the second one.

Lemma~\ref{lem:positivity} implies that the signs of all paths starting from boundary sources and reaching the same black vertex are equal. If one starts from the boundary source $i_2$ and moves horizontally till the first black vertex, then either both weights of the horizontal edge starting at $i_2$ and of the vertical edge connected to this black vertex are positive or both are negative. If they are positive, one proceeds to the next black vertex in this row. If they are negative, one applies the following weight gauge transformation assigning weight $-1$ to this black vertex and to the white vertex to its left. Following the gauge rule transformation in \cite{Postnikov2006} these two edges become positive, the weight of the horizontal edge connecting these black and the white vertices remains $+1$, and the other two edges starting at this white vertex change the signs. (See Figure~\ref{fig:gauge1} for an example). Moreover, if the vertical edge from this white vertex reaches the boundary, it is positive due to Lemma~\ref{lem:positivity}.
  
\begin{figure}[H]
  \centering\includegraphics[width=0.7\textwidth]{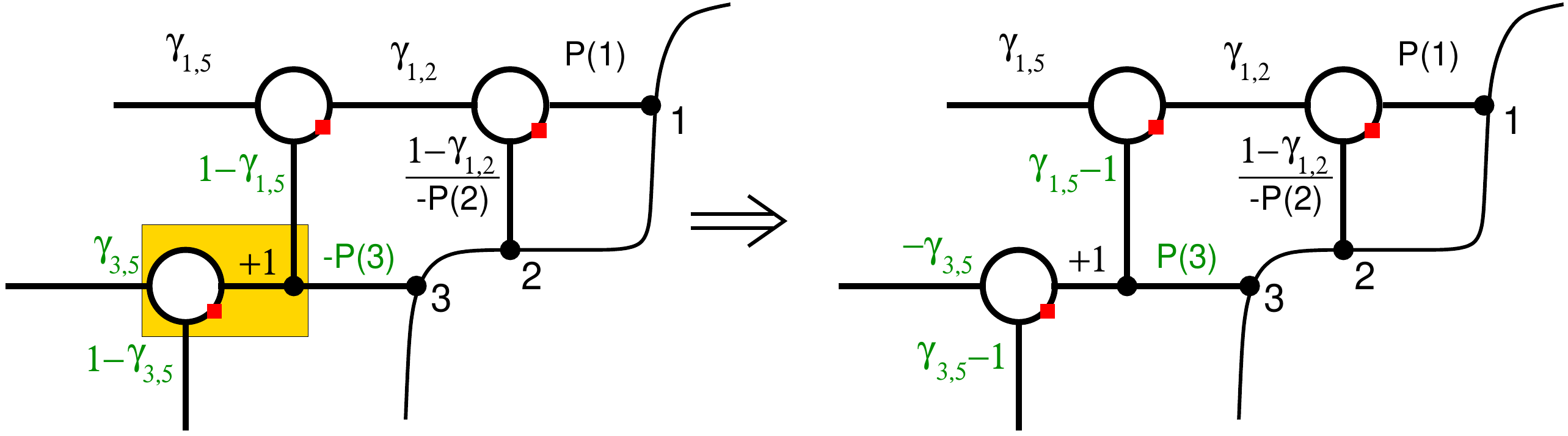 }
  \caption{\small{\sl We illustrate the results of the applying the weight gauge transformation to a pair of black and white vertices inside the yellow box. The weights changing sign after the weight gauge are marked green.}}\label{fig:gauge1}
  \end{figure}   

  Then one moves along the second row to the next black vertex and checks the signs of the edges ending at this black vertex. Again, by Lemma~\ref{lem:positivity} these signs coincide. If they are positive, one proceeds to the next pair. If they are negative, one applies the same weight gauge transformation to the corresponding pair of black and white vertices. When one reaches the black vertex connected to the boundary sink $\kappa_n$, all weights of the second row are positive.

  Then one proceeds to the next row and repeats the procedure until the last row is completed.

The proof of the last statement of the Lemma is straightforward.
  
\end{proof}
We illustrate Lemma~\ref{lem:pos_weights1} by an example in Figures~\ref{fig:network4},~\ref{fig:network5}.

In the following Lemma we prove that the map from the Dubrovin-Natanzon Cartier divisors to the positroid cell $S$ is well defined and regular also for singular divisors. 
  
\begin{lemma}\label{lem:pos_weights2}
Under proper weight gauge the edge weights remain positive and finite also for singular Dubrovin-Natanzon Cartier divisors. 
\end{lemma}
\begin{proof}
  The first step is to check that at any singular configuration one can define a weight gauge transformation such that the weights are finite, non-zero and continuously dependent on the divisor in the Cartier sense. There are only three cases to analyze:

\begin{enumerate}  
\item If three divisor points approach a triple singularity which involves two white components
  $\Gamma_{i_{l-1},j}$,  $\Gamma_{i_{l},j}$ and the boundary source $\kappa_{i_l}$, one has: 
\begin{equation}
\gamma_{i_{l-1},j}\rightarrow 1, \ \ \gamma_{i_{l},j}\rightarrow \infty, \ \ P(\kappa_{i_l}) \rightarrow 0.  
\end{equation}
Therefore the edge weights defined in Lemma~\ref{lem:weights} become singular (see Figure~\ref{fig:gauge2} left.) If one applies the weight gauge transformation assigning to the black vertex corresponding to this triple point and to the white vertex  $\Gamma_{i_{l},j}$ the weight $(-1)^{k+l+1}  P(\kappa_{i_l})$, then the resulting weights are as in  Figure~\ref{fig:gauge2} right. Therefore, in agreement with the standard blow-up procedure defined in Equation~(\ref{eq:cartier3}), the new weights have finite limit and are regular.
\begin{figure}[H]
  \centering\includegraphics[width=0.8\textwidth]{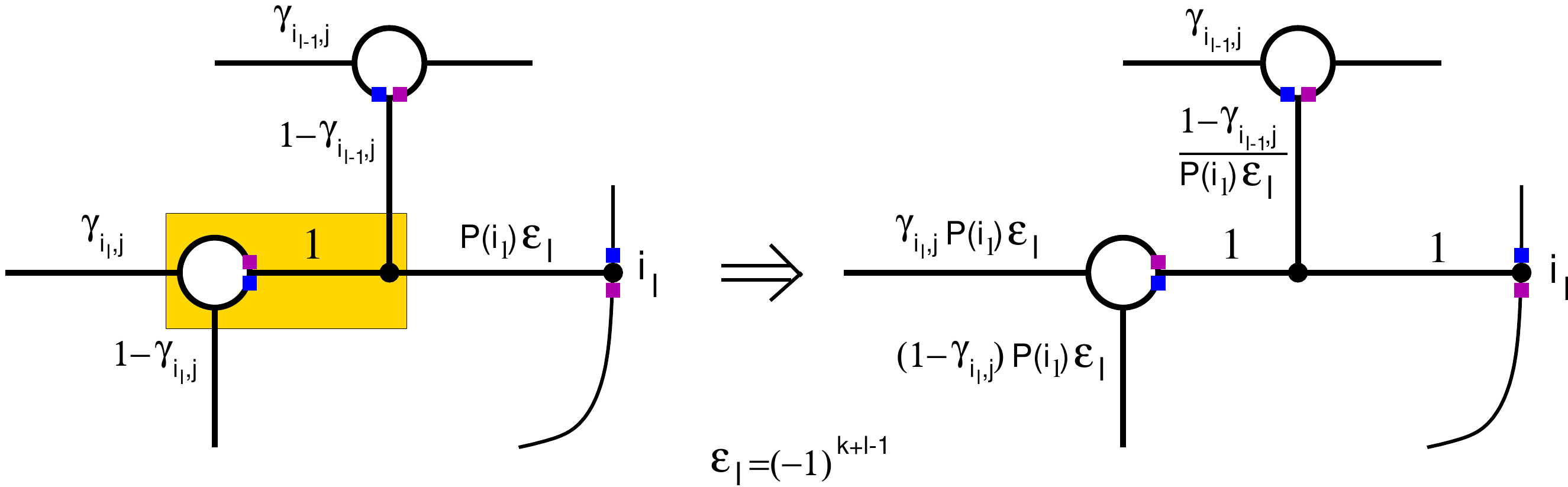 }
  \caption{\small{\sl We illustrate the effect of the weight gauge transformation at a triple singularity involving a boundary vertex. The yellow box marks the two vertices to which we apply the transformation.}}\label{fig:gauge2}
  \end{figure} 
   
\item If three divisor points approach a triple singularity which involves three white components
  $\Gamma_{i_{l-1},j}$,  $\Gamma_{i_{l},j}$,  $\Gamma_{i_{l},j-1}$, then
\begin{equation}
\gamma_{i_{l-1},j}\rightarrow 1, \ \ \gamma_{i_{l},j}\rightarrow \infty, \ \ \gamma_{i_{l},j-1} \rightarrow 0.  
\end{equation}
Again the edge weights defined in Lemma~\ref{lem:weights} become singular (see Figure~\ref{fig:gauge3} left.) If one applies the weight gauge transformation assigning to the black vertex corresponding to this triple point and to the white vertex  $\Gamma_{i_{l},j}$ the weight $\gamma_{i_{l},j-1}$, then the resulting weights are as in  Figure~\ref{fig:gauge3} right. Again, in agreement with the standard blow-up procedure defined in Equation~(\ref{eq:cartier3}), the new weights have finite limit and are regular. 
\begin{figure}[H]
  \centering\includegraphics[width=0.95\textwidth]{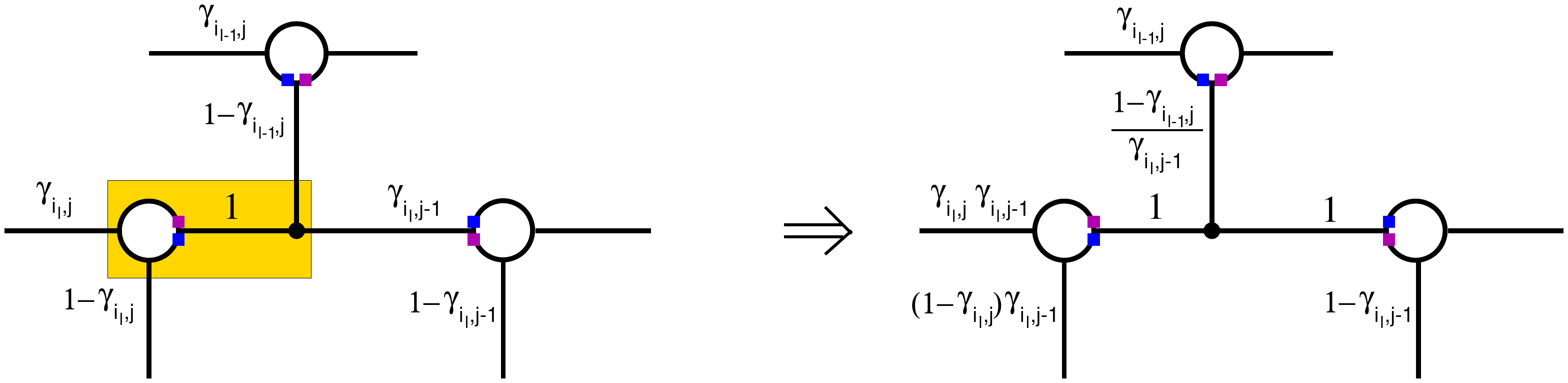 }
  \caption{\small{\sl  We illustrate the effect of the weight gauge transformation at a triple singularity involving three white vertices. The yellow box marks the two vertices to which we apply the transformation.}}\label{fig:gauge3}
  \end{figure} 

\item  If two divisor points approach a double singularity which involves the white component 
$\Gamma_{i_{l},j}$ and the boundary vertex $\kappa_j$ then
 \begin{equation}
\gamma_{i_{l},j}\rightarrow 1,  \ \ \ P(\kappa_{j}) \rightarrow 0.   
\end{equation}
In this case the edge weights defined in Lemma~\ref{lem:weights} are finite and regular in agreement with the standard blow-up procedure defined in Equation~(\ref{eq:cartier2}) 

\begin{figure}[H]
  \centering\includegraphics[width=0.3\textwidth]{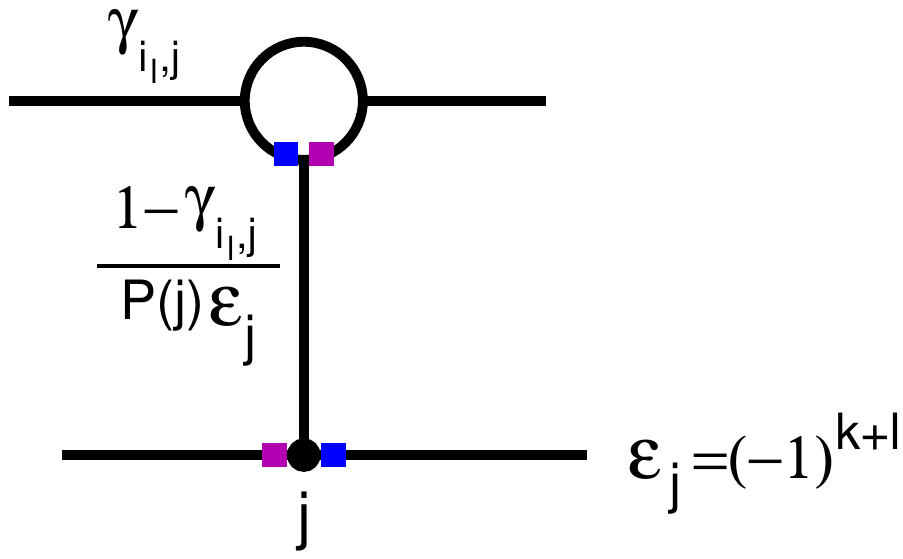 }
  \caption{\small{\sl  We illustrate the case of double singularity involving a boundary vertex.}}\label{fig:gauge4}
  \end{figure} 

\end{enumerate}
After applying the above procedure to all points of the divisor located at double and triple points of the curve, all edge weights of the network become finite. Let us remark that the above gauge transformation preserves the positive sign of every directed path from a source to a sink. Therefore, proceeding as in the proof of Lemma~\ref{lem:pos_weights1} one may transform all edge weights to positive ones. This completes the proof.

\end{proof}  

\begin{appendices}

\section{Totally non-negative Grassmannians and Le-networks}

In our paper we represent points of the Grassmannian $Gr(k,n)$, $k<n$, by $k\times n$ matrices of maximal rank $k$. We denote the equivalence class of the matrix $A$ with respect to the standard left action of $GL(k)$ by $[A]$.

Let us recall
\begin{definition}\label{def:tp} 
  A  point $[A]$ of the \textbf{real} Grassmannian $Gr(k,n)$ belongs to its totally non-negative part  $\GTNN$, if all non-zero maximal minors of $A$ (Pl\"ucker coordinates) share the same sign. Without loss of generality we may assume that all minors of $A$ are non-negative.
\end{definition}

\subsection{Positroid cells}

Our construction of ${\mathtt{MM}}$-curves is based on the positroid cell decomposition of $\GTNN$, constructed by A. Postnikov in \cite{Postnikov2006}.

\begin{definition}
Positroid cells are the sets of points of a totally non-negative Grassmannian sharing the same sets of strictly positive maximal minors (Pl\"ucker coordinates). Equivalently, the positroid cells are the intersections of Gelfand-Serganova strata \cite{GelfandSerganova1987} with the totally non-negative part of the Grassmannians.
\end{definition}

Let us recall the combinatorial description of the positroid cells introduced in \cite{Postnikov2006}. 

Each point of the Grassmannian can be uniquely represented by a matrix in reduced row echelon form using Gauss reduction, and it is represented by a Young tableau. The submatrix associated with the pivot columns of the  reduced row echelon form is the $k\times k$ identity matrix. In the l-th row of the matrix the elements which are located in the non-pivot columns to the right of the l-th pivot column may take non-zero values. To such point of the Grassmannian we associate a Young diagram in the English notation. The number of boxes in the l-th row of this Young diagram is equal to the number of non-pivot columns to the right of the $l$-th pivot one. If the south-eastern boundary of the Young diagram is enumerated from $1$ to $n$ in the increasing order from north-est to south-west, then the vertical edges label the pivot columns and the horizontal edges label the non-pivot columns.

Each Schubert cell is the set of all Grassmannian points sharing the same Young diagram, and the dimension of the cell is equal to the number of boxes in this diagram.

\begin{figure}[H]
\centering{
\resizebox{10cm}{!}{%
$
\begin{bmatrix} 
1 & a_{1\,2} & 0  &a_{1\,4} & a_{1\,5} & 0 & a_{1\,7} & 0  & a_{1\,9} & a_{1\,10}\\
0 & 0 & 1 & a_{2\,4} & a_{25} & 0 & a_{2\,7} & 0 & a_{2\,9} & a_{2\,10} \\
0 & 0 & 0 & 0 & 0 & 1 & a_{3\,7} & 0& a_{3\,9} & a_{3\,10}\\
0 & 0 & 0 & 0 & 0 & 0 & 0 & 1 &a_{4\,9} & a_{4\,10}
\end{bmatrix}
$
}
\hspace{2mm}\parbox{4.4cm}{\vphantom{t}\\ \includegraphics[width=3.8cm]{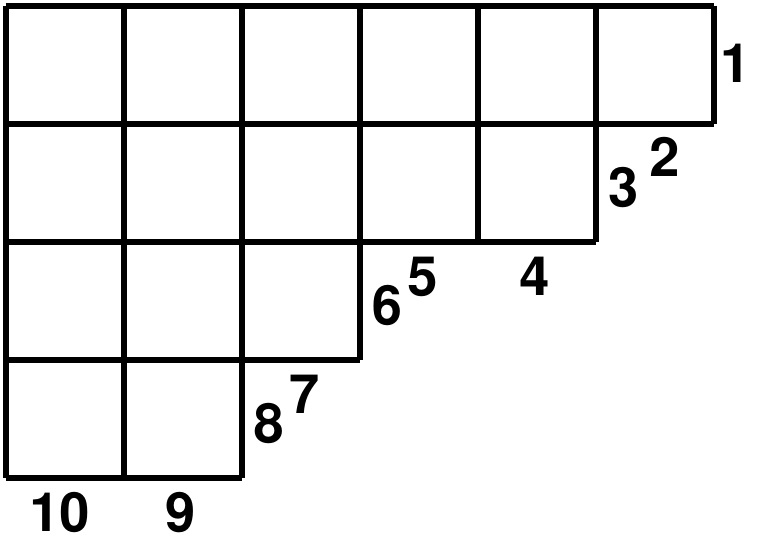}}
\caption{\small{\sl On the left: A matrix representing a point in $Gr(4,10)$ written in reduced row echelon form. On the right: the corresponding  Young diagram (English notation). Pivot columns are 1, 3, 6, 8; non-pivot columns are 2, 4, 5, 7, 9, 10. \label{fig:appendix_1}} }}
\end{figure}

In Figure~\ref{fig:appendix_1} we show an example of a 16-dimensional Schubert cell in  $Gr(4,10)$ and its Young diagram.

\begin{figure}[H]
\centering\includegraphics[width=4cm]{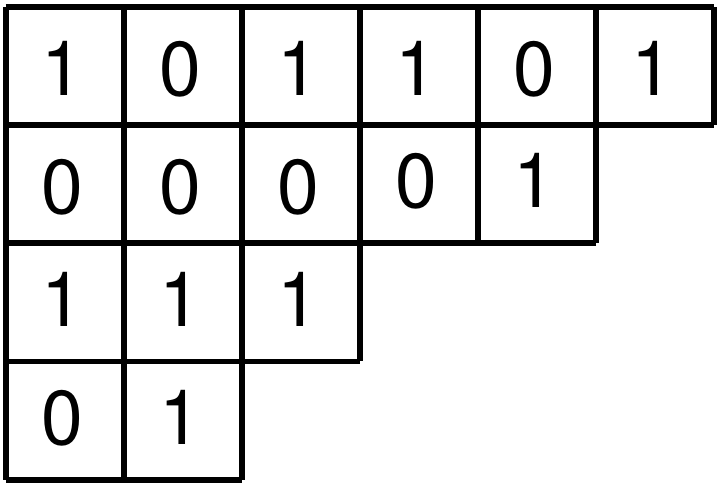}\hspace{1cm}\includegraphics[width=4.5cm]{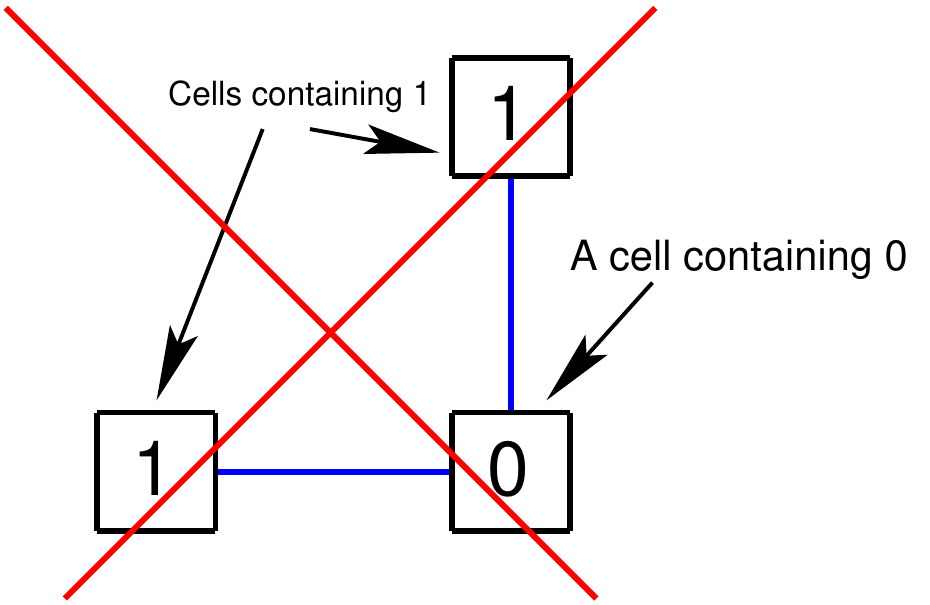}
\caption{\small{\sl On the left: a Young diagram filled by zeroes and ones complying the Le-rule; on the right: a configuration forbidden by the Le-rule.}} \label{fig:appendix_2}
\end{figure}

Following \cite{Postnikov2006}, a positroid cell is represented by a \textbf{Le-diagram}, which is a Young diagram filled by zeroes and ones fulfilling the \textbf{Le-rule:} if a box contains a zero then either all boxes in the same row to the left are filled by zeros or all boxes in the same column above are filled by zeros, see Figure~\ref{fig:appendix_2}.

\begin{definition}\label{def:5}
A positroid cell is called irreducible if in the representing Le-diagram every row and every column has a box filled by 1. 
\end{definition}

\begin{figure}[H]
\centering\includegraphics[width=4cm]{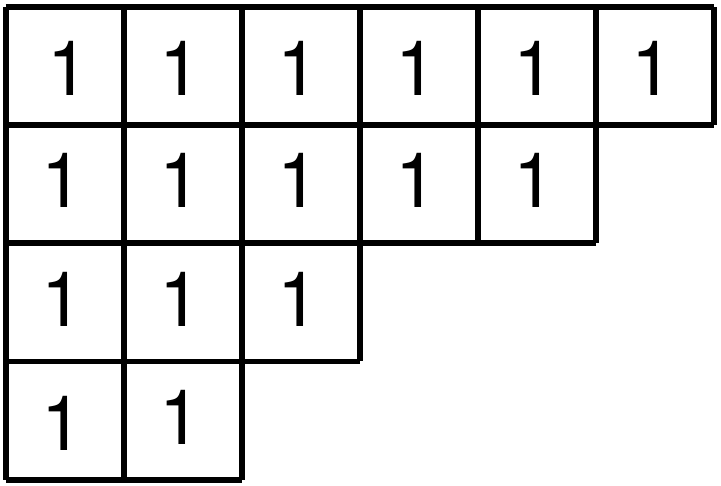}
\caption{\small{\sl An example of a  totally positive Schubert cell.}} \label{fig:appendix_3}
\end{figure}

\begin{definition}\label{def:6}
Irreducible positroid cells represented by Young diagrams filled only by ones are called \textbf{totally positive Schubert cells.} For an example see Figure~\ref{fig:appendix_3}.
\end{definition}
In particular, a totally positive Schubert cell in $\GTNN$ has exactly $n-k$ filled boxes in the upper line and $k$ filled boxes in the leftmost column. 

\begin{figure}[H]
\centering\includegraphics[width=4cm]{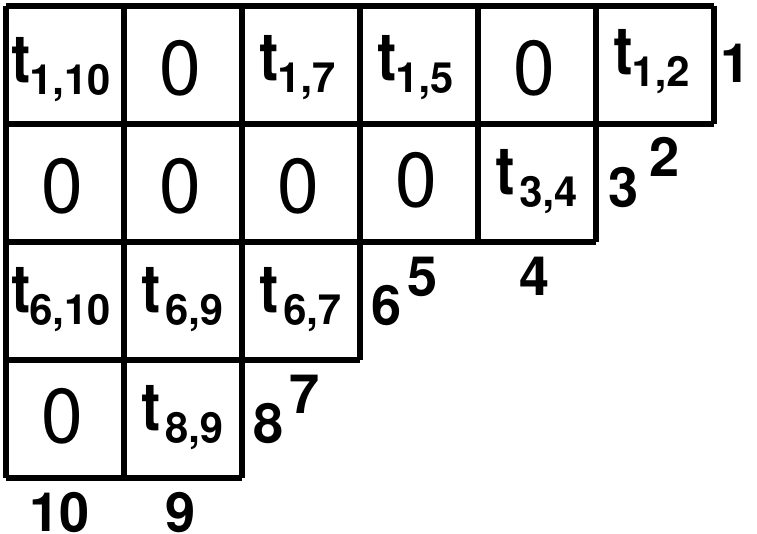}
\caption{\small{\sl A Le-tableau representing a point on $Gr^{\mbox{\tiny TNN}} (4,10)$}} \label{fig:appendix_4}
\end{figure}

Each point in the positroid cell is represented by a Young tableau where ones are substituted by positive weights \cite{Postnikov2006} (for an example see Figure~\ref{fig:appendix_4}). Such tableaux are called \textbf{Le-tableaux}.  The construction in \cite{Postnikov2006} provides a birational parametrization  of positroid cells in terms of Le-networks.

\begin{figure}[H]
  \centering\includegraphics[height=5cm]{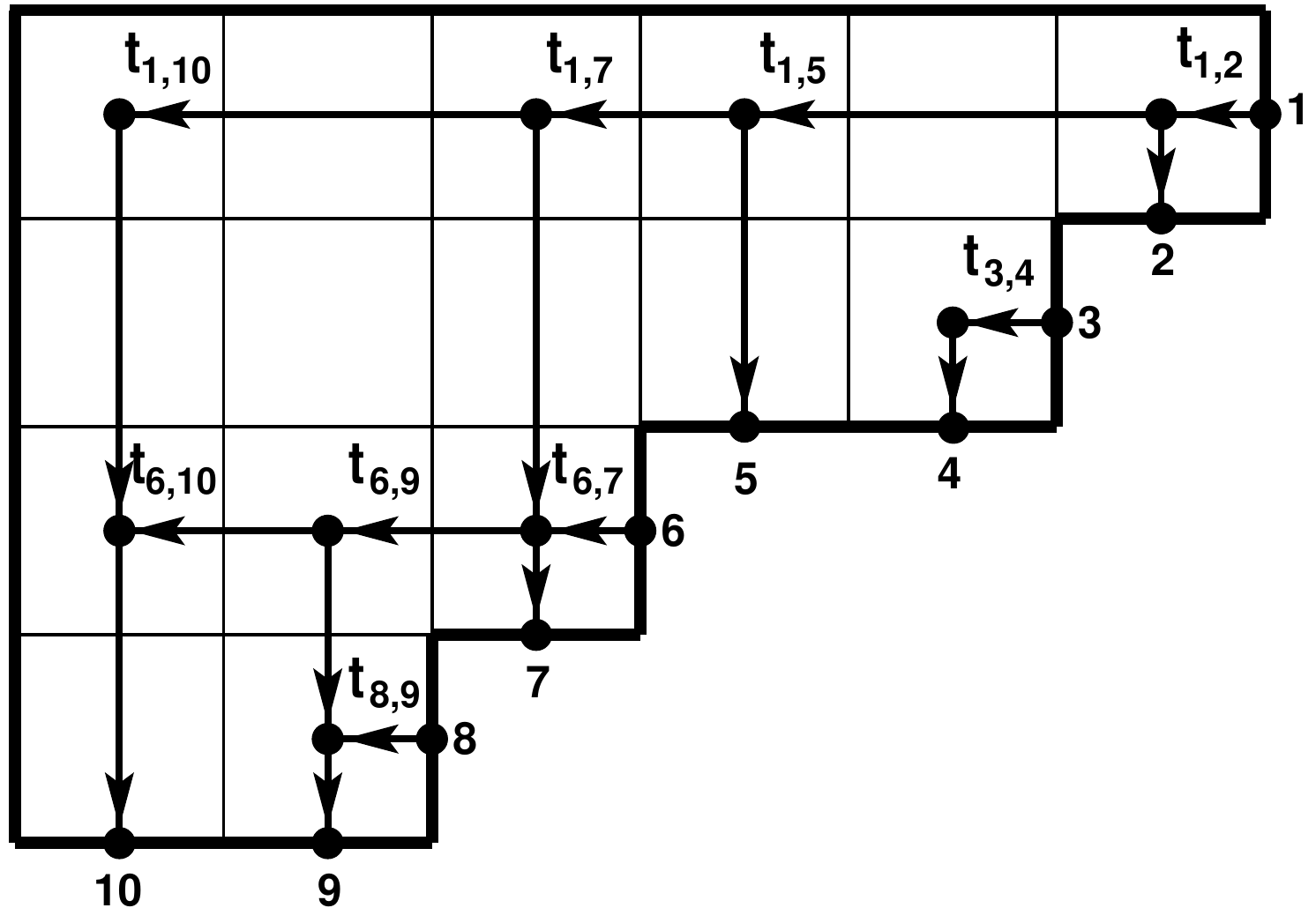}
    \caption{\small{\sl The planar network associated to the Le-tableau of Figure~\ref{fig:appendix_4}}}\label{fig:appendix_5}    
\end{figure}

Let us recall now the construction of the \textbf{Le-network} associated with the Le-tableau:
\begin{enumerate}
\item We erase the zeroes from the boxes of the Le-tableau.
\item We put a boundary vertex in the middle of each edge of the south-eastern boundary.
\item We put an internal vertex in the middle of each box containing a positive weight.  
\item From each internal vertex we draw an edge connecting it to the next vertex to the right, oriented from right to left, and assign the weight from this box to such edge.
\item From each internal vertex we draw an edge connecting it to the next vertex below, oriented downward, and assign weight $1$ to such edge. In the Figures we do not draw unit weights on vertical edges.
\end{enumerate}
In this way we construct an acyclically oriented planar network in the disk. The boundary of the disk corresponds to the boundary of the Young diagram (see Figure~\ref{fig:appendix_5}). 
The boundary vertices corresponding to the pivot columns are called \textbf{sources}, and the boundary vertices corresponding to the non-pivot columns are called \textbf{sinks}. In Figure~\ref{fig:appendix_5} the boundary vertices $1,3,6,8$ are sources and the boundary vertices $2,4,5,7,9,10$ are sinks.

\begin{figure}[H]
  \centering\includegraphics[height=2.5cm]{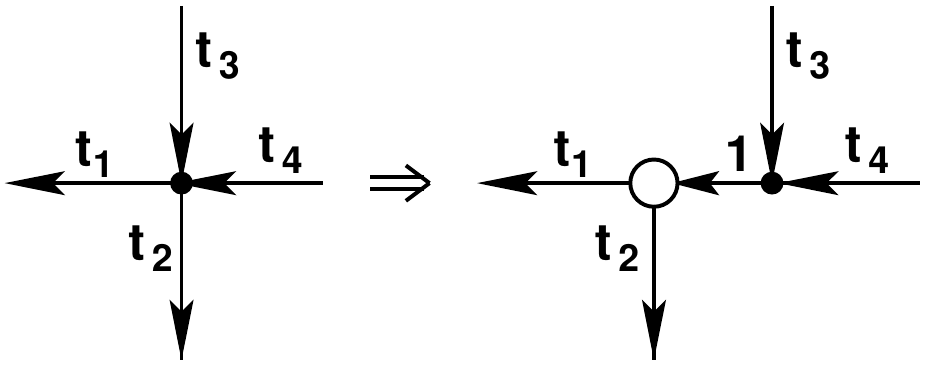}
    \caption{\small{\sl The transformation of the network at a four-valent vertex.}} \label{fig:appendix_6}     
\end{figure}

\begin{figure}[H]
  \centering\includegraphics[height=5cm]{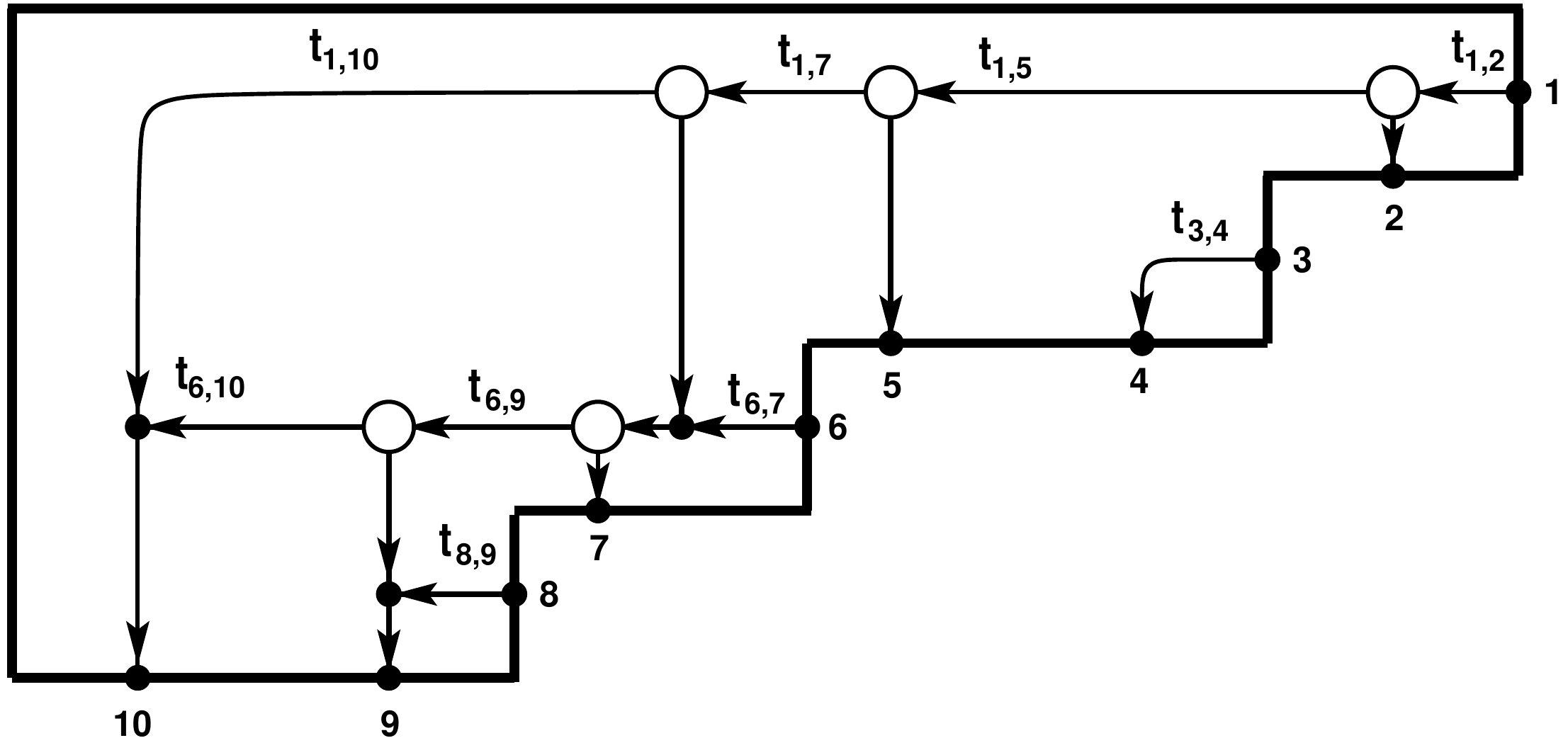}
    \caption{\small{\sl The result of transforming the planar network in Figure~\ref{fig:appendix_5}.}} \label{fig:appendix_7}  
\end{figure}

The next step is to transform this Le-network into a perfectly oriented planar bicolored network in the disk with internal vertices of valency 3.
\begin{enumerate}
\item We eliminate all internal bivalent vertices, and we assign to the new edge the product of weights on the initial edges.
\item If a trivalent vertex has exactly one incoming edge, we color it white.
\item If a trivalent vertex has exactly one outgoing edge, we color it black. 
\item We replace each four-valent vertex by a pair of black and white ones as in Figure~\ref{fig:appendix_6}; we keep the weights at the old edges and assign the unit weight to the new one. 
\end{enumerate}
An example of a perfectly oriented planar bicolored network is presented in Figure~\ref{fig:appendix_7}.

Let us recall the formulas from \cite{Postnikov2006} expressing the elements of the reduced row echelon form matrix $A$ represented by the canonically perfectly oriented Le-network introduced above. Let the pivot columns be $i_r$, $1\le r \le k$, $i_1<i_2\ldots<i_k$. 
For the pivot columns we have 
$$
A_{l,i_r}=\delta_{l,r}.
$$
For the non-pivot columns we have 
\begin{equation}\label{eq:bmm}
A_{lj} :=  (-1)^{\sigma(i_l,j)}  \sum\limits_{P:b_{i_l}\mapsto b_j} w(P),
\end{equation}
where the sum is over all directed paths from the source $b_{i_l}$ to the sink $b_j$, $w(P)$ is the product of the edge weights of $P$, and  $\sigma(i_l,j)$ is the number of sources strictly between $i_l$ and $j$.

\begin{remark}
It is well-know that Schubert cells provide CW decomposition of the classical Grassmannians.  In \cite{PostnikovSpeyerWilliams2009} it was shown that the positroid cells also form a CW-complex; in \cite{GalashinKarpLam2018} it was shown that $\GTNN$  is homeomorphic to a ball of dimension $k\times(n-k)$. 
\end{remark}

\end{appendices}

\section{Competing interests}
No competing interest is declared.

\section{Author contributions statement}

Both authors (S.A and P.G.G) equally contributed to the paper.

\section{Acknowledgments}

The authors would like to thank Alexander Zheglov for useful discussions. 

\section{Funding}
The work of S. Abenda was supported by HORIZON-MSCA-2022-SE-01-01 CaLIGOLA,  COST Action CaLISTA CA21109, GNFM-INdAM,  INFN projects MMNLP and GAST.
The work of P.G. Grinevich was performed at the Steklov International Mathematical Center and supported by the Ministry of Science and Higher Education of the Russian Federation (agreement no. 075-15-2025-303).





\end{document}